\documentclass[11pt]{article}
\usepackage{amsmath, upgreek}
\usepackage{amssymb}
\usepackage{amsfonts}
\usepackage{amsmath,tabu}
\usepackage{cite}
\usepackage{color}

\usepackage{xcolor}
\usepackage{amscd}
\usepackage{xspace}
\usepackage{verbatim}
\usepackage{graphicx}
\newcommand{\mysection}[1]{
\section{#1}\setcounter{equation}{0}}
\title{\bf Local behaviour of the solutions of the Chipot-Weissler equation
}
\author{{\bf Marie-Fran\c{c}oise Bidaut-V\'eron\footnote{\noindent Laboratoire de Math\'{e}matiques et Physique Th\'{e}orique,
UMR 7013, Universit\'e de Tours, 37200 Tours, France. E-mail: veronmf@univ-tours.fr}} \\
 {\bf Laurent V\'eron \footnote{\noindent
Laboratoire de Math\'{e}matiques et Physique Th\'{e}orique, UMR 7013,  Universit\'e de Tours, 37200 Tours, France. E-mail: veronl@univ-tours.fr}}\\[2mm]
}

\date{}
\begin{document}
 \maketitle


\newcommand{\txt}[1]{\;\text{ #1 }\;}
\newcommand{\tbf}{\textbf}
\newcommand{\tit}{\textit}
\newcommand{\tsc}{\textsc}
\newcommand{\trm}{\textrm}
\newcommand{\mbf}{\mathbf}
\newcommand{\mrm}{\mathrm}
\newcommand{\bsym}{\boldsymbol}
\newcommand{\scs}{\scriptstyle}
\newcommand{\sss}{\scriptscriptstyle}
\newcommand{\txts}{\textstyle}
\newcommand{\dsps}{\displaystyle}
\newcommand{\fnz}{\footnotesize}
\newcommand{\scz}{\scriptsize}
\newcommand{\be}{\begin{equation}}
\newcommand{\bel}[1]{\begin{equation}\label{#1}}
\newcommand{\ee}{\end{equation}}
\newcommand{\eqnl}[2]{\begin{equation}\label{#1}{#2}\end{equation}}
\newcommand{\barr}{\begin{eqnarray}}
\newcommand{\earr}{\end{eqnarray}}
\newcommand{\bars}{\begin{eqnarray*}}
\newcommand{\ears}{\end{eqnarray*}}
\newcommand{\nnu}{\nonumber \\}
\newtheorem{subn}{\name}
\renewcommand{\thesubn}{}
\newcommand{\bsn}[1]{\def\name{#1}\begin{subn}}
\newcommand{\esn}{\end{subn}}
\newtheorem{sub}{\name}[section]
\newcommand{\dn}[1]{\def\name{#1}}   
\newcommand{\bs}{\begin{sub}}
\newcommand{\es}{\end{sub}}
\newcommand{\bsl}[1]{\begin{sub}\label{#1}}
\newcommand{\bth}[1]{\def\name{Theorem}
\begin{sub}\label{t:#1}}
\newcommand{\blemma}[1]{\def\name{Lemma}
\begin{sub}\label{l:#1}}
\newcommand{\bcor}[1]{\def\name{Corollary}
\begin{sub}\label{c:#1}}
\newcommand{\bdef}[1]{\def\name{Definition}
\begin{sub}\label{d:#1}}
\newcommand{\bprop}[1]{\def\name{Proposition}
\begin{sub}\label{p:#1}}

\newcommand{\R}{\eqref}
\newcommand{\rth}[1]{Theorem~\ref{t:#1}}
\newcommand{\rlemma}[1]{Lemma~\ref{l:#1}}
\newcommand{\rcor}[1]{Corollary~\ref{c:#1}}
\newcommand{\rdef}[1]{Definition~\ref{d:#1}}
\newcommand{\rprop}[1]{Proposition~\ref{p:#1}}
\newcommand{\BA}{\begin{array}}
\newcommand{\EA}{\end{array}}
\newcommand{\BAN}{\renewcommand{\arraystretch}{1.2}
\setlength{\arraycolsep}{2pt}\begin{array}}
\newcommand{\BAV}[2]{\renewcommand{\arraystretch}{#1}
\setlength{\arraycolsep}{#2}\begin{array}}
\newcommand{\BSA}{\begin{subarray}}
\newcommand{\ESA}{\end{subarray}}
\newcommand{\BAL}{\begin{aligned}}
\newcommand{\EAL}{\end{aligned}}
\newcommand{\BALG}{\begin{alignat}}
\newcommand{\EALG}{\end{alignat}}
\newcommand{\BALGN}{\begin{alignat*}}
\newcommand{\EALGN}{\end{alignat*}}
\newcommand{\note}[1]{\textit{#1.}\hspace{2mm}}
\newcommand{\Proof}{\note{Proof}}
\newcommand{\qeda}{\hspace{10mm}\hfill $\square$}
\newcommand{\qed}{\\
${}$ \hfill $\square$}
\newcommand{\Remark}{\note{Remark}}
\newcommand{\modin}{$\,$\\[-4mm] \indent}
\newcommand{\forevery}{\quad \forall}
\newcommand{\set}[1]{\{#1\}}
\newcommand{\setdef}[2]{\{\,#1:\,#2\,\}}
\newcommand{\setm}[2]{\{\,#1\mid #2\,\}}
\newcommand{\mt}{\mapsto}
\newcommand{\lra}{\longrightarrow}
\newcommand{\lla}{\longleftarrow}
\newcommand{\llra}{\longleftrightarrow}
\newcommand{\Lra}{\Longrightarrow}
\newcommand{\Lla}{\Longleftarrow}
\newcommand{\Llra}{\Longleftrightarrow}
\newcommand{\warrow}{\rightharpoonup}
\newcommand{
\paran}[1]{\left (#1 \right )}
\newcommand{\sqbr}[1]{\left [#1 \right ]}
\newcommand{\curlybr}[1]{\left \{#1 \right \}}
\newcommand{\abs}[1]{\left |#1\right |}
\newcommand{\norm}[1]{\left \|#1\right \|}
\newcommand{
\paranb}[1]{\big (#1 \big )}
\newcommand{\lsqbrb}[1]{\big [#1 \big ]}
\newcommand{\lcurlybrb}[1]{\big \{#1 \big \}}
\newcommand{\absb}[1]{\big |#1\big |}
\newcommand{\normb}[1]{\big \|#1\big \|}
\newcommand{
\paranB}[1]{\Big (#1 \Big )}
\newcommand{\absB}[1]{\Big |#1\Big |}
\newcommand{\normB}[1]{\Big \|#1\Big \|}
\newcommand{\produal}[1]{\langle #1 \rangle}

\newcommand{\thkl}{\rule[-.5mm]{.3mm}{3mm}}
\newcommand{\thknorm}[1]{\thkl #1 \thkl\,}
\newcommand{\trinorm}[1]{|\!|\!| #1 |\!|\!|\,}
\newcommand{\bang}[1]{\langle #1 \rangle}
\def\angb<#1>{\langle #1 \rangle}
\newcommand{\vstrut}[1]{\rule{0mm}{#1}}
\newcommand{\rec}[1]{\frac{1}{#1}}
\newcommand{\opname}[1]{\mbox{\rm #1}\,}
\newcommand{\supp}{\opname{supp}}
\newcommand{\dist}{\opname{dist}}
\newcommand{\myfrac}[2]{{\displaystyle \frac{#1}{#2} }}
\newcommand{\myint}[2]{{\displaystyle \int_{#1}^{#2}}}
\newcommand{\mysum}[2]{{\displaystyle \sum_{#1}^{#2}}}
\newcommand {\dint}{{\displaystyle \myint\!\!\myint}}
\newcommand{\q}{\quad}
\newcommand{\qq}{\qquad}
\newcommand{\hsp}[1]{\hspace{#1mm}}
\newcommand{\vsp}[1]{\vspace{#1mm}}
\newcommand{\ity}{\infty}
\newcommand{\prt}{\partial}
\newcommand{\sms}{\setminus}
\newcommand{\ems}{\emptyset}
\newcommand{\ti}{\times}
\newcommand{\pr}{^\prime}
\newcommand{\ppr}{^{\prime\prime}}
\newcommand{\tl}{\tilde}
\newcommand{\sbs}{\subset}
\newcommand{\sbeq}{\subseteq}
\newcommand{\nind}{\noindent}
\newcommand{\ind}{\indent}
\newcommand{\ovl}{\overline}
\newcommand{\unl}{\underline}
\newcommand{\nin}{\not\in}
\newcommand{\pfrac}[2]{\genfrac{(}{)}{}{}{#1}{#2}}

\def\ga{\alpha}     \def\gb{\beta}       \def\gg{\gamma}
\def\gc{\chi}       \def\gd{\delta}      \def\ge{\epsilon}
\def\gth{\theta}                         \def\vge{\varepsilon}
\def\gf{\phi}       \def\vgf{\varphi}    \def\gh{\eta}
\def\gi{\iota}      \def\gk{\kappa}      \def\gl{\lambda}
\def\gm{\mu}        \def\gn{\nu}         \def\gp{\pi}
\def\vgp{\varpi}    \def\gr{\rho}        \def\vgr{\varrho}
\def\gs{\sigma}     \def\vgs{\varsigma}  \def\gt{\tau}
\def\gu{\upsilon}   \def\gv{\vartheta}   \def\gw{\omega}
\def\gx{\xi}        \def\gy{\psi}        \def\gz{\zeta}
\def\Gg{\Gamma}     \def\Gd{\Delta}      \def\Gf{\Phi}
\def\Gth{\Theta}
\def\Gl{\Lambda}    \def\Gs{\Sigma}      \def\Gp{\Pi}
\def\Gw{\Omega}     \def\Gx{\Xi}         \def\Gy{\Psi}

\def\CS{{\mathcal S}}   \def\CM{{\mathcal M}}   \def\CN{{\mathcal N}}
\def\CR{{\mathcal R}}   \def\CO{{\mathcal O}}   \def\CP{{\mathcal P}}
\def\CA{{\mathcal A}}   \def\CB{{\mathcal B}}   \def\CC{{\mathcal C}}
\def\CD{{\mathcal D}}   \def\CE{{\mathcal E}}   \def\CF{{\mathcal F}}
\def\CG{{\mathcal G}}   \def\CH{{\mathcal H}}   \def\CI{{\mathcal I}}
\def\CJ{{\mathcal J}}   \def\CK{{\mathcal K}}   \def\CL{{\mathcal L}}
\def\CT{{\mathcal T}}   \def\CU{{\mathcal U}}   \def\CV{{\mathcal V}}
\def\CZ{{\mathcal Z}}   \def\CX{{\mathcal X}}   \def\CY{{\mathcal Y}}
\def\CW{{\mathcal W}} \def\CQ{{\mathcal Q}}
\def\BBA {\mathbb A}   \def\BBb {\mathbb B}    \def\BBC {\mathbb C}
\def\BBD {\mathbb D}   \def\BBE {\mathbb E}    \def\BBF {\mathbb F}
\def\BBG {\mathbb G}   \def\BBH {\mathbb H}    \def\BBI {\mathbb I}
\def\BBJ {\mathbb J}   \def\BBK {\mathbb K}    \def\BBL {\mathbb L}
\def\BBM {\mathbb M}   \def\BBN {\mathbb N}    \def\BBO {\mathbb O}
\def\BBP {\mathbb P}   \def\BBR {\mathbb R}    \def\BBS {\mathbb S}
\def\BBT {\mathbb T}   \def\BBU {\mathbb U}    \def\BBV {\mathbb V}
\def\BBW {\mathbb W}   \def\BBX {\mathbb X}    \def\BBY {\mathbb Y}
\def\BBZ {\mathbb Z}

\def\GTA {\mathfrak A}   \def\GTB {\mathfrak B}    \def\GTC {\mathfrak C}
\def\GTD {\mathfrak D}   \def\GTE {\mathfrak E}    \def\GTF {\mathfrak F}
\def\GTG {\mathfrak G}   \def\GTH {\mathfrak H}    \def\GTI {\mathfrak I}
\def\GTJ {\mathfrak J}   \def\GTK {\mathfrak K}    \def\GTL {\mathfrak L}
\def\GTM {\mathfrak M}   \def\GTN {\mathfrak N}    \def\GTO {\mathfrak O}
\def\GTP {\mathfrak P}   \def\GTR {\mathfrak R}    \def\GTS {\mathfrak S}
\def\GTT {\mathfrak T}   \def\GTU {\mathfrak U}    \def\GTV {\mathfrak V}
\def\GTW {\mathfrak W}   \def\GTX {\mathfrak X}    \def\GTY {\mathfrak Y}
\def\GTZ {\mathfrak Z}   \def\GTQ {\mathfrak Q}

\font\Sym= msam10 
\def\SYM#1{\hbox{\Sym #1}}
\newcommand{\bdw}{\prt\Gw\xspace}
\date{}
\maketitle\medskip

\noindent{\small {\bf Abstract} We study the local properties of positive solutions of the equation $-\Gd u=u^p-m\abs{\nabla u}^q$ in a punctured domain $\Gw\setminus\{0\}$ of $\BBR^N$ or in a exterior domain $\BBR^N\setminus B_{r_0}$ in the range $\min\{p,q\}>1$ and $m> 0$. We prove a series of a priori estimates depending $p$ and $q$, and of the sign of $q-\frac {2p}{p+1}$ and $q-p$. Using various techniques we obtain removability results for singular sets and we give a precise description of behaviour of solutions near an isolated singularity or at infinity in $\BBR^N$.
}\medskip

\tableofcontents
\vspace{1mm}
\hspace{.05in}
\medskip

\mysection{Introduction}
The aim of this paper is to study the local properties of solutions of 
\bel{I-1}
\CL_{m,p,q}u:=-\Gd u +m|\nabla u|^q-|u|^{p-1}u=0\qquad\text{in }\Gw,
\ee
where $m$ is a nonnegative real number, $p,q\geq 1$ and $\Gw$ is either a punctured domain if we are interested in isolated singularities, or an exterior domain if we study the asymptotic behaviour of solutions. This equation  has been introduced by Chipot and Weissler \cite{ChWe1} in connection with the associated evolution problem
\bel{I-1+1}
\prt_tu+\CL_{m,p,q}u=0\qquad\text{in }\Gw\ti (0,T).
\ee
Its study has been developed in the radial case in \cite{ChWe2} and completed in \cite{Vo}. A very deep research of radial ground states has been carried on by Serrin and Zou in \cite{SeZo} and \cite{SeZo2}. Several non-existence results of positive, not necessarily radial, supersolutions in an exterior domain have been obtained in \cite{AGMQ1} and  \cite{ABGMQ}.\smallskip

The interest of the operator $\CL_{m,p,q}$ lies in the presence of two reaction terms which are acting in opposite directions and are of a different nature. The following exponents play a key role in the study of asymptotics of solutions of $(\ref{I-1})$,
\bel{I-2}
\ga=\frac{2}{p-1}\,,\;\gb=\frac{2-q}{q-1}\,,\; \gg=\frac{q}{p-q}\,\text{if }q\neq p\;\text{ and }\;\gs=(p+1)q-2p.
\ee
When $q=\frac{2p}{p+1}$ the equation $(\ref{I-1})$ is invariant under the transformation $T_\ell$ defined by 
\bel{I-2*}
T_\ell[u](x)=\ell^\ga u(\ell x).
\ee
This critical value of $q$ plays a fundamental role in the analysis of the solutions. If $1<q<\frac{2p}{p+1}$, the source term is dominant for large values of $u$ e.g. near a singular point, and the behaviour of singular solutions is modelled by the 
{\it Lane-Emden equation}
 \bel{I-3}  
-\Gd u-u^p=0.
 \ee
 If $\frac{2p}{p+1}<q<p$, the diffusion is negligeable and the behaviour of singular solutions is modelled by an {\it eikonal  equation}
 \bel{I-4}  
u^p-m|\nabla u|^q=0.
 \ee
 Notice that in this equation the sign of $p-q$ is fundamental and makes the distinction between the existence or the non-existence of singular solutions. 
 Another equation which plays a crucial  role is the {\it Riccatti equation}
 \bel{I-5}  
-\Gd u+m|\nabla u|^q=0.
 \ee
For this equation the value of $q$ with respect to $2$ is the key element. 
Finally, if $q=\frac{2p}{p+1}$ no reaction term is dominant and the value of $m$ becomes fundamental as the following result   proved in \cite{BiGaVe2} shows it:\medskip
  
  \nind {\bf Theorem A} {\it Let $N\geq 2$, $1<p<\frac{N+2}{N-2}$ and $q=\frac{2p}{p+1}$. Then there exist two positive constants $c=c_(N,p)$ and $m_0$ such that for any real number $m$ verifying  $|m|\leq m_0$, any positive solution $u$ of $(\ref{I-1})$ in $\Gw$ satisfies
   \bel{I-6}  
u(x)+|\nabla u(x)|^{\frac{2}{p+1}}\leq c\left(\dist(x,\prt\Gw)\right)^{-\ga}\qquad\text{for all }x\in\Gw.
 \ee
 As a consequence there exists no positive solution (called ground state) in $\BBR^N$.
}\medskip

An a priori estimate holds by a perturbation method for positive solutions, for all values of $m$ whenever $1<p<\frac{N+2}{N-2}$, and the following result is obtained in  \cite{PoQuSo}.\medskip

\nind{\bf Theorem B} {\it Let $N\geq 2$, $1<p<\frac{N+2}{N-2}$ and $1<q<\frac{2p}{p+1}$. For any $m\in\BBR$ there exists a positive constant $c=c(N,p,q,m)$ such that any positive solution $u$ of $(\ref{I-1})$ in $\Gw$ satisfies
   \bel{I-6*}  
u(x)+|\nabla u(x)|^{\frac{2}{p+1}}\leq c\left(1+\left(\dist(x,\prt\Gw)\right)^{-\ga}\right)\qquad\text{for all }x\in\Gw.
 \ee
}

Up to now, these two results were the only ones known concerning a priori estimates for general nonnegative solutions when $m>0$. In the present article we prove new upper estimates for positive solutions $u$ of $(\ref{I-1})$ either in a punctured domain $B_{r_0}\setminus\{0\}$ or in an exterior domain $\Gw=B_{r_0}^c$.\smallskip

The next statements extend previous results concerning positive supersolutions  proved in \cite{AGMQ1}. If $u$ is a positive continuous function defined either in $B_{r_0}\setminus\{0\}$ or in $B_{r_0}^c$,  we set
    \bel{I-7}  \dsps
\gm(r)=\inf_{|x|=r}u(x),
 \ee
and we prove the following estimates valid in the case $1<q<p$.\medskip

\bth{Theorem C} Let $N\geq 1$,  $p,q>1$ and $m>0$.\smallskip
 
 \nind 1- Let $u$ be a $C^2$ positive supersolution of $(\ref{I-1})$ in $B^c_{r_0}$, then 
 
\nind 1-(i) If $\frac{2p}{p+1}< q<p$ there exists $C=C(N,p,q,u)>0$ such that 
   \bel{as5}\dsps
\gm(r)\leq Cr^{-\ga}\quad\text{for all }r\geq 2r_0.
\ee

\nind 1-(ii) If $1<q\leq \frac{2p}{p+1}$ there exists $C=C(N,p,q,u)>0$ such that 
   \bel{as6}\dsps
\gm(r)\leq Cr^{-\gg}\quad\text{for all }\; r\geq 2r_0.
\ee

\nind 1-(iii) If $1<p\leq q$ and $\gm(|x|)$ is bounded, then $(\ref{as6})$ is still satisfied.\smallskip

\nind 2- Let $u$ be a positive supersolution of $(\ref{I-1})$ in $B_{r_0}\setminus\{0\}$, then\smallskip

\nind 2-(i) If $\frac{2p}{p+1}\leq q<p$ there exists $C=C(N,p,q,u)>0$ such that 
  \bel{as7}\dsps
\gm(r)\leq Cr^{-\gg}\quad\text{for all }0<r\leq\tfrac {r_0}{2}.
\ee

\nind 2-(ii) If $1<q< \frac{2p}{p+1}$ there exists $C=C(N,p,q,u)>0$ such that 

   \bel{as8}\dsps
\gm(r)\leq Cr^{-\ga}\quad\text{for all }0<r\leq\tfrac {r_0}{2}.
\ee
\es
All the estimates on $\gm(r)$ will play a crucial role for the study of radial solutions of $(\ref{I-1})$ see \cite{BiVe4}.\smallskip

In the case $q\geq p$, the upper estimates are no more satisfied. The next result points out a dichotomy for estimates of positive supersolutions in an exterior domain when $q\geq p$. \medskip

\bth{Theorem D} Let $N\geq 2$ and $1<p\leq q$. If $u$ is any  positive supersolution of $(\ref{I-1})$ in $B_{r_0}^c$, then 
  for any $\gr>{r_0}$ there exists $c_\gr $, $C_\gr$, $C'_{\gr}$, $C''_{\gr}>0$ such that, for $|x|\geq \gr$, \\
  (i) either 
   \bel{I-9}  
u(x)\geq \left\{\BA{lll} X_m|x|^{\frac{q}{q-p}}\left(1-\frac{C_\gr}{|x|}\right)_+&\qquad\text{if }q>p\\[2mm]
c_\gr e^{m^{-\frac 1m}|x|}&\qquad\text{if }q=p,
\EA\right.
 \ee
where $X_m=\left(m|\gg|^q\right)^{\frac{1}{p-q}}$,\smallskip

\nind (ii) or $p>\frac N{N-2}$ and 
   \bel{I-10}  \BA{lll}
(a) \qquad\qquad\qquad\qquad\qquad\qquad&\gm(|x|)\leq C'_{\gr}[x]^{-\ga}\qquad\quad\quad\quad\qquad\qquad\qquad\qquad\qquad\qquad\quad\\[2mm]
(b) \qquad\qquad &u(x)\geq C''_{\gr}[x]^{2-N}.
 \EA\ee
\es

When $q>p$, the function $U(x)=X_m|x|^{|\gg|}$ is a $C^1$ subsolution of $(\ref{I-1})$ in $\BBR^N$, a fact which shows the optimality of the lower estimate.  \smallskip

In the case $q>p$ we prove a series of new estimates of {\it solutions}, by a delicate combination of Bernstein,  Keller-Osserman methods and Moser iterative scheme. The general Bernstein estimates will play a fundamental role in the description of the behaviour of positive solutions near an isolated singularity or at infinity in $\BBR^N$.
\medskip

 \bth {Theorem E} Let $q>p>1$, $m>0$ and $u$ be a nonnegative solution of $(\ref{I-1})$ in a domain $G\subset\BBR^N$. Then\smallskip

\nind 1- If $G= B_{r_0}\setminus\{0\}$, there exists $c>0$ depending on $N,p,q$ and $\norm u_{L^\infty(B_{r_0}\setminus B_{\frac{3r_0}{4}})}$ such that 
\bel{Est-N1}
\abs {\nabla u(x)}\leq c|x|^{-\frac 1{q-1}}\quad\text{for all }0<|x|\leq \frac {r_0}2.
\ee
\nind 2- If $G= B^c_{r_0}$, there exists $c>0$ depending on $N$, $p$, $q$ and $\norm u_{L^\infty(B_{2r_0}\setminus B_{r_0})}$  such that 
\bel{Est-N2}
\abs {\nabla u(x)}\leq c|x|^{\frac p{q-p}}\quad\text{for all }|x|\geq 2{r_0}.
\ee
\es

Note that in $B_{r_0}\setminus\{0\}$ the dominant effect comes from the Riccatti equation, while it comes from the eikonal equation in $B^c_{r_0}$. However it concerns solutions which may blow-up at infinity. When $q<p$, the {\it eikonal equation} plays a fundamental role in the proof of the next result which uses all the previous techniques involved in the proof of \rth{Theorem E} above combined with the doubling Lemma method of \cite{Hu}.

\bth {Theorem E'}  Let $p>1$, $m>0$ and $r_0>0$. \smallskip

\nind 1- Let $1<q<\frac{2p}{p+1}$. If $u$ is a positive solution of $(\ref{I-1})$ in $B^c_{r_0}$ satisfying 
\bel{b3}
\dsps \lim_{|x|\to\infty}u(x)=0, 
\ee
then there exists a positive constant $C=C(N,p,q,u,r_0,m)$ such that 
\bel{b4}
u(x)\leq C|x|^{-\frac{q}{p-q}}\,\text{and  }\; |\nabla u(x)|\leq C|x|^{-\frac{p}{p-q}}
\ee
 for all $x\in B^c_{2r_0}$.\smallskip

\nind 2- Let $\frac{2p}{p+1}<q<p$. Any $u$ positive solution $u$ of $(\ref{I-1})$ in $B_{r_0}\setminus\{0\}$  satisfies $(\ref{b4})$
for all $x\in B_{\frac{r_0}{2}}\setminus\{0\}$ for some constant $C=C(N,p,q,u,r_0,m)>0$.
\es

In a forthcoming article \cite {BiVe4} we prove the existence of infinitely many  different radial solutions satisfying the decay estimate $(\ref{b4})$ by a combination of  ODE and dynamical systems approach.\smallskip
 
 The following result is the counterpart at infinity  Theorems A and 
  B.\medskip

  \bth {Theorem F} Let $1<p<\frac{N+2}{N-2}$, $m>0$ and $u$ be a positive solution of $(\ref{I-1})$ in $B^c_{{r_0}}$ ($r_0>0$) satisfying 
\bel{AS1}
\dsps \lim_{|x|\to\infty}u(x)=0. 
\ee
Assume \smallskip

\nind (i) either $\frac{2p}{p+1}<q\leq 2$ and $m$ is arbitrary, \\
\nind (ii) or $q=\frac{2p}{p+1}$ and $m\leq \ge_0$ for some $\ge_0>0$ depending on $N$ and $p$.\smallskip

\nind Then there exists a positive constant $C=C(N,p,q,u,r_0,m)$ such that 
\bel{AS2}
u(x)\leq C|x|^{-\frac{2}{p-1}}\,\text{and  }\; |\nabla u(x)|\leq C|x|^{-\frac{p+1}{p-1}}\quad\text{for all }x\in B^c_{2r_0}.
\ee
\es

Thanks to the estimates of \rth {Theorem E} we can prove removability results for singularities of positive solutions of $(\ref{I-1})$.\medskip

  \bth{Theorem G} Let $N\geq 2$, $\Gw\subset\BBR^N$  be a bounded smooth domain containing $0$. If $1\leq p< q$ and $q\geq \frac{N}{N-1}$, any  nonnegative  solution $u\in C^2(\Gw\setminus\{0\})$ of  $(\ref{I-1})$ in $\Gw\setminus\{0\}$ can be extended as a weak solution of the same equation in $\Gw$ and it belongs to $L^\infty_{loc}(\Gw)\cap W^{1,q}_{loc}(\Gw)\cap H^1_{loc}(\Gw)$.\es
  
  This result admits extensions for removability of more general sets included in a domain $\Gw\subset\BBR^N$ in two completely different directions. Using a geometric construction as in \cite{Ve1} we prove: \medskip
  
  \bth {Theorem H} Let $N\geq 3$, $\Gw\subset\BBR^N$  be a bounded domain, $\Gs\subset\Gw$ a $k$-dimensional compact complete submanifold ($0\leq k\leq N-2$), $m>0$ and $1\leq p<q$ such that $q\geq \frac{codim (\Gs)}{codim (\Gs)-1}$. Then any positive solution of $(\ref{I-1})$ in $\Gw\setminus \Gs$ is locally bounded and can be extended as a weak solution in $\Gw$.\es
  
  Using capacitary estimates we extend to the case $q> 2$ a previous removability result due to Brezis and Nirenberg \cite{BrNi} obtained in the case $q=2$.
\medskip

  \bth {Theorem I} Assume $p>0$, $q\geq\max\{2,p\}$ and  $m>0$. If $K$ is a compact subset of $\Gw$ such that $cap_{1,q'}(K)=0$, then any positive solution of $(\ref{I-1})$ in $\Gw\setminus K$ is locally bounded and can be extended as a weak solution in $\Gw$.\es
 
  The last Section is devoted to the study of asymptotics of positive solutions, either near a singularity or at infinity. In the case $q<\frac{2p}{p+1}$ the dominant equation for the study of isolated singularity is the Lane-Emden one, and the techniques involved combine energy methods and Fourier analysis. The description of the singular behaviour depends upon the value of $p$ with respect to 
  $\frac {N}{N-2}$ and  $\frac {N+2}{N-2}$, and we obtain the complete classification of the possible behaviours of a positive solution near an isolated singularity:
  
  \bth {Theorem J} Let $N\geq 2$, $m>0$, $1<p<\frac{N+2}{N-2}$ and $1<q<\frac{2p}{p+1}$. If $u$ is a nonnegative solution of  $(\ref{I-1})$ in $B_{r_0}\setminus\{0\}$, then either $u$ is a classical solution of $(\ref{I-1})$ in $B_{r_0}$, or\smallskip

\nind 1- when  $N\geq 3$ and $1<p<\frac{N}{N-2}$ (resp. $N=2$ and $p>1$) there exists $k> 0$ such that $|x|^{N-2}u(x)$   (resp. $-u(x)/\ln |x|$) converges to $k$  when $x\to 0$. Furthermore $u$ satisfies
\bel{SI1}
-\Gd u+m|\nabla u|^q-u^p=c_Nk\gd_0\quad\text{in }\CD'(B_{r_0});
\ee

\nind 2- when $N\geq 3$ and $p=\frac{N}{N-2}$, $|x|^{N-2}(-\ln |x|)^{\frac{N-2}{2}} u(x)$ converges to $\left(\frac{N-2}{\sqrt 2}\right)^{N-2}$ when $x\to 0$;\smallskip

\nind 3- when $N\geq 3$ and $\frac{N}{N-2}<p<\frac{N+2}{N-2}$, $|x|^\ga u(x)$ converges to $\gw_0:=\left(\ga\frac{(N-2)p-N}{p-1}\right)^{\frac{1}{p-1}}$ when $x\to 0$.
\es
 
 In the case $q>p$ the dominant equation near an isolated singularity is the Riccatti equation; the removability result of \rth{Theorem G} is no more valid if $1<q<\frac{N}{N-1}$, and we mainly use a scaling method. 
  
     \bth {Theorem K} Let $N\geq 3$, $1<p<q<\frac{N}{N-1}$, $m>0$ and $u$ be a nonnegative solution of  $(\ref{I-1})$ in $B_{r_0}\setminus\{0\}$. Then either $u$ is a classical solution,  \\
(i) or $|x|^\gb u(x)$ converges to $\xi_m:=\frac 1\gb\left(\frac{(N-1)q-N}{m(q-1)}\right)^{\frac 1{q-1}}$  when $x\to 0$,\smallskip

\nind(ii) or there exists $k>0$ such that $|x|^{N-2}u(|x|,.)\to c_Nk$ in $L^1(S^{N-1})$ when $x\to 0$ and $u$ satisfies 
$$-\Gd u+m|\nabla u|^q-u^p=k\gd_0\qquad\text{in }\CD'(B_{{r_0}}).$$
\es

The asymptotic behaviour of solutions in an exterior domain exhibits also the two types of underlying dominant equations: either the Lane-Emden equation, or the eikonal equation. This depends on the value of $q$ with respect to $\frac{2p}{p+1}$, see \rth{AST1}, \rth{AST2}. The techniques are similar to the ones used in the analysis of isolated singularities but the range of values of $q$ are reversed; a phenomenon which is easily understandable when considering the scaling transformations leaving the underlying equations invariant. \medskip

\mysection{Estimates on supersolutions}
\subsection{Some preliminary results}
In the sequel we denote by $c$ or $C$ a generic positive constant the value of which may vary from one occurence to another. When needed we introduce 
the constants $c_i$, $C_i$ with $i=1,2,...$, in particular within the development  of the proof of a statement. If it is important we precise the parameters ($N$, $p$, $q$, $m$ etc.) on which the various constants depend. In the next result we precise a bootstrap argument some variants of which have already been used in \cite{BiGr}, \cite {BiGaYa} and \cite{Bi2}.

\blemma{ST0} Let $d,\,h \in\BBR$ with $0<d<1$ and $y$, $\Phi$ be two positive continuous functions defined on $(0,{r_0}]$ (resp. $[{r_0},\infty)$). We assume that there exist $C^*,M>0$ and $\ge_0\in (0,\tfrac 18]$ such that for any $\ge\in (0,\ge_0]$ and $0<r\leq\tfrac {r_0}2$ (resp. any $r\geq 2{r_0}$),
\bel{as0}\dsps
y(r)\leq C^*\ge^{-h}\Phi(r)y^d(r(1-\ge))\quad\text{and }\max_{\tfrac{
r}{2}\leq \gt\leq r}\Phi(\gt)\leq M\Phi(r),
\ee
respectively
\bel{as1}\dsps
 y(r)\leq C^*\ge^{-h}\Phi(r)y^d(r(1+\ge))\quad\text{and }\max_{ r\leq \gt\leq\tfrac{3r}{2}}\Phi(\gt)\leq M\Phi(r).
\ee
Then there exists $c_1=c_1(C^*,M,d,h,\ge_0)>0$ such that 
\bel{as2}\dsps
y(r)\leq c_1\left(\Phi(r)\right)^{\tfrac 1{1-d}},\ee
in $(0,\tfrac {{r_0}}2]$ (resp. in $[2{r_0},\infty)$).
\es
\Proof The result is obvious when $h\leq 0$, so we can suppose $h>0$. Consider the sequence 
$\ge_n=2^{-n}\ge_0$, $n\geq 0$. Then the series $\sum\ge_n$ is convergent and 
$$S=\sum_{j=1}^\infty\ge_j\leq\frac{1}{4}.
$$ 
For $n\geq 1$ we denote $P_n=(1-\ge_1)...(1-\ge_j)...(1-\ge_n)$ and $Q_n=(1+\ge_1)...(1+\ge_j)...(1+\ge_n)$. Clearly the sequence $\{P_n\}$ is decreasing while  the sequence $\{Q_n\}$ is increasing. Furthermore
$$Q_n\leq\prod_{j=1}^\infty(1+\ge_j):=Q\leq e^{S}\leq e^{\frac 14}<\frac 32.
$$
Concerning $P_n$, we have $1-\ge_n>\frac{1}{1+2\ge_n}$. Therefore
$$P_n\geq \prod_{j=1}^n(1+2\ge_j)^{-1}\geq e^{-2S}\geq e^{-\frac 12},
$$
which implies $\frac 12<P_n<1$. Then , for any $r\in (0,\frac {r_0}2]$ (resp. $r\geq 2{r_0}$) we have that $rP_n\in [\frac r2,r]$ 
(resp. $rQ_n\in [r,\frac{3r}{2}]$). First we assume $(\ref{as0})$ and use $P_n$. Then 
$$y(rP_{n-1})\leq c_2\ge_n^{-h}\Phi(rP_{n-1})y^d(rP_n).
$$
In particular
$$\left\{\BA{lll}
y(r)\leq c_2\ge_1^{-h}\Phi(r)y^d(rP_1))\\[1mm]
y^d(rP_1)\leq c_2^d\ge_2^{-hd}\Phi^d(rP_1)y^{d^2}(rP_2))\\
[-2pt] \vdots  \\
y^{d^{n-1}}(rP_{n-1})\leq c_2^{d^{n-1}}\ge_n^{-hd^{n-1}}\Phi^{d^{n-1}}(rP_{n-1})y^{^{d^{n}}}(rP_n)).
\EA\right.$$
By the assumption on $\Phi$, this implies
$$y(r)\leq c_2^{1+d+d^2+\dot+d^{n-1}}\ge_1^{-h}\ge_2^{-hd}\dots\ge_n^{-hd^{n-1}}
\Phi(r)\phi^d(rP_1)\dots\Phi^{d^{n-1}}(rP_{n-1})y^{d^n}(rP_n),
$$
for any $n\geq 2$. Hence for any $n\geq 2$, 
\bel{as3-}\BA{lll}\dsps y(r)\leq (c_2\ge_0^{-h})^{1+d+\dots+d^{n-1}}2^{h(1+2d+...+nd^{n-1})}
\Phi(r)\Phi^{d}(rP_1)...\Phi^{d^{n-1}}(rP_{n-1})y^{d^n}(rP_n)\\[2mm]
\phantom{y(r)}\dsps\leq
(c_2\ge_0^{-h})^{1+d+\dots+d^{n-1}}2^{h(1+2d+...+nd^{n-1})}M^{d+d^2+...d^{n-1}}\Phi^{1+d+d^2+...d^{n-1}}(r).
\EA\ee
Letting $n\to\infty$ and using the fact that $P_n\to P>0$ and $y^{d^n}(rP_n)\to 1$ as $n\to\infty$, since $0<d<1$, we obtain  
\bel{as4-1}
y(r)\leq (c_2\ge_0^{-h})^{\frac {1}{1-d}}2^{\frac{h}{(1-d)^2}}M^{\frac d{1-d}}
\left(\Phi(r)\right)^{\frac{1}{1-d}}.
\ee
If we assume $(\ref{as1})$, the proof of $(\ref{as2})$ in $[2{r_0},\infty)$ is similar.\qeda
\smallskip

Next we recall and extend the monotony property dealing with supersolutions of Riccatti equation proved in \cite{AGMQ1}.

\blemma{ST2} Let $N\geq 2$, $q>1$ and $u\in C^2(B_{r_0}\setminus\{0\})$ (resp. $u\in C^2(B^c_{r_0})$) be a positive function such that 
$$-\Gd u+|\nabla u|^q\geq 0\quad\text{in } B_{r_0}\setminus\{0\}\quad \text{(resp. in }B^c_{r_0}).$$
Then the function $\gm$ defined by $(\ref{I-7})$ is nonincreasing on $(0,r_0]$ (resp. there exists $r_1\geq r_0$ such that $\gm$ is monotone 
on $[r_1,\infty)$).
\es
\Proof The case of an exterior domain is treated in \cite[Lemma 5]{AGMQ1}. In the first case, then for any $r_1\in (0,r_0)$ and $\gd>0$ there exists 
$r_d\in (0, r_1]$ such that for any $0<r\leq r_\gd$ such that $\gm(r_1)\leq \gd r^{2-N}$ if $N\geq 3$ or $\gm(r_1)\leq \gd |\ln r|$ if $N=2$. Let 
$h(x)=\gm(r_1)-\gd|x|^{2-N}$ if $N\geq 3$ (resp. $h(x)=\gm(r_1)-\gd| |\ln |x||$ if $N=2)$. Then $u\geq h$ on $\prt B_{r_1}\cup \prt B_{r}$. By the standard comparison principle \cite{AGMQ1},\cite{Ng1}, $u\geq h$ in $\overline B_{r_1}\setminus B_{r}$. If we let $r\to 0$ we derive 
$u\geq h$ in $\overline B_{r_1}\setminus\{0\}$, and by letting $\gd\to 0$ we finally obtain $u\geq \gm(r_1)$  in $\overline B_{r_1}\setminus\{0\}$. In particular this inequality implies $\gm(r)\geq \gm(r_1)$ if $0<r\leq r_1$.\qeda

\subsection{Estimates of the spherical minimum. Proof of \rth{Theorem C}}
In this Section we consider non-necessarily radial supersolutions $u$ of $(\ref{I-1})$, either in a punctured or in an exterior domain. We give estimates of the minimum of $u$ on spheres with center $0$ $\dsps\gm(r)=\min_{|y|=r}u(y)$.\smallskip

We first consider supersolutions of the exterior problem
 \bel{as3-}\dsps
-\Gd u+m|\nabla u|^q-f(u)= 0\quad\text{in }B^c_{r_0},
\ee
where $m>0$ and $f$ satisfies\smallskip

\nind (F)  {\it $\phantom{--}f$ is a continuous nondecreasing  function on $\BBR_+$ verifying $f(0)=0$ and $f>0$ on $(0,\infty)$}.\smallskip

We recall the following result of \cite[Theorems 1, 3, 4]{AGMQ1}.
\medskip

\nind {\bf Theorem C} {\it (1) If $\dsps\liminf_{r\to 0}r^{-p}f(r)>0$ and $1<p\leq \tfrac{N}{N-2}$, $q>\frac{2p}{p+1}$, there exists no positive supersolution $u\in C^2(B_{r_0}^c)$ of $(\ref{as3-})$ such that $\dsps\liminf_{|x|\to \infty}u(x)<\infty$.\smallskip

\nind (2) If $\dsps\liminf_{r\to \infty}r^{-p}f(r)>0$ and $1<q<p$, there exists no positive supersolution $u\in C^2(B_{r_0}^c)$ of $(\ref{as3-})$ such that $\dsps\lim_{|x|\to \infty}u(x)=\infty$.}\medskip

Here we combine a technique developed in \cite[Lemma 6]{AGMQ1} in order to prove \rth{Theorem C} with the bootstrap argument of \rlemma{ST0}.\smallskip

\blemma{LX} Let $m>0$, $N\geq 1$, $q>1$ and $f$ satisfying (F). Let $u\in C^2(B_{r_0}^c)$ (resp $u\in C^2(B_{r_0}\setminus\{0\})$) be any positive function
satisfying 
\bel{LX1}
-\Gd u+m|\nabla u|^q\geq f(u)\quad\text{in }B_{r_0}^c\quad\left(\text{resp. in } B_{r_0}\setminus\{0\})\right).  
\ee
\smallskip

\nind 1- Then for any $R\geq 2r_0$ (resp. for any $0<R\leq \frac{r_0}{2}$) and for any $0<\ge\leq\frac 12$,
\bel{LX2}
\dsps \min_{(1-\ge)R\leq r\leq (1+\ge)R } f(u(r))\leq c_1\left(\frac{\gm(R)}{\ge^2R^2}+\frac{\gm^q(R)}{\ge^qR^q}\right),
\ee
where $c_1=c_1(N,q,m)>0$.\smallskip

\nind 2- As a consequence, any positive $C^2$ supersolution $u$ of $(\ref{as3-})$ in $B_{r_0}^c$ satisfies \\
(i) either $\dsps\lim_{|x|\to\infty}u(x)=\infty$, \\
(ii) or $\dsps\liminf_{|x|\to\infty}u(x)=0$.
\es
\Proof 1- Let $R\geq 2r_0$ (resp. $0<R\leq \frac{r_0}2$) and $\ge\in (0,\frac 12]$.
Let $\phi_\ge$ be a smooth nonnegative radial cut-off function defined on $\BBR_+$, vanishing on $[0,1-\ge]\cup[1+\ge,\infty)$ with value $1$ on $[1-\tfrac\ge 2,1+\tfrac\ge 2]$, such that 
 $|\phi'_\ge|\leq\tfrac C\ge\chi_{_{I_\ge}}$ and $|\phi''_\ge|\leq\tfrac C{\ge^2}\chi_{_{I_\ge}}$ where $\chi_{_{I_\ge}}=[1-\ge,1-\tfrac\ge 2]\cup[1+\tfrac\ge 2,1+\ge]$. We set
 $$v(x)=u(x)-\gm(R)\phi_\ge (\tfrac{|x|}{R}).
 $$
 There exists $x_{R,\ge}$ such that $|x_{R,\ge}|=R$ and $u(x_{R,\ge})=\gm(R)$, thus $v(x_{R,\ge})=0$. If $u$ is defined in $B_{r_0}^c$, we have that 
 $v=u>0$ in $(B_{R(1-\ge)}\cap B_{r_0}^c)\cup B_{R(1+\ge)}^c$. If $u$ is defined in $B_{r_0}\setminus\{0\}$, then 
 $v=u>0$ in $(B_{R(1-\ge)}\setminus\{0\})\cup \left(B_{r_0}\cap B^c_{R(1+\ge)}\right)$. Then $v$ achieves its nonpositive minimum at some 
 $\widetilde x_{R,\ge}\in B_{R(1+\ge)}\cap \overline B_{R(1-\ge)}^c$, where $\nabla v(\widetilde x_{R,\ge})=0$ and $\Gd v(\widetilde x_{R,\ge})\geq 0$. Since $v(\widetilde x_{R,\ge})\leq 0$ there holds $\gm(|\widetilde x_{R,\ge}|)\leq\gm(R)$ and
 $$\BA {lll}f(u(\widetilde x_{R,\ge}))=-\Gd v(\widetilde x_{R,\ge})+m|\nabla v(\widetilde x_{R,\ge})|^q\\[2mm]
 \phantom{f(u(\widetilde x_{R,\ge}))}
 =-\gm(R)\Gd \left(\phi_\ge(\tfrac{|x|}{R})\right)+m\gm^q(R)\left|\nabla \left(\phi_\ge(\tfrac{|x|}{R})\right)\right|^q\\[2mm]
 \phantom{f(u(\widetilde x_{R,\ge}))}\dsps
 \leq c_1\left(\frac{\gm(R)}{\ge^2R^2}+\frac{\gm^q(R)}{\ge^qR^q}\right),
 \EA$$
 where $c_1=c_1(N,p,q,m)>0$. Because 
$\dsps u(\widetilde x_{R,\ge})\geq \min_{(1-\ge)R\leq r\leq (1+\ge)R} \gm(r)$, $(\ref{LX2})$ follows from the monotonicity of $f$. \smallskip

\nind 2- From \rlemma {ST2}, $\gm(r)$ is monotone for large $r$. \\
If $\gm$ is bounded, then 
$$\dsps \min_{\frac R2\leq r\leq 2R}f(\gm(r))\leq c_3\left(\frac{1}{R^2}+\frac{1}{R^q}\right).
$$
Hence $\dsps \lim_{R\to\infty}\min\left\{f(\gm(\tfrac R2)),f(\gm(2R))\right\}=0$ which implies that $\gm(R)\to 0$ when $R\to\infty$, since $f$ is continuous and vanishes only at $0$.\\
If $\gm$ is unbounded, then $\dsps \lim_{r\to\infty}\gm(r)=\infty$ which implies $\dsps \lim_{|x|\to\infty}u(x)=\infty$.\qeda 
 \medskip
 
 Now we assume that $f(u)=u^p$, $p>1$, and prove \rth{Theorem C}. We recall that the exponents $\ga$, $\gb$ and $\gg$ have been defined at $(\ref{I-2})$.\smallskip

 \nind{\it Proof of \rth{Theorem C}}. Let $p,q>1$ and $u$ be a positive supersolution of $(\ref{I-1})$ in $B_{r_0}^c$ (resp. $B_{r_0}\setminus\{0\}$). Let $R\geq 2r_0$ (resp. $0<R\leq\frac R2$). From \rlemma {LX}, we have that: \smallskip
 
 \nind if $\gm$ is nonincreasing on $[R-\ge,R+\ge]$, then $\gm(R)\geq u(\widetilde x_{R,\ge})\geq\gm(|\widetilde x_{R,\ge}|)\geq\gm(R(1+\ge))$,
then
\bel{hach}\BA{lll}\dsps\gm^p(R(1+\ge))\leq c_4\left(\frac{\gm(R)}{\ge^2R^2}+\frac{\gm^q(R)}{\ge^qR^q}\right)\leq c_4\ge^{-h}\left(\frac{\gm(R)}{R^2}+\frac{\gm^q(R)}{R^q}\right)\quad\text{with }h=\max\{2,q\},
\EA
\ee
  
  \nind if $\gm$ is nondecreasing on $[R-\ge,R+\ge]$, then  $\gm(R)\geq u(\widetilde x_{R,\ge})\geq \gm(|\widetilde x_{R,\ge}|)\geq\gm(R(1-\ge))$, then 
\bel{hoch}\BA{lll}\dsps\gm^p(R(1-\ge))\leq c_4\ge^{-h}\left(\frac{\gm(R)}{R^2}+\frac{\gm^q(R)}{R^q}\right).
\EA\ee
 Note that for any $c,R>0$ there holds
\bel{mas}
\frac{\gm^q(R)}{R^q}\leq c\frac{\gm(R)}{R^2}\Longleftrightarrow \gm(R)\leq c^{-\frac{1}{q-1}}R^{-\gb},
\ee
since $\gb=\frac{2-q}{q-1}$. 
\smallskip

\nind {\it 1- The exterior problem}. From \rlemma {ST2}, $\gm(r)$ is monotone for $R\geq r_1\geq r_0$ large enough, so we assume $R>r_1$, and either 
$\gm$ is decreasing or it increases to $\infty$. In our cases, we claim that $\gm$ is decreasing. It holds by assumption if $q\geq p$. When $q<p$ and if $\gm$
were increasing, then 
$$\gm((1-\ge)R)\leq c_5\ge^{-\frac hp}R^{-\frac hp}\gm^\frac qp(R),
$$
and by \rlemma{ST0},
$$\gm(R)\leq c_6r^{-\frac {h}{p-q}}\quad\text{ for } R\geq r_2,
$$
contradiction.\\
 Hence $\gm$ is decreasing and tends to $0$ at infinity by $(\ref{hoch})$. Furthermore $(\ref{hoch})$ implies
\bel{as*}\gm^p((1+\ge)R)\leq C\ge^{-h}R^{-\tilde h}\gm(R)\,\text{ and thus }\;\gm((1+\ge)R)\leq C\ge^{-\frac{h}{p}}R^{-\frac{\tilde h}{p}}\gm^\frac{1}{p}(R)
\ee
with $\tilde h=\min\{2,q\}$. Applying again \rlemma{ST0} we deduce
    \bel{as**}\dsps
\gm(R)\leq c_7R^{-\frac{\tilde h}{p-1}}.
\ee
Note that if $q\geq 2$, $\frac{\tilde h}{p-1}=\ga$ and we obtain $(\ref{as5})$. If $1<q< 2$, then $\tilde h=q$  and $\frac{\tilde h}{p-1}=\frac{q}{p-1}$ and we encounter two possibilities:\smallskip
 
\nind (a) if $\frac{q}{p-1}\geq \gb$, then (\ref{as7}) implies 
$$\gm(R)\leq c_8R^{-\gb},$$
 and by the equivalence in (\ref{mas}) 
 $$\frac{\gm^q(R)}{R^q}\leq c_8^{1-q}\frac{\gm(R)}{R^2},$$ 
which in turn implies 
$$\gm^p(R(1+\ge))\leq 2c_8\ge^{-2}\frac{\gm(R)}{R^2}.$$
By \rlemma{ST0} we obtain $(\ref{as5})$. \smallskip

This holds in particular when $1<p\leq q<2$ which completes the proof of 1-(iii).\smallskip

\nind (b) Let   $A_0=\frac{q}{p-1}< \gb$. For any $0<A\leq \gb$ and $\gm(R)\leq c_9A^{-A}$ we have that
$$\gm^p(2R)\leq c_{10}\left(R^{-(A+2)}+R^{-(A+1)q}\right)=c_{10}R^{-(A+1)q}\left(1+R^{A(q-1)-(2-q)}\right)\leq 2c_{10}R^{-(A+1)q},
$$
so $\gm(2R)\leq c_{11}R^{-\frac{(A+1)q}{p}}$. We define a sequence $\{A_n\}$ by $A_0=\frac{q}{p-1}$ and 
\bel{suit}A_n=\frac{(A_{n-1}+1)q}{p}\quad\text{for }n\geq 1.
\ee
Then,  as long as $A_{n-1}\leq \gb$, we have
$$\gm(2^nR)\leq C_nR^{-A_n}.
$$
Furthermore $A_1-A_0=\frac{q(q-1)}{p(p-1)}$ and $A_n-A_{n-1}=\frac{q(A_{n-1}-A_{n-2})}{p}$. Therefore the sequence $\{A_n\}$ is increasing.\smallskip

\nind {\it Proof of 1-(i)}. For $q>\frac{2p}{p+1}$ we have $\gb<\ga<\gg$. If $A_{n-1}<\gb$ for any $n\geq 1$ the sequence $\{A_n\}$ converges to $\gg$, contradiction. Therefore there exists $n_0\geq 1$ such that $A_{n_0+1}\geq\gb$, so we conclude as in case (a).\smallskip

\nind {\it Proof of 1-(ii)}. If $1<q\leq \frac{2p}{p+1}$, then $\gg<\ga<\gb$, and $A_0<\gg\leq\gb$ since $q>1$. So the sequence $\{A_n\}$ is still increasing 
and it converges to $\gg$.
This implies that for any $\gth>0$, there exists $C_\gth$ such that 
$$\gm(R)\leq C_\gth R^{-\gg+\gth} \quad\text{for }R\geq 2r_0.
$$
Set $g(r)=r^{-\gg}$, then
$$g^p(R(1+\ge))\leq R^{-p\gg}\leq\ge^{-q}\frac{g^q(R)}{R^q},
$$
since $\gg=\frac q{p-q}$. Recalling that
$$\gm^p(R(1+\ge))\leq c_4\ge^{-q}\left(\frac{\gm(R)}{R^2}+\frac{\gm^q(R)}{R^q}\right),
$$
and putting $\gf(R)=\max\{g(R),\gm(R)\}$ we obtain
$$\gf(R(1+\ge))\leq c_{12}\ge^{-q}\left(\frac{\gm(R)}{R^2}+\frac{\gm^q(R)}{R^q}+\frac{g^q(R)}{r^q}\right)\leq c_{13}\ge^{-q}\left(\frac{\gf(R)}{R^2}+\frac{\gf^q(R)}{R^q}\right).
$$
Because $\gf(R)\geq g(R)\geq R^{-\gb}$ as $\gg\leq\gb$, we have $\frac{\gf(R)}{R^2}\leq\frac{\gf^q(R)}{R^q}$, hence
$$\gf(R(1+\ge))\leq c_{14}\ge^{-\frac qp}R^{-\frac qp}\gf^{\frac qp}(R).
$$
It follows from \rlemma{ST0}-$(\ref{as2})$-$(\ref{as3})$ that $\gf(R)\leq c_{15}R^{-\gg}$. This is $(\ref{as6})$.
\smallskip

\nind {\it 2- The problem in $B_{r_0}\setminus\{0\}$}. By \rlemma{ST2}, $\gm$ is nonincreasing and $(\ref{hach})$ holds. If 
$\gm$ is bounded, then it admits a positive limit at $0$ and the two estimates in {\it 2} hold. Hence we assume that $\gm(R)\to\infty$ as 
$R\to 0$. From $(\ref{hoch})$
$$\gm^p(R(1-\ge))\leq c_4\ge^{-h}\left(\frac{\gm(R)}{R^2}+\frac{\gm^q(R)}{R^q}\right),
$$
where, we recall it, $h=\max\{2,q\}$. We notice that if $(\ref{mas})$ holds, then 
$$\gm((1+\ge)R)\leq c_{16}R^{-\frac 2p}\gm^\frac{1}{p}(R)\Longrightarrow \gm(R)\leq C'R^{-\ga},
$$
which is the desired estimate in the case $1<q\leq \frac {2p}{p+1}$. We notice also that the fact that $\gm(R)\to\infty$ as $R\to 0$ implies 
$$\gm^p(R(1+\ge))\leq c_4\ge^{-h}\left(\frac{1}{R^2}+\frac{1}{R^q}\right)\gm^q(R)\leq 2c_4\ge^{-h}R^{-h}\gm^q(R),
$$
which in turn yields
\bel{N1}
\gm(R)\leq c_{17} R^{-\frac{h}{p-q}} \quad\text{for }0<R\leq r_1<r_0.
\ee
Hence, if $h=q$, we obtain $(\ref{as7})$.\smallskip

\nind {\it Proof of 2-(i)}. Let  $2>q\geq \frac{2p}{p+1}$. Then $\gb\leq\ga\leq\gg$, then we start with $\gm(R)\leq R^{-A_0}$ with $A_0=\frac{2}{p-q}>\gg$. For any $A>0$ larger than $\gg$ and such that $\gm(R)\leq c_{18} R^{-A}$, there holds
$$\gm^p(\tfrac R2)\leq c_{19} R^{-(1+A)q},
$$
as above since $A>\gb$. The sequence $\{A_n\}$ still defined by $(\ref{suit})$ satisfies 
$$\gm\left(\tfrac R{2^n}\right)\leq c_nR^{-A_n}
$$
as long as $A_{n-1}>\gb$. We have $A_1-A_0=\frac{q-(p-q)A_0}{p}<0$. Since $A_{n+1}-A_n=\frac{q}{p}(A_{n}-A_{n-1})$, the sequence $\{A_n\}$ is decreasing and it converges to $\gg$. We adapt the technique developed in {\it 1-(ii)}: for any $\gth>0$ there exists $C_\gth>0$ such that 
$$\gm(R)\leq C_\gth R^{-\gg-\gth}\quad\text{for }0<R\leq\frac {r_0}{2}.
$$
Defining $g(R)=R^{-\gg}$ and $\gf(R)=\max\{g(R),\gm(R)\}$, then we obtain
$$\gf^p(R(1-\ge))\leq c_{20} \ge^{-h}\left(\frac{\gm(R)}{R^2}+\frac{\gm^q(R)}{R^q}+\frac{g^q(R)}{R^q}\right)\leq c_{21} \ge^{-h}\left(\frac{\gf(R)}{R^2}+\frac{\gf^q(R)}{R^q}\right)
$$
Because $\gg>\gb$ we have $R^{-\gb}\leq R^{-\gg}\leq \gf(R)$ for $0<R\leq 1$ which implies that $\frac{\gf(R)}{R^2}\leq \frac{\gf^q(R)}{R^q}$ and 
$$\gf^p(R(1-\ge))\leq 2 c_{21}\ge^{-h}\frac{\gf^q(R)}{R^q}.
$$
It follows by \rlemma {ST0} that $\gf(R)\leq c_{22} R^{-\gg}$ and $(\ref{as7})$.\smallskip

\nind {\it Proof of 2-(ii)}. If $1<q< \frac{2p}{p+1}$. Then $\gg<\gb<\ga$. We proceed as in case {\it 2-(i)} with the same sequence $\{A_n\}$. We notice that $A_0=\frac{2}{p-q}>\ga>\gg$ since $q>1$. Then $A_1<A_0$ and as above $\{A_n\}$ is nonincreasing and converges to $\gg$. As in the proof of {\it 1-(i)}
 there exists an integer $n_0$ such that $A_{n_0}\leq\gb$ which in turn implies $(\ref{mas})$, and finally $(\ref{as8})$ holds.
\qeda\medskip

\nind\Remark From \rth{Theorem C} we recover easily the result of {Theorem C}-({\it 2}). Indeed, if $\dsps f(r)>cr^p$ for $c>0$ and $r\geq r_1$ and $1<q<p$, any positive supersolution $u$ of $(\ref{as3-})$ in $B_{r_1}^c$ such that $\dsps \lim_{|x|\to\infty}u(x)=\infty$ is a supersolution of 
$$-\Gd u+m|\nabla u|^q=cu^p
$$
in this domain. Then $\dsps \lim_{r\to\infty}\gm(r)=0$ from the upper estimates of \rth{Theorem C}, contradiction.
 \subsection{Construction of radial minorant solutions in the exterior problems}


The next result extends the construction of  \cite[Theorem 1.3]{Bi} and brings precisions to \cite[Lemma 4]{ABGMQ} that we recall below.\smallskip

\nind{\it Assume $N \geq 2$, $q>1$ and let $f:(0,\infty)\mapsto\BBR$ be positive, nondecreasing and continuous. Suppose there exists a positive supersolution $u$ of problem $(\ref{as3})$ below. Then there exists a positive radial supersolution $v$ of $(\ref{as3})$. In addition, if $u$ does not blow up at infinity, then $v$ is bounded, while if $u$ blows up at infinity, $v$ is bounded from below.}\medskip

Our result is the following.

\bth{Subth} Let $q>1$, $m>0$ and $f:\BBR_+\mapsto\BBR_+$ be a Lipschitz continuous function satisfying assumption (F). Suppose that there exists a positive $C^2(\overline B^c_{r_0})$ function $u$ satisfying 
 \bel{as3}\dsps
-\Gd u+m|\nabla u|^q-f(u)\geq 0\quad\text{in }B^c_{r_0},
\ee
then there exists a positive radial and monotone function $v\in C^2(\overline B^c_{r_0})$ smaller than $u$ satisfying
 \bel{as4}\dsps
-\Gd v+m|\nabla v|^q-f(v)= 0\quad\text{in }B^c_{r_0},
\ee
such that:\smallskip

\nind 1- $\dsps v(r_0)=\min_{|x|=r_0}u(x)$ and $\dsps\lim_{r\to\infty}v(r)=\infty$, when $\dsps\lim_{|x|\to\infty}u(x)=\infty$.\smallskip

\nind 2- $\dsps 0<v(r_0)=a\leq \min_{|x|=r_0}u(x)$ and $\dsps\lim_{r\to\infty}v(r)=0$, when $\dsps\liminf_{|x|\to\infty}u(x)=0$, under the additional condition 
when $q>2$, 
  \bel{sub1***}\dsps
  a<\Gth:=\left(\frac{q(N-1)-N}{m(q-1)}\right)^{\frac{1}{q-1}}\!\!\!\!r_0^{2-N}\int_1^{\frac{\gt}{r_0}}\!t^{1-N}\left(1-t^{N-q(N-1)}\right)^{-\frac{1}{q-1}}dt.
  \ee
  \es
 \Proof The proof is based upon an iterative process reminiscent of a method used in \cite{Bi}. However the technicalities are much more involved and developed in the Appendix. By \rlemma {LX} a positive supersolution $u$ in an exterior domain either tends to $\infty$ at $\infty$ or satisfies $\dsps\liminf_{|x|\to\infty}u(x)=0$.\smallskip

\nind For $\gt>r_0$ we set $\dsps \tilde b_0=\inf_{|x|=r_0}u(x)$ and $\dsps b_\gt=\inf_{|x|\geq\gt}u(x)$. If $0<a\leq \tilde b_0$ and $0\leq b\leq b_\gt$ we consider the sequence of radially symmetric functions defined in $B_\gt\cap B_{r_0}^c$ functions $\{v_{k,\gt}\}_{k\in\BBN}$ such that  $v_{0,\gt}\equiv 0$, and for $k\geq 1$
\bel{E5}\BA{lll}
-\Gd v_{k,\gt}+m|\nabla v_{k,\gt}|^q=f(v_{k-1,\gt})\quad&\text{in } B_\gt\cap B_{r_0}^c\\
\phantom{-\Gd +m|\nabla v_{k,\gt}|^q}
v_{k,\gt}=b&\text{in } \prt B_\gt\\
\phantom{-\Gd +m|\nabla v_{k,\gt}|^q}
\dsps v_{k,\gt}=a&\text{in } \prt B_{r_0}.
\EA
\ee
If $1<q\leq 2$, the function $v_{1,\gt}$ exists without any restriction on $a$ and $b$.\\
If $q>2$ we have existence if $a\leq b$ provided $\gt\geq\gt^*$ where $\gt^*$ is defined in \rlemma{Ap1} (2), and if $a>b$ the condition for existence is 
$$b<a<b+\Gth.
$$
In both case, the function $v_{1,\gt}$ is positive, monotone. and dominated by $u$. \\
 Next for $k=2$ we apply the extension \cite[Corollary 1.4.5]{Vebook} of the classical result 
\cite[Th\'eor\`eme 2.1]{BMP}. The function $v_{1,\gt}$ satisfies equation $(\ref{E5})$ with right-hand side $0$ instead of $f(v_{k-1,\gt})$. By the maximum principle it is dominated by the supersolution $u$, thus $f(u)\geq f(v_{1,\gt})$. Then there exists a function $v_{2,\gt}$ which satisfies $(\ref{E5})$ with $k=2$ and 
$$v_{1,\gt}\leq v_{2,\gt}\leq u.
$$
Note that this function is unique by the maximum principle. We introduce there the spherical coordinates $(r,\gth)\in\BBR_+\ti S^{N-1}$ in $\BBR^N$. Let $\bar v_{2,\gt}(r)$ be the spherical average of $v_{2,\gt}(r,.)$ on $S^{N-1}$. Since $f(v_{1,\gt})$ is radial, by convexity, $\bar v_{2,\gt}$ satisfies 
$$\BA{lll}
-\Gd \bar v_{2,\gt}+m|\nabla \bar v_{2,\gt}|^q\leq f( v_{1,\gt})\quad&\text{in } B_\gt\cap B_{r_0}^c\\
\phantom{-\Gd +m|\nabla\bar v_{k,\gt}|^q}
\bar v_{2,\gt}=b&\text{in } \prt B_\gt\\
\phantom{-\Gd +m|\nabla\bar v_{k,\gt}|^q}
\dsps\bar v_{2,\gt}=a&\text{in } \prt B_{r_0}.
\EA
$$
By the maximum principle we have $\bar v_{2,\gt}(r)\leq v_{2,\gt}(r,\gth)$ for any $r$ and any $\gth$, which implies that $\bar v_{2,\gt}= v_{2,\gt}$, hence $v_{2,\gt}$ is spherically symmetric. Iterating this process, we construct the increasing the
sequence $\{v_{k,\gt}\}_{k\in\BBN}$ of positive spherically symmetric solutions of $(\ref{E5})$ dominated by $u$ in $B_\gt\cap B_{r_0}^c$. For $k\geq 2$ the function  $v_{k,\gt}$ cannot have a local minimum, hence if $a\leq b$ it is monotone increasing (as a function of $|x|$) and if   $a>b$, it is decreasing for $|x|$ close to $\gt$.\\
Since the sequence $\{v_{k,\gt}\}_{k\in\BBN}$ is increasing and $v_{k,\gt}\leq u$, it converges to some radial positive function $v_{\infty,\gt}:=v_\gt$ by Ascoli theorem and $v_\gt$ is a positive $C^2$ solution of 
\bel{E5+1}\BA{lll}
-\Gd v_{\gt}+m|\nabla v_{\gt}|^q=f(v_{\gt})\quad&\text{in } B_\gt\cap B_{r_0}^c\\
\phantom{-\Gd +m|\nabla v_{\gt}|^q}
v_{\gt}=b&\text{in } \prt B_\gt\\
\phantom{-\Gd +m|\nabla v_{\gt}|^q}
\dsps v_{\gt}=a&\text{in } \prt B_{r_0}.
\EA
\ee
If $a\geq b$ then necessarily $v_{k,\gt}\leq v_{k,\gt'}$ in $B_\gt\cap B_{r_0}^c$ otherwise 
$v_{k,\gt'}$ would have a local minimum in $B_{\gt'}\cap B_{r_0}$. \smallskip

\nind {\it Assertion 1}.  Here $\gm(r)\to\infty$ when $r\to\infty$. Let $r_1>r_0$ such that $\dsps b_\gt>\min_{|x|=r_0}u(x)$ for all $\gt\geq r_1$.
Let  $v_{\infty,\gt}:=v_\gt$ be the solution of $(\ref{E5+1})$ with $a=\min_{|x|=r_0}u(x)$ and $b=\gb_{r_1}$ and $\gt>\gt^*$ if $q>2$, which is not a restriction since we aim to let $\gt\to\infty$. Since $v_\gt$ cannot have any local minimum in 
 $B_{\gt}\cap B^c_{r_0}$, we have 
 $$a\leq v_\gt(|x|)\leq u(x)\quad\text{for all }x\in B_{\gt}\cap B^c_{r_0}.
$$
By standard ODE techniques, for any $T>r_1$, $v_\gt$ is bounded in $C^3(\overline B_{T}\cap B^c_{r_0})$ uniformly with respect to $\gt\geq T+1$. Hence there exists a sequence $\{\gt_n\}$ tending to infinity and a radially symmetric positive function $v\in C^2{B^c_{r_0}}$ such that 
\bel{E5+2}\BA{lll}
-\Gd v+m|\nabla v|^q=f(v)\quad&\text{in } B_{r_0}^c\\
\phantom{-\Gd +m|\nabla v|^q}
\dsps v=a&\text{in } \prt B_{r_0}.
\EA
\ee
Furthermore $a\leq v\leq u$. By \rlemma {LX} $v(r)\to\infty$ when $r\to\infty$ which proves 1.\smallskip

\nind {\it Assertion 2}.  We solve $(\ref{E5+1})$ with $b=0$ and $a\leq \min_{|x|=r_0}u(x)$ with the additional condition $a<\Gth$ if $q>2$ and we set 
$v_{\infty,\gt}:=v_\gt$. Then $0\leq v_\gt\leq a$ and since the function $v_\gt$ cannot have a local minimum in $(r_0,\gt)$, we have also that 
$$v_\gt(|x|)\leq v_{\gt'}(|x|)\leq u(x)\quad\text{for all }\gt'>\gt\,\text{ and }x\in B_\gt\cap B_{r_0}^c.
$$
Letting $\gt\to\infty$ we obtain that $v_\gt$ converges in the local $C^2(B_{r_0}^c)$-topology to some $v\in C^2(B_{r_0}^c)$, which satisfies $(\ref{E5+2})$ and 
$v(|x|)\leq u (x)$ for $x\in B_{r_0}^c$. Therefore $v(r)\to 0$ as $r\to\infty$ and we complete the proof of  {\it 2}.
\qeda
\bcor{Subcor} Let $N\geq 2$, $m>0$, $q> \frac N{N-1}$ and $f$ be as in \rth{Subth}. Then any positive $C^2(\overline B^c_{r_0})$ function $u$ verifying 
$(\ref{as3})$ satisfies
\bel{crit1S}
u(x)\geq c|x|^{2-N}\quad\text{for all }\;x\in B_{r_0}^c
\ee
for some $c>0$.
\es
\Proof For $r_0<\gt$, we introduced the function $v_{1,\gt}$  which satisfies 
$$\BA{lll}\dsps - v''_{1,\gt}-\frac{N-1}{r} v'_{1,\gt}+m|v_{1,\gt}|^q=0\qquad\text{in }(r_0,\gt)\\
\phantom{\dsps-v_{1}-\tfrac{N.,,.}{r},v_{1\gt}+m||^q}
\! v_{1,\gt}(r_0)=a\\[1mm]
\phantom{\dsps-v_{1}-\tfrac{N.,,.}{r},,v_{1\gt}+m||^q}
\! v_{1,\gt}(\gt)=0
\EA$$
with $\dsps 0<a\leq \min_{|x|=r_0}u(x)$. We have seen therein that $ v_{1,\gt}(|x|)\leq u(x)$ for $x\in B_\gt\setminus B_\gr$. If $q>2$ we choose $a\leq\Gth$. When $\gt\to\infty$, 
$ v_{1,\gt}\uparrow v_{1,\infty}$ and $v:= v_{1,\infty}(|x|)\leq u(x)$ in $B_{r_0}^c$. Since $v'\leq 0$, we have
$$v''+v^p=m|v'|^q-\frac{N-1}{r}v'\geq 0.
$$
then 
$$E(r):=\left(\frac{v'(r)^2}{2}+\frac{v(r)^{p+1}}{p+1}\right)'\leq 0.
$$
Therefore $E(r)$ admits a limit when $r\to\infty$. Because $v(r)\to 0\geq 0$, this implies that $v'(r)$ admits also a limit $\ell\leq 0$ when $r\to\infty$ and this, limit is necessarily $0$ since $v$ is bounded. \\
Set $w(r)=-r^{N-1}v'$, then $w\geq 0$ and
$$w'+mr^{(1-q)(n-1)}w^q\geq 0.
$$
Integrating this equation as it is done in Appendix, we obtain
$$(w^{1-q})'(r)+\frac{m(q-1)}{q(N-1)-N}m(r^{(N-q(N-1)})'\leq 0,
$$
which implies by integration 
$$w^{1-q}(r)-w^{1-q}(r_1)\leq \frac{m(q-1)}{q(N-1)-N}\left(r_1^{N-q(N-1)}-r^{N-q(N-1)}\right).
$$
Therefore $w(r)\geq c_1>0$ and $v'(r)\geq -c_1r^{1-N}$ and thus $v(r)\geq \frac {c_1}{N-2}r^{2-N}$. Because $u(x)\geq v(r)$ for $|x|=r\geq r_0$ this yields 
$(\ref{crit1S})$.
\qeda\medskip

\nind\Remark As  a consequence we recover Theorem C-(1) in the case $q>\frac N{N-1}$. Indeed, suppose that $f(s)\geq Cs^p$ near $s=0$ and $1<p\leq\frac N{N-2}$. Then if there exists a positive supersolution of $(\ref{as3-})$ which is bounded at infinity, then $\dsps\liminf_{|x|\to\infty}u(x)=0$ by \rlemma{LX}. 
Since $u$ is a supersolution of 
$$-\Gd u+m|\nabla u|^q= Cu^q\quad\text{in }B_{r_1}^c
$$
for some $r_1>r_0$, by \rth{Theorem C} and \rcor{Subcor} there exists a positive radially symmetric solution $v$ of the above equation such that
$$u(x)\geq v(|x|)\geq c|x|^{2-N}\quad\text{for all  }x\in B_{r_1}^c.
$$
By \rth{Theorem C} we have also $\gm(|x|)\leq C|x|^{-\ga}$ in $B_{r_1}^c$. This is a contradiction when $p>\frac{N}{N-2}$. When $p=\frac{N}{N-2}$ we set 
$v(r)=r^{2-N}X(t)$ with $t=\ln r$.
Then 
$c_1\leq X(t)\leq c_2$ for $t\geq t_1=\ln r_1.$ Hence $X$ is a bounded solution of 
$$X''-(N-2)X'+CX^p-me^{(N-q(N-1))t}\left(|(N-2)X-X'|\right)^{q}=0,
$$
and it is straightforward to verify that the $\gw$-limit set of the trajectory $\dsps\CT_+[v]=\bigcup_{t\geq t_1}\{X(t)\}$ is reduced to $\{0\}$, which is still a contradiction.


\subsection{Dichotomy result when $q\geq p$. Proof of \rth{Theorem D}}
In this Section we suppose $q\geq p>1$. Then there exist supersolutions of $(\ref{I-1})$ such that $\dsps\lim_{|x|\to\infty}u(x)=\infty$, e.g. $u(x)=e^{\gl|x|}$ for any $\gl>0$ if $q>p$ or $\gl$ large enough if  $q=p$. 

\nind {\it Proof of \rth{Theorem D}}. Our proof is based upon \rth{Subth} with $f(u)=u^p$.  Let $u$ be a positive supersolution of $(\ref{I-1})$. From \rlemma{LX}, either 
$u(x)\to \infty$ or $\gm(|x|)\to 0$  when $|x|\to\infty$. \smallskip

\nind (i) Suppose that $\dsps\lim_{|x|\to\infty}u(x)=\infty$. By \rth{Subth} there exists a radial and increasing function $v$ below $u$ in $B_{r_1}^c$ satisfying
\bel{crit12}\BA{lll}\dsps
-v''-\frac{N-1}{r}v'+mv'^q=v^p\quad\text{in }\,(r_1,\infty)\\[2mm]
\phantom{----aa--;-}\dsps v(r_1)=\min_{|x|=r_1}u(x)\\
[2mm]
\phantom{------a}\dsps \lim_{r\to\infty}v(r)=\infty.
\EA\ee
For $\ge>0$ we set $F_\ge(r)=v^p(r)-(1+\ge)m(v'(r))^q$. This type of function introduced by \cite{SeZo} is fundamental in the study of  radial soutions. Then
$$F'_\ge(r)=pv'v^{p-1}-q(1+\ge)mv''v'^{q-1}=pv'v^{p-1}+q(1+\ge)mv'^{q-1}\left(\frac{N-1}{r}v'+v^p-mv'^q\right).
$$
If there exists some $r_2>r_1$ such that $F_\ge(r_2)=0$, then 
$$F'_\ge(r_2)=pv'v^{p-1}+q(1+\ge)mv'^{q-1}\left(\frac{N-1}{r_2}v'+\ge mv'^q\right)>0.
$$
This implies that $F_\ge(r)>0$ for all $r>r_2$. As a consequence, $F_\ge(r)$ has a constant sign for $r$ large enough.  \\
When $N\geq 3$ we can take $\ge=0$. If $F_0\leq 0$ for $r>r_2>{r_0}$, then $v^p(r)\leq m(v'(r))^q$ which implies
\bel{crit13}\BA{lll}\dsps
v(r)\geq\left(m|\gg|^q\right)^{\frac{1}{p-q}}(r-r_2)^{|\gg|}\quad\text {for all }r>r_2,
\EA\ee
in the case $q>p$ and 
\bel{crit14}\BA{lll}\dsps
v(r)\geq v(r_2)e^{m^{-\frac 1m}(r-r_2)}\quad\text {for all }r>r_2,
\EA\ee
when $q=p$. This yields $(\ref{I-9})$.\\
If $F_0\geq 0$ for $r>r_2>{r_0}$, then $\Gd v\leq 0$ if $|x|>r_2$, and the function $r^{N-1}v'r)$ is nonincreasing on $[r_2,\infty)$, thus $v'(r)\leq cr^{1-N}$. If 
$N\geq 3$, it implies that $v(r)$ remains bounded, which is a contradiction. \\
When $N=2$ we take $\ge=1$. If $F_1(r_3)=0$ for some $r_3$, then either $F_1$ is positive  for $r\geq r_3$, which implies 
$$-2v''=\frac {1}{r}v'+v^p+F_2(r)\geq v^p\quad\text{for }r\geq r_2.
$$
In such a case, we deduce by multiplying by $v'\geq 0$ that the function $r\mapsto \left(v'^2+\frac{v^{p+1}}{p+1}\right)(r)$ is nonincreasing, hence bounded, contradiction. If this does not hold, then $F_1$ is nonpositive  for $r\geq r_3$, which yields
\bel{crit15}v(r)\geq \left\{\BA{lll}\left(2m|\gg|^q\right)^{\frac{1}{p-q}}(r-r_2)^{-\gg}\quad&\text{if $r\geq r_2$ when $N\geq 3$}\\[2mm]
v(r_2)e^{(2m)^{-\frac 1{2m}}(r-r_2)}\quad&\text{if $r\geq r_2$ when $N=2$}.
\EA\right. 
\ee
If we have now $F_0(r)>0$, then $v'(r)\leq cr^{-1}$ which implies $v(r)\leq c\ln r+d$, which is not compatible with $(\ref{crit15})$. Therefore $F_0(r)\leq 0$ which again implies that $(\ref{I-9})$ holds.\smallskip

\nind \nind (ii) Assume now that $\dsps\lim_{r\to\infty}\gm(r)=0$. Inequality $(\ref{I-10})$-(a)  follows from \rth{Theorem C} (1-iii). Since $q>p>\frac{N}{N-2}$ we have 
$q>\frac{N}{N-1}$. Thus $(\ref{I-10})$-(b)  is a consequence of \rcor{Subcor}.\qeda



\mysection{Estimates on solutions}


\subsection{General estimates}

A major tool for proving a priori estimates either near an isolated singularity or at infinity is the Keller-Osserman combined with Bernstein method applied to the function $z=|\nabla u|^2$. We recall the variant of Keller-Osserman a priori estimate that we proved in \cite{BiGaVe2}.
\blemma{KOB} Let  $q>1$ $d\geq 0$ and $P$ and $Q$ two continuous functions defined in $B_\gr(a)$ such that $\inf\{P(y):y\in B_\gr(a)\}>0$ and 
$\sup\{Q(y):y\in B_\gr(a)\}<\infty$. If $z$ is a positive $C^1$ function defined in $B_\gr(a)$ and such that 
\bel{Es1}
-\Gd z+P(y)z^q\leq Q(y)+d\frac{|\nabla z|^2}{z}\quad\text{in }B_\gr(a),
\ee
then there exists a positive constant $C=C(N,q,d)>0$ such that 
\bel{Es2}\dsps
z(x)\leq C\left(\left(\frac{1}{\gr^2}\frac{1}{\dsps\inf_{B_\gr(a)}P}\right)^{\frac1{q-1}}+\left(\sup_{B_\gr(a)}\frac{Q}{P}\right)^{\frac1{q}}\right)\quad\text{for all }x\in B_{\frac\gr2}(a).
\ee
\es

 In the next statement we show how an upper estimate on $u(x)$ by a power of $|x|$ implies a precise estimate on $|\nabla u(x)|$. 
 \bth{pre} Let $p,q>1$, $m>0$ and $r_0>0$.\smallskip
 
 \nind 1- If $u$ is a positive solution of $(\ref{I-1})$ in $B_{r_0}\setminus\{0\}$ where it satisfies
 \bel{Es5*}
 |x|^\gl u(x)\leq c
 \ee
for some constant $c>0$ and some exponent $\gl>0$, then there exists $c_1=c_1(N,p,q,\gl,c)>0 $ such that 
 \bel{Es5}\dsps
|\nabla u(x)|\leq c_1\left(|x|^{-\frac{1}{q-1}}+|x|^{-\frac{\gl p}{q}}+|x|^{-\frac{\gl(p-1)}{2(q-1)}}\right)\quad\text{for all }x\in B_{\frac{r_0}2}\setminus\{0\}.
\ee
Furthermore, when $1<q\leq 2$, one has an improvement of $(\ref{Es5})$ under the form
 \bel{Es6}\dsps
 |\nabla u(x)|\leq c'_1|x|^{-(\gl+1)}\quad\text{for all }x\in B_{\frac{r_0}2}\setminus\{0\},
\ee
for any $\gl>0$ such that $\gl\leq\min\{\ga,\gb\}$.\smallskip

\nind 2- If $u$ is a positive solution of $(\ref{I-1})$ in $B^c_{r_0}$, then 
 \bel{Es8}\dsps
 \limsup_{|x|\to\infty}u(x)<\infty\Longrightarrow  \limsup_{|x|\to\infty}|\nabla u(x)|<\infty,
 \ee
  \bel{Es9}\dsps
 \lim_{|x|\to\infty}u(x)=0\Longrightarrow  \lim_{|x|\to\infty}|\nabla u(x)|=0.
 \ee
 If $u$ satisfies $(\ref{Es5*})$  in $B^c_{r_0}$ for some $c>0$ and $\gl>0$, then there exists 
 $c_1:=c_1(N,p,q,\gl,c)>0$  such that 
  \bel{Es10}\dsps
|\nabla u(x)|\leq c_1\left(|x|^{-\frac{1}{q-1}}+|x|^{-\frac{\gl p}{q}}+|x|^{-\frac{\gl(p-1)}{2(q-1)}}\right)\quad\text{for all }x\in B^c_{2r_0}.
\ee
 Furthermore, if $1<q\leq 2$, one has an improvement of $(\ref{Es10})$ under the form 
   \bel{Es11*}\dsps
 |\nabla u(x)|\leq c_2|x|^{-(\gl+1)}\quad\text{for all }x\in B^c_{2r_0},
\ee
 for $c_2:=c_2(N,p,q,\gl,c)>0$ for any $\gl\geq\max\{\ga,\gb\}$.
   \es
  \Proof We use Bernstein method, setting $z(x)=|\nabla u(x)|^2$ and Weitzenb\"ock's formula
  $$-\frac 12\Gd z=|D^2u|^2+\langle\nabla (\Gd u),\nabla u\rangle. $$  
Using  the inequality $|D^2u|^2\geq\frac 1N(\Gd u)^2$ and the equation satisfied by $u$ we obtain
$$\dsps
-\frac 12\Gd z+\frac 1N(mz^{\frac q2}-u^p)^2+\langle\nabla (mz^{\frac q2}-u^p),\nabla u\rangle\leq 0.
$$
Developing this inequality yields
$$-\frac 12\Gd z+\frac {m^2}Nz^q+\frac 1Nu^{2p}\leq\frac {2m}Nu^pz^{\frac q2}+pu^{p-1}z+\frac {mq}2z^{\frac q2-1}\langle\nabla z,\nabla u\rangle.
$$
Now for $\ge>0$
$$z^{\frac q2-1}\langle\nabla z,\nabla u\rangle= z^{\frac q2-\frac12}\langle\frac{\nabla z}{\sqrt z},\nabla u\rangle\leq z^{\frac q2}\frac{|\nabla z|}{\sqrt z}\leq \ge z^q+\frac 1\ge \frac{|\nabla z|^2}{z},
$$
$$u^{p-1}z\leq \ge z^q+\ge^{-\frac 1{q-1}}u^{\frac{q(p-1)}{q-1}},
$$
and
$$u^pz^{\frac q2}\leq \ge z^q+\frac 1\ge u^{2p}.
$$
We choose $\ge$ small enough and get
  \bel{Es11**}\dsps
-\Gd z+\frac{m^2}{N}z^q\leq c_{3}\frac{|\nabla z|^2}{z}+c_{4}u^{2p}+c_{5}u^{\frac{q(p-1)}{q-1}}
\ee
where $c_i=c_i(N,p,q,m)>0$, $i=3,4,5$. We Apply \rlemma{KOB} in $\overline B_{2\gr}(a)$, with $\overline B_{2\gr}(a)\subset B_{r_0}\setminus\{0\}$ in case 1, or  $\overline B_{2\gr}(a)\subset \overline B^c_{r_0}$ in case 2, we obtain for some positive constant $c_{6}:=c_{6}(N,q,m)>0$,
  \bel{Es12*}\dsps
\sup_{B_{\gr}(a)}z(y)\leq c_{6}\left(\gr^{-\frac 2{q-1}}+\sup_{B_{2\gr}(a)}\left(u^{2p}+u^{\frac{q(p-1)}{q-1}}\right)^{\frac 1q}\right), 
\ee
which is equivalent to 
  \bel{Es13}\dsps
\sup_{B_{\gr}(a)} |\nabla u(z)|\leq c_{7}\left(\gr^{-\frac 1{q-1}}+\sup_{B_{2\gr}(a)}\left(u^{\frac pq}+u^{\frac{p-1}{2(q-1)}}\right)\right),
\ee
where $c_{7}=c_{7}(N,q,m,c_{6})>0$. \smallskip

\nind 1- Next we assume that $u(x)\leq c_{8}|x|^{-\gl}$ in $B_{r_0}\setminus\{0\}$. Then $(\ref{Es13})$ yields
  exactely $(\ref{Es5})$ with $c_9=c_9(N,m,p,q,\gl,c_{8})>0$.\smallskip
  
 In some cases we can obtain a different estimate which requires $1<q\leq 2$. For $k>0$ we set
$$u_k(x)=k^\gl u(kx).
$$
Then $u_k$ satisfies
  \bel{Es15-k}\dsps
-\Gd u_k+mk^{\gl+2-q(\gl+1)}|\nabla u_k|^q-k^{\gl+2-\gl p}u_k^p=0\quad\text{in } B_{k^{-1}r_0}.
\ee
The function $u_k$ is uniformly bounded in the spherical shell $\Gg_{\frac {r_0}8,\frac {2r_0}3}:=\left\{x:\frac {r_0}8\leq |x|\leq \frac {r_0}2\right\}$. If we assume that 
  \bel{Es17}\dsps
\gl+2-q(\gl+1)\geq 0\Longleftrightarrow \gl\leq\tfrac{2-q}{q-1}=\gb\quad\text{and }\gl+2-\gl p\geq 0\Longleftrightarrow \gl\leq\tfrac{2}{p-1}=\ga,
\ee
then we deduce from standard regularity estimates \cite{GT} (this is why we need $1<q\leq 2$) that 
  \bel{Es18}\dsps
|\nabla u_k(x)|\leq c_{9}\Longleftrightarrow |\nabla u(kx)|\leq c_{9}k^{-\gl-1}\quad\text{for all }x\in \Gg_{\frac {r_0}4,\frac {r_0}2}.
\ee
This implies in particular 
  \bel{Es19}\dsps
|\nabla u(x)|\leq c_{9}|x|^{-\gl-1}\quad\text{for all }x\in B_{\frac{r_0}4}\setminus\{0\}.
\ee
Now, this estimate is better than the one in $(\ref{Es5})$ if and only if $\gl\leq\min\{\ga,\gb\}$ and 
  \bel{Es20*}\gl+1\leq\max\left\{\frac 1{q-1},\frac{\gl p}{q},\frac{\gl(p-1)}{2(q-1)}\right\},
\ee
that means 
\bel{Es20}
\gl\leq\gb, \text{ or }\left(q<p\text{ and } \gl>\gg\right),\text{ or } \left(q<\tfrac{p+1}{2}\text{ and } \gl>\tfrac{2(q-1)}{p+1-2q}\right).\ee 
Hence it is an improvement for any $\gl\leq\min\{\ga,\gb\}$.
\smallskip

\nind 2- We apply $(\ref{Es13})$ for $|a|>\gr/2$ with $\gr=\frac {|a|}4$, then we get 
$$\dsps |\nabla u(a)|\leq c_{10}\left(|a|^{-\frac{1}{q-1}}+\max_{|x|\geq \frac{|a|}{2}}\left( u^{\frac pq}+u^{\frac {p-1}{2(q-1)}}\right)\right).
$$
Clearly $(\ref{Es8})$ and $(\ref{Es9})$ follow.

\nind Next we assume $1<q\leq 2$ and $u(x)\leq c_{10}|x|^{-\gl}$ in $B_{r_0}^c$, then $(\ref{Es13})$ yields precisely $(\ref{Es10})$.\smallskip

\nind Again the function $u_k$ defined previously is uniformly bounded in the spherical shell $\Gg_{\frac {3r_0}{2},4r_0}$. In order to apply the standard elliptic equations regularity results to $(\ref{Es15-k})$, we need again $1<q\leq 2$ and
 \bel{Es23*}\dsps
\gl+2-q(\gl+1)\leq 0\Longleftrightarrow \gl\geq\gb\quad\text{and }\,\gl+2-\gl p\leq 0\Longleftrightarrow \gl\geq\ga,
\ee
This yields 
  \bel{Es24*}\dsps
|\nabla u(x)|\leq c_{11}|x|^{-\gl-1}\quad\text{for all }x\in B^c_{2r_0}.
\ee
This estimate is an improvement of $(\ref{Es10})$ if $\gl\geq\max\{\ga,\gb\}$ and 
  \bel{Es25}\gl+1\geq\min\left\{\frac 1{q-1},\frac{\gl p}{q},\frac{\gl(p-1)}{2(q-1)}\right\}.
\ee
That means
\bel{Es26-}
\gl\leq\gb, \text{ or }\left(q\geq p\text{ and } \gl\leq \gg\right),\text{ or } \left(q<\tfrac{p+1}{2}\text{ and } \tfrac{\gl (p+1-2q)}{2(q-1)}<1\right).\ee 
Hence it is an improvement for any $\gl\geq\max\{\ga,\gb\}$.\qeda
\subsection{Upper estimates on solutions when $q>p$. Proof of \rth {Theorem E}}


\nind
\nind{\it Proof of \rth {Theorem E}}. We apply \rlemma{KOB}.\smallskip

\nind 1- {\it Proof of 1-} By change of scale we can assume that ${r_0}=1$. For $0<\gth <\frac14$ we set $\Gw_{\gth}=B_{1-\gth}\setminus B_{\gth}$. For $0<\ge<\frac12$, we have by $(\ref{Es13})$
\bel{Es9*}\dsps
\max_{\overline\Gw_{\gth}}|\nabla u|\leq C\left(\left(\frac1{\gth\ge}\right)^{\frac1{q-1}}+\max_{\overline\Gw_\frac{\gth}{1+\ge}}
\left(u^{\frac pq}+u^{\frac{p-1}{2(q-1)}}\right)\right),\ee
and  $u^{\frac{p-1}{2(q-1)}}\leq u^{\frac pq}+1$ since $q>\frac{2p}{p+1}$. Hence
$$\dsps\max_{\overline\Gw_{\gth}}|\nabla u|\leq c_{1}\left(\left(\frac1{\gth\ge}\right)^{\frac1{q-1}}+1+\max_{\overline\Gw_\frac{\gth}{1+\ge}}
u^{\frac pq}\right).
$$
Next we estimate $u$ in function of its gradient: for any $x\in \overline\Gw_\frac{\gth}{1+\ge}$,
$$\dsps
u(x)\leq u\left((1-\gth)\frac x{|x|}\right)+\left|x-(1-\gth)\frac x{|x|}\right|\max_{y\in [x,(1-\gth)\frac x{|x|}]}|\nabla u(y)|.
$$
Therefore
$$\dsps\max_{\overline \Gw_{\frac\gth{1+\ge}}}u\leq \max_{\overline B_1\setminus B_{\frac 12}}u+\max_{\overline \Gw_{\frac\gth{1+\ge}}}|\nabla u|\leq c'_{1}+\max_{\overline \Gw_{\frac\gth{1+\ge}}}|\nabla u|.$$
Since $1\leq\frac1{\gth\ge}$, we deduce
$$\dsps
\max_{\overline\Gw_\gth}|\nabla u|\leq c_{2}\left(\left(\gth\ge\right)^{-\frac{1}{q-1}}+\left(\max_{\overline \Gw_{\frac\gth{1+\ge}}}|\nabla u|\right)^{\frac pq}\right).
$$
We set
$$A(\gth)=\gth^{\frac{1}{q-1}}\max_{\Gw_{\gth}}|\nabla u|,
$$
then $A(\frac{\gth}{1+\ge})\leq A((1-\frac\ge 2)\gth)$ since $\ge,\gth\leq \frac 12$, hence 
$$A(\gth)\leq c_{4}\left(\ge^{-\frac{1}{q-1}}+\gth^{\frac{q-p}{q(q-1)}}(1+\ge)^{\frac{p}{q(q-1)}}\left(A((1-\tfrac\ge 2)\gth)\right)^{\frac pq}\right).
$$
If we set $F(\gth)=1+A(\gth)$ there holds
\bel{Es9+}\BA{lll}
F(\gth)\leq c_{5}\ge^{-\frac{1}{q-1}}F^{\frac pq}(A(1-\tfrac\ge 2)\gth),
\EA\ee
and we can apply the bootstrap result of \rlemma{ST0} with $\Gf=1$, $h=\frac{1}{q-1}$ and $d=\frac pq$. We deduce that $F$ is bounded, hence 
\bel{Es10+}\BA{lll}\dsps
\max_{\overline\Gw_\gth}|\nabla u|\leq c_{6}\gth^{-\frac{1}{q-1}}.
\EA\ee
Thus $(\ref{Est-N1})$ holds.\\


\medskip

\nind 2- {\it Proof of 2-} By change of scale we assume again that ${r_0}=1$.
For  $T>3$ and $0<\ge<1/2$ we set 
$$\Gw_T=B_{T}\setminus \overline B_{1}\;\text{ and }\; \Gw_{T,\ge}=B_{T-\ge}\setminus \overline B_{1+\ge}.$$
 By $(\ref{Es13})$, for any $\gr>0$ and $x\in B^c_{1+2\gr}$ we have
$$\dsps
|\nabla u(x)|\leq c_{7}\left(\gr^{-\frac{1}{q-1}}+1+\max_{\overline B_{2\gr}(x)}u^{\tfrac pq}\right).
$$
Taking $\gr=\frac\ge 2$ we get 
\bel{Es11}\dsps
\max_{\overline \Gw_{T,\ge}}|\nabla u|\leq c_{8}\left(\ge^{-\frac{1}{q-1}}+1+\max_{\overline \Gw_{T}}u^{\frac pq}\right).
\ee
It is clear that
$$\dsps
\max_{\overline \Gw_{T}}u\leq \max_{|x|=1}u(x)+T\max_{\overline \Gw_{T}}|\nabla u|.
$$
reporting this inequality in $(\ref{Es11})$ we obtain that for any $T\geq 1$, 
\bel{Es12}\BA{lll}\dsps
1+\max_{\overline \Gw_{T,\ge}}|\nabla u|\leq c_{9}\ge^{-\frac 1{q-1}}T^{\frac pq}\left(1+\max_{\overline \Gw_{T}}|\nabla u|\right)^{\frac pq}.
\EA\ee
We set $F(T)=1+\max_{\overline \Gw_{T}}|\nabla u|$, then
\bel{Es13*}\BA{lll}\dsps
F(T(1-\ge))\leq 1+\max_{1\leq |x|\leq 1+\ge}|\nabla u(x)|+\max_{\overline \Gw_{T,\ge}}|\nabla u|\\[2mm]
\phantom{f(T(1-\ge))}
\leq 1+\max_{1\leq |x|\leq 2}|\nabla u(x)|++\max_{\overline \Gw_{T,\ge}}|\nabla u|\\[2mm]
\phantom{f(T(1-\ge))}
\leq c_{10}\left(\ge^{-\frac 1{q-1}}+1+\left(\max_{|x|=1}u(x)+T\max_{\overline \Gw_{T}}|\nabla u|\right)^{\frac pq}\right)\\[2mm]
\phantom{f(T(1-\ge))}
\leq c_{11}\ge^{-\frac{1}{q-1}}T^{\frac pq}F^{\frac pq}(T).
\EA\ee
Using again the bootstrap result of \rlemma{ST0} with $d=\frac pq$ we obtain in particular for $T\geq 2$,
\bel{Es14}
F(T)\leq c_{12}T^{\frac pq\frac{1}{1-\frac{p}{q}}}=c_{12}T^{\frac p{q-p}}.
\ee
This implies 
\bel{Es15}|\nabla u(x)|\leq c_{13}|x|^{\frac p{q-p}}.\ee
Using $(\ref{Es15})$ we get
$$\dsps
\max_{\overline\Gw_T}u\leq \max_{|x|=1}u(x)+T\max_{\overline\Gw_T}|\nabla u|\leq c_{14}T^{1+\frac{p}{q-p}}=c_{14}T^{\frac{q}{q-p}},
$$
which leads to 
\bel{Es16}
u(x)\leq c_{14}|x|^{\frac q{q-p}}\quad\text{for all }x\in B^c_3.
\ee
\qeda\medskip

By integrating the inequalities $(\ref{Est-N1})$ and $(\ref{Est-N2})$, we obtain:
\bcor{TC1} Under the assumption of \rth {Theorem E},  any nonnegative solution $u$ of $(\ref{I-1})$ in $G$ satisfies:\smallskip

\nind 1- If $G=B_{r_0}\setminus\{0\}$.\\ 
1-(i) If $q>\max\{2,p\}$, then $u$ can be extended as a continuous function in $B_{r_0}$.\\ 
1-(ii) If $q=2>p$, then there exists a constant $C_1>0$ such that 
\bel{Es26}\dsps
u(x)\leq C_1(|\ln  |x||+1)\quad\text{for all }x\in B_{\frac {r_0}2}\setminus\{0\}.
\ee
1-(iii) If $2>q>p$, then there exists a constant $C_2>0$ such that 
\bel{Es27}\dsps
u(x)\leq C_2|x|^{-\frac{2-q}{q-1}}\quad\text{for all }x\in B_{\frac {r_0}2}\setminus\{0\}.
\ee

\nind 2- If $G= B^c_{r_0}$, then there exists a constant $C_3>0$ such that 
\bel{Es28}\dsps
u(x)\leq C_3|x|^{\frac{q}{q-p}}\quad\text{for all }x\in B^c_{2{r_0}}\setminus\{0\}.
\ee
\es
\medskip

\nind \Remark The constants $C_i$ in $(\ref{Es26})$-$(\ref{Es27})$ (resp. $(\ref{Es28})$) depend on $\dsps\sup_{B_{r_0}\setminus B_{\frac {3r_0}4}} u(y)$ (resp. $\dsps\sup_{B_{2{r_0}}\setminus B_{{r_0}}}u(y)$). Up to modifying $\gth$ it is possible to reduce that domain of dependance of the constant with respect to $u$ 
 to $\dsps\sup_{B_{r_0}\setminus B_{(1-\gt)r_0}} u(y)$ (resp. $\dsps\sup_{B_{(1+\gt)r_0}\setminus B_{r_0}} u(y)$ for any $\gt\in (0,1)$.
\subsection{Upper estimates on solutions when $q < p$. Proof of \rth{Theorem E'}}
We recall the doubling Lemma \cite {Hu},\cite{PoQuSo}.
\bth{dlem} Let $(X,d)$ be a complete metric space, $D$ a non-empty subset of $X$, $\Gs$ a closed subset of $X$ containing $D$ and 
$\Gg=\Gs\setminus D$. Let $M:D\mapsto (0,\infty)$ be a map which is bounded on compact subsets of $D$ and let $k>0$ be a real number. If $y\in D$ is such that 
$$M(y)\dist(y,\Gg)>2k,
$$
there exists $x\in D$ such that 
$$\BA{lll}
&M(x)\dist(x,\Gg)>2k\\[1mm]
&M(x)\geq M(y)\\
&M(z)\leq 2M(x)\quad\text{for all }z\in D\text{ s.t. }d(z,x)\leq \myfrac{k}{M(x)}.
\EA$$
\es

\nind {\it Proof of \rth{Theorem E'}-(1).} We can assume that $r_0=1$. By $(\ref{Es9})$, $(\ref{AS1})$ implies that $|\nabla u(x)|\to 0$ when $|x|\to\infty$. The estimate $(\ref{b4})$ is equivalent to 
\bel{b4*}
u(x)\leq C|x|^{-\frac{q}{p-q}}= C |x|^{-\gg}
\ee
for all $x\in B^c_{2}$ by $(\ref{Es5})$, hence also to 
\bel{b5}
u^\frac{1}{\gg}x)+ |\nabla u(x)|^{\frac{1}{\gg+1}}\leq \frac C{|x|}
\ee
for all $x\in B^c_{2}$. We set 
\bel{b5}
M(x):=u^\frac{1}{\gg}(x).
\ee
Then $M(x)\to 0$ when $|x|\to\infty$. Let us assume that $|x|^\gg u(x)$ is unbounded in $B^c_{2r_0}$. Then by \rth{dlem} applied with $\Gs=B_2^c$, $D=\overline B_2^c$, thus $\Gg=B_2^c\setminus \overline B_2^c=\prt B_2$, and $k=n$, there exists a sequence  $\{y_n\}\subset\overline B_2^c$ such that $(|y_n|-2)M(y_n)\to\infty$ when $n\to\infty$. There exists a sequence $\{x_n\}\subset \overline B_2^c$ such that 
\bel{AS3*}\BA{lll}
&|x_n|M(x_n)>(|x_n|-2)M(x_n)>2n\\[1mm]
&M(x_n)\geq M(y_n)\\[1mm]
&M(z)\leq 2M(x_n)\quad\text{for all }z\in \overline B_2^c\text{ s.t. }|z-x_n|\leq \myfrac{n}{M(x_n)}.
\EA\ee
Clearly $\{x_n\}$ is unbounded since $M$ is bounded on bounded subsets of $B_2^c$ and, up to extracting a sequence, we can assume that 
$|x_n|\to\infty$ as $n\to\infty$. 
We now define
\bel{b6}u_n(x)=\frac{u(z(x,n))}{M^\gg(x_n)}\quad\text{with }z(x,n)=x_n+\frac{x}{M(x_n)}.
\ee
Then 
\bel{b6*}\dsps u_n(0)=1\quad\text{and }\,u_n(x)\leq 2^\gg\quad\text{for }x\in B_n.
\ee
The main point is to use estimate $(\ref{Es13})$ in order to obtain a uniform estimate on $\nabla u_n$. We apply this inequality in $B_{\frac{n}{M(x_n)}}(x_n)$ which yields
\bel{b6**}\dsps
\max_{z\in B_{\frac{n}{2M(x_n)}}(x_n)}|\nabla u(z)|\leq c_7\left(\left(\frac{n}{2M(x_n)}\right)^{-\frac 1{q-1}}+
\max_{z\in B_{\frac{n}{M(x_n)}}(x_n)}\left(u^{\frac{p}{q}}(z)+u^{\frac{p-1}{2(q-1)}}(z)\right)\right)
\ee
Furthermore
$z\in B_{\frac{n}{M(x_n)}}(x_n)$ is equivalent to $ |x|\leq n$. Similarly, $z\in B_{\frac{n}{2M(x_n)}}(x_n)$ is equivalent to $ |x|\leq\frac n2$.
If $u_n$ is defined by  $(\ref{b6})$, then 
 $$\nabla u_n(x)=\frac{\nabla u(z(x,n))}{M^{\gg+1}(x_n)}.$$
 We have that $\frac pq< \frac{p-1}{2(q-1)}$ since $q<\frac{2p}{p+1}$. Combined with the decay estimate $(\ref{b3})$ we infer that 
\bel{b6x}\max_{z\in B_{\frac{n}{M(x_n)}}(x_n)}\left(u^{\frac{p}{q}}(z)+u^{\frac{p-1}{2(q-1)}}(z)\right)\leq c_8\max_{z\in B_{\frac{n}{M(x_n)}}(x_n)}u^{\frac{p}{q}}(z).
\ee
We now replace $u(z)$ and $\nabla u(z)$ by their respective value with respect to $u_n(x)$ and $\nabla u_n(x)$ and we get
\bel{b7}\dsps
\max_{|x|\leq\frac n2}|\nabla u_n(x)|\leq c_9\left(n^{-\frac 1{q-1}} \left(M(x_n)\right)^{\frac{1}{q-1}-\gg-1}+
\max_{|x|\leq n}u^\frac pq_n(x)\right).
\ee
Because $1<q<\frac {2p}{p+1}$, $\frac{1}{q-1}-\gg-1>0$. Since $M(x_n)\to 0$ when $n\to\infty$ it follows that 
\bel{b7*}\dsps
|\nabla u_n(x)|\leq c_{10}\quad\text{for all }x\in B_{\frac n2}.
\ee
Therefore the new constraints are
\bel{b7**}\dsps u_n^\frac{1}{\gg}0)=1\quad\text{and }\,u_n(x)+|\nabla u_n(x)|\leq 2^\gg+c_{10}\quad\text{for }x\in B_{\frac n2}.
\ee
We have also
$$-\Gd u_n(x)=-\frac{\Gd u(z(x,n))}{M^{\gg+2}(x_n)},
$$
hence
$$\BA{lll}\dsps
-\Gd u_n(x)=\frac{u^p(z(x,n))-m|\nabla u(z(x,n))|}{M^{\gg+2}(x_n)}
\\[5mm]\phantom{-\Gd u_n(x)}\dsps
=\frac{M^{\gg p}(x_n)u_n^p(x)-mM^{(\gg +1)q}(x_n)|\nabla u_n(x)|}{M^{\gg+2}(x_n)}
\\[5mm]\phantom{-\Gd u_n(x)}
\dsps =M^{\gg( p-1)-2}(x_n)u_n^p-mM^{(\gg(q-1)-2+q)q}(x_n)|\nabla u_n(x)|^q.
\EA$$
There holds
$$\gg(p-1)-2=\gg(q-1)-2+q=\frac{\gs}{p-q},
$$
and by assumption, $\gs<0$. Therefore $u_n$ satisfies 
\bel{b16}
-\ge_n\Gd u_n(x)=u_n^p-m|\nabla u_n|^q\quad\text{with }\, \ge_n=M^{-\frac{\gs}{p-q}}(x_n)\to 0\,\,\text{as }n\to\infty.
\ee
Jointly with the conditions $(\ref{b7**})$ there exists a subsequence 
of  $\{u_n\}$ still denoted by $\{u_n\}$ and a function $v\in W^{1,\infty}(\BBR^N)$ such that $u_n$ converges to $v$ locally uniformly in $\BBR^N$
  and $\nabla u_n\rightharpoonup\nabla v$ for the weak topology of $L_{loc}^\infty(\BBR^N)$. By a classical viscosity result \cite[Proposition IV.1]{CL}, $v$ is a bounded viscosity solution of 
\bel{b17}m|\nabla v]^q=v^p\quad\text{in }\BBR^N.
\ee
By \cite[Proposition 4.3]{CL} $(\ref{b7})$ has a unique viscosity solution which is zero which is not compatible with $v(0)=1$ by $(\ref{b7**})$, which ends the proof.\medskip

\nind {\it Proof of \rth{Theorem E'}-(2).} We can take that $r_0=1$. The proof is still based upon \rth{dlem} with 
$\Gs=\overline B_{\frac12}$, $D=\overline B_{\frac12}\setminus\{0\}$ and $\Gg=\{0\}$. Thus we assume that there exists a solution $u\in C(\overline B_{1}\setminus\{0\})$,  solution of $(\ref{I-1})$ in $B_{1}\setminus\{0\}$
and a sequence of points $\{y_n\}\subset \overline B_{1}\setminus\{0\}$ such that 
\bel{xx0}
|y_n|M(y_n)\geq 2n
\ee
where we have set
$$M(x)=u^\frac 1\gg(x).
$$
There exists a sequence $\{x_n\}\subset B_{1}\setminus\{0\}$ such that 
\bel{xx1}\BA{lll}&|x_n|M(x_n)>2n\\[1mm]
 &M(x_n)\geq M(y_n)\\[1mm]
 & M(z)\leq 2M_n(x_n)\quad\text{for all }z\in B_{\frac{n}{M(x_n)}}(x_n).
 \EA\ee
 Clearly $x_n\to 0$ as $n\to\infty$. We define $u_n$ by $(\ref{b6})$ and $(\ref{b6*})$ holds. The gradient estimate $(\ref{b6**})$ is verified and 
 if $z\in B_{\frac{n}{M(x_n)}}(x_n)$ , we have $|z|\leq |x_n|+|z-x_n|\leq |x_n|+\frac{n}{M(x_n)}$ which tends to $0$ as $n\to\infty$. If we replace $u(z)$ by $u_n(x)=\frac{u(z(x,n))}{M^\gg(x_n)}$, $(\ref{b6**})$ becomes
 \bel{xx2}\dsps
\max_{|x|\leq\frac n2}|\nabla u_n(x)|\leq c_{11}\left(n^{-\frac 1{q-1}} \left(M(x_n)\right)^{\frac{1}{q-1}-\gg-1}+\max_{|x|\leq n}\left(u^{\frac pq}_n(x)+\left(M(x_n)\right)^{-\frac{\gs}{2(q-1)(p-q)}}u^{\frac {p-1}{2(q-1)}}_n(x)\right)
\right).
\ee
Notice that $M(x_n)\to\infty$ and $\frac{1}{q-1}-\gg-1=\frac{-\gs}{(q-1)(p-q)}<0$. Using 
 $(\ref{b6*})$ we obtain
  \bel{xx3}\dsps
\max_{|x|\leq\frac n2}|\nabla u_n(x)|\leq c_{11}\left(o(1)+2^{\frac{p}{p-q}}+o(1)
\right)\leq c_{12}.
\ee
 Hence $(\ref{b7**})$ holds with a new constant $c_{13}$. Equation $(\ref{b16})$ is verified, but now $\gs>0$. Hence $\ge_n\to 0$ as $n\to\infty$.
We conclude by the same argument as the one used in (1).
\qeda\medskip

\nind \Remark In \rth{Theorem E'}-(2) It is possible to obtain a constant $C$ in estimate $(\ref{b4})$ independent $u$ provided the functions under consideration
are uniformly locally bounded from above in $\overline B_{r_0}\setminus\{0\}$ in the sense that for any $\ge>0$ there exists $C_\ge>0$ independent of $u$ such that
  \bel{xx'}\dsps
u(x)\leq C_{\ge}\quad\text{for all }x\in B_{r_0}\setminus B_\ge.
\ee
This assumption implies  that  in the proof of  \rth{Theorem E'}-2), $M(x_n)\to \infty$ independently of $u$.

\subsection{Asymptotic estimates on decaying solutions in the case $q>\frac {2p}{p+1}$}
Using \rth{dlem}, we prove \rth {Theorem F}.\smallskip

\nind {\it Proof of \rth {Theorem F}.} We can assume that $r_0=1$. By $(\ref{Es9})$, $\nabla u(x)$ tends to $0$ as $|x|\to\infty$. Estimate $(\ref{AS2})$ is equivalent to 
\bel{AS3}
M(x):=u^\frac{p-1}{2}(x)+ |\nabla u(x)|^\frac{p-1}{p+1}\leq C|x|^{-1}\quad\text{for all }x\in B^c_{2}.
\ee
Using $(\ref{AS1})$ jointly with $(\ref{Es9})$ we have that $M(x)\to 0$ as $|x|\to\infty$. Let us assume that for any $C>0$ inequality $(\ref{AS3})$ does not hold; then there exists a sequence $\{y_n\}\subset B_{2}^c$ such that $\lim_{n\to\infty}(|y_n|-2)M(y_n)=\infty$. There exists a sequence $\{x_n\}\subset \overline B_2^c$ such that $\ref{AS3*}$ holds.
Clearly $\{x_n\}$ is unbounded since $M$ is bounded on bounded subset of $B_2^c$ and, up to extracting a sequence, we can assume that 
$|x_n|\to\infty$ as $n\to\infty$. We set
\bel{AS3**}u_n(x)=\frac{u(z(x,n)}{M^\ga(x_n)}\quad\text{with }\,z(x,n)=x_n+\frac{x}{M(x_n)}. 
\ee
Then we have $M(x_n)|x_n|>2n$ and for any $x\in B_n$,
\bel{AS4}
M(z(n,x))=u^{\frac{p-1}{2}}(z(n,x)+|\nabla u|^{\frac{p-1}{p+1}}(z(n,x)\leq 2M(x_n).
\ee
Then
$$\nabla u_n(x)=\frac{\nabla u(z(x,n))}{M^{\ga+1}(x_n)}\,,\; \Gd u_n(x)=\frac{\Gd u(z(x,n))}{M^{\ga+2}(x_n)},
$$
which implies
$$\BA{lll}\dsps
\Gd u_n(x)=\frac{u^p(z(x,n))-m|\nabla u|^q(z(x,n))}{M^{\ga+2}(x_n)}
\\[4mm]\phantom{\Gd u_n(x)}
\dsps =\frac{M^{\ga+2}(x_n)u_n(x)-mM^{(\ga+1)q}(x_n)|\nabla u(z(x,n))|^q}{M^{\ga+2}(x_n)}.
\EA$$
Hence $u_n$ satisfies 
$$-\Gd u_n=u_n^p-m(M(x_n))^{(\ga+1)q-\ga p}|\nabla u_n|^q\quad\text{in }B_n,
$$
with the additional condition
$$u_n^{\frac{p-1}{2}}(0)+|\nabla u_n(0)|^\frac{p-1}{p+1}=1.
$$
Observe that 
 $$ (\ga+1)q-\ga p=\frac{(p+1)q-2p}{p-1}\geq 0,$$ 
with equality if $q=\frac{2p}{p+1}$ and strict inequality otherwise. 
Furthermore
$$u_n^{\frac{p-1}{2}}(x)+|\nabla u_n(x)|^\frac{p-1}{p+1}\leq 2\quad\text{for all }x\in B_n.
$$
By standard elliptic equations regularity results \cite{GT}, the sequence $\{u_n\}$ is eventually locally compact in the $C^1_{loc}(\BBR^N)$-topology, thus, up to extracting a subsequence, $\{u_n\}$ converges in this topology to some nonnegative $C^1(\BBR^N)$ function $v$ which satisfies 
 \bel{b1}-\Gd v=v^p\quad\text{in }\BBR^N
\ee
 if $q>\frac{2p}{p+1}$ since  $M(x_n)\to 0$ as $n\to\infty$,  and 
 \bel{b2}-\Gd v+m|\nabla v|^q=v^p\quad\text{in }\BBR^N
 \ee
 if $q=\frac{2p}{p+1}$. Furthermore $v^{\frac{p-1}{2}}(0)+|\nabla v(0)|^\frac{p-1}{p+1}=1$. Since $1<p<\frac{N+2}{N-2}$, by Gidas and Spruck
  result \cite {GS} equation $(\ref{b1})$ admits no global positive solution. Concerning $(\ref{b2})$, if $m\leq\ge_0$ satisfies  no global positive solution can exist by Theorem B.
  This ends the proof.\qeda\medskip
  
  \nind\Remark In the case $q=\frac{2p}{p+1}$, the assumption $(\ref{AS1})$ can be relaxed and replaced by
  \bel{AS1*}
  \dsps\limsup_{|x|\to\infty}u(x)<\infty. 
  \ee
  Actually, if this holds we have by $(\ref{Es8})$
    \bel{AS1**}
  \dsps\limsup_{|x|\to\infty}|\nabla u(x)|<\infty. 
  \ee
  The function $u_n$ defined by $(\ref{AS3**})$ satisfies the same equation $(\ref{I-1})$ as $u$ and the limit $v$ also. We end the proof as in \rth {Theorem F}.\medskip

.
\mysection{Removable singularities}
In this Section we give partial extensions to $(\ref{I-1})$ of previous results dealing with removability of singularities for equations
$$-\Gd u+m|\nabla u|^q= 0
$$
and 
$$-\Gd u+m|\nabla u|^2-u^p\leq 0,
$$
obtained respectively in \cite{NgV} and \cite{BrNi}.
\subsection{Removable isolated singularities. Proof of \rth{Theorem G}} 
\nind{\it Proof of \rth {Theorem G}.} We can assume that $\overline B_{r_0}\subset\Gw$ with ${r_0}\geq 1$ and $a=0$. Since $(\ref{Est-N1})$ holds we have 
\bel{Rem1}|\nabla u(x)|\leq c|x|^{-\frac1{q-1}}\quad\text{and }\;u(x)\leq c_1+c_2\left\{\BA{lll}|x|^{\frac{q-2}{q-1}} &\text{ if }q>2\\[2mm]
|\ln|x||&\text{ if }q=2\EA\right.\quad\text{ for  }0<|x|\leq {r_0}.
\ee
Since $q>p$ and $q\geq \frac{N}{N-1}$, we have that $\nabla u\in L^p(B_{r_0})$, which implies $u^p\in L^1(B_{r_0})$. \smallskip

\nind {\it Step 1: We claim that $\nabla u\in L^q(B_{r_0})$ and the equation holds in $\CD'(B_{r_0})$}. Let $\eta_n\in C^{\infty}_0(B_{r_0}\setminus\{0\})$ such that $\eta_n=1$ on $B_{{r_0}/2}\setminus B_{1/n}$, $\eta_n=0$ if 
$|x|\leq 1/{2n}$ and if $|x|\geq 2{r_0}/3$ and $0\leq \eta_n\leq 1$. We construct $\eta_n$ such that $|\nabla\eta_n|\leq cn{\bf 1}_{B_{1/n}\setminus B_{1/2n}}$. Then 
$$\int_{B_{r_0}}\nabla u.\nabla\eta_n dx+m\int_{B_{r_0}}|\nabla u|^q\eta_n dx=\int_{B_{r_0}}u^p\eta_n dx.
$$
By Holder's inequality and using $(\ref{Est-N1})$ there holds with $q'=\frac q{q-1}$,
$$\left|\int_{B_{r_0}}\nabla u.\nabla\eta_n dx\right|=\left|\int_{B_{1/n}\setminus B_{1/2n}}\nabla u.\nabla\eta_n dx\right|\leq c_2n^{q'-N}.
$$
Since $q\geq \frac{N}{N-1}$, then $q'-N\leq 0$, and the right-hand side is bounded, hence $|\nabla u|^q\in L^1(B_{\frac {r_0}2})$ by Fatou's theorem and the first statement follows. \\
Next consider $\gz\in C^{\infty}_0(B_{{r_0}/2})$ and take $\gz\eta_n$ as a test function, then
$$\int_{B_{r_0}}\left(\gz\nabla u.\nabla\eta_n+ \eta_n\nabla u.\nabla\gz\right)dx+m\int_{B_{r_0}}|\nabla u|^q\gz\eta_n dx=\int_{B_{r_0}}u^p\gz \eta_n dx.
$$
Since 
\bel{+1}\left|\int_{B_{r_0}}\gz\nabla u.\nabla\eta_n dx\right|\leq c_3n^{1-\frac N{q'}}\norm{\gz}_{L^\infty}\left(\int_{B_{1/n}\setminus B_{1/2n}}|\nabla u|^q\right)^{\frac{1}{q}},
\ee
and the left-hand side tends to $0$ as $n\to\infty$, we conclude by the dominated convergence theorem that 
$$\int_{B_{r_0}}\nabla u.\nabla\gz dx+m\int_{B_{r_0}}|\nabla u|^q\gz dx=\int_{B_R}u^p\gz dx,
$$
which proves the second statement.  \smallskip

\nind {\it Step 2: $u$ is bounded}. 
For proving the boundedness assertion we can assume that $\frac N{N-1}\leq q<2$.
 As a test function we take $\gz=\eta^q_n$, then 
$$q\int_{B_{r_0}}\eta_n^{q-1}\nabla u.\nabla\eta_n dx+m\int_{B_{r_0}}\eta_n^{q}|\nabla u|^q dx=\int_{B_{r_0}}\eta_n^{q}u^p dx.
$$
We have 
$$\BA{lll}
\dsps\int_{B_{r_0}}\eta_n^{q}|\nabla u|^q dx=\int_{B_{r_0}}|\eta_n\nabla u|^q dx=\int_{B_{r_0}}|\nabla (\eta_nu)-u\nabla \eta_n|^q dx\\
\phantom{\dsps\int_{B_{r_0}}\eta_n^{q}|\nabla u|^q dx}\dsps
\geq 2^{1-q}\int_{B_{r_0}}|\nabla (\eta_nu)|^q dx-\int_{B_{r_0}}u^q|\nabla \eta_n|^qdx.
\EA$$
By $(\ref{Rem1})$
$$\int_{B_{r_0}}u^q|\nabla \eta_n|^qdx\leq c_4n^{q'-N}\leq c'
$$
 as we have already seen it and, from $(\ref{+1})$ there holds
$$\left|\int_{B_{r_0}}\eta_n^{q-1}\nabla u.\nabla\eta_n dx\right|\to 0\text{ as }n\to\infty.
$$
It follows that $\nabla (\eta_nu)$ is bounded in $L^q(B_{r_0})$ independently of $n$, and by Sobolev inequality, 
$$\norm{\eta_n u}_{L^{q^*}(B_{r_0})}\leq c''\quad\text{with }q^*=\frac{Nq}{N-q},$$
which in turn implies that $\norm{ u}_{L^{q^*}(B_{r_0})}\leq c_1$. Set
\bel {Rem2}r_1=\frac{Nq}{N-q}-p.
\ee
Taking $\eta^{q+r_1}_n(T_k(u))^{r_1}$ as a test function, where $T_k(r)=\min\{r,k\}$ for $r,k>0$, we obtain
$$\BA{lll}
\dsps
{r_1}\int_{B_{r_0}\cap\{u<k\}}(T_k(u))^{r_1-1}\eta_n^{q+r_1}|\nabla u|^2 dx+(q+r_1)\int_{B_{r_0}}(T_k(u))^{r_1}\eta_n^{q+r_1-1}\nabla\eta_n.\nabla u dx\\[4mm]\phantom{-----------}
\dsps +m\int_{B_{r_0}} T_k(u^{r_1})|\nabla u|^q\eta^{q+r_1}_n dx=\int_{B_{r_0}}T_k(u^{r_1})u^p\eta^{q+r_1}_n dx.
\EA$$
From {\it Step 1} $|\nabla u|\in L^q(B_{r_0})$, thus 
$$\int_{B_{r_0}}(T_k(u))^{r_1}\eta_n^{q+r_1-1}\nabla\eta_n.\nabla u dx\to 0\quad\text{as }n\to\infty,
$$
hence
$$o(1)+m\int_{B_{r_0}} T_k(u^{r_1})|\nabla u|^q\eta^{q+r_1}_n dx\leq \int_{B_{r_0}}T_k(u^{r_1})u^p\eta^{q+r_1}_n dx.
$$
Letting successively $n\to\infty$ and $k\to\infty$, we deduce by Fatou's lemma and the monotone convergence theorem that 
\bel {Rem3}
m\int_{B_{r_0}} u^{r_1}|\nabla u|^q \tilde\eta^{q+r_1} dx\leq \int_{B_{r_0}}u^{\frac {Nq}{N-q}} \tilde\eta^{q+r_1}dx,
\ee
where $\dsps \tilde\eta^{q+r_1}=\lim_{n\to\infty}\eta^{q+r_1}_n$ belongs to $C^\infty_0(B_{r_0})$ and takes value $1$ in $B_{\frac {r_0}2}$ and $0\leq\tilde\eta\leq 1$.
Since
$$\BA{lll}
\dsps\int_{B_{r_0}} u^{r_1}|\nabla u|^q \tilde\eta^{q+r_1} dx=\left(\frac{q}{q+r_1}\right)^q\int_{B_{r_0}} |\tilde\eta^{1+\frac{r_1}{q}}\nabla (u^{1+\frac{r_1}{q}})|^q dx\\[4mm]
\phantom{\dsps\int_{B_{r_0}} u^{r_1}|\nabla u|^q \tilde\eta^qdx}
\dsps
\geq \left(\frac{q}{r_1+q}\right)^q2^{1-q}\int_{B_{r_0}} |\nabla (\tilde\eta u)^{1+\frac {r_1}q}|^q dx
-\left(\frac{q}{r_1+q}\right)^q\int_{B_{r_0}}u^{q+r_1} |\nabla \tilde\eta|^q dx
\\[4mm]
\phantom{\dsps\int_{B_{r_0}} u^{r_1}|\nabla u|^q\eta^{_*\,q}  dx}\dsps
\geq c_{N,q}\left(\frac{q}{r_1+q}\right)^q\left(\int_{B_{r_0}} (\tilde\eta u)^{\frac {N(q+r_1)}{N-q}} dx\right)^{\frac{N-q}{N}}-K_1\left(\frac{q}{r_1+q}\right)^q,
\EA$$
where
$$K_1={r_0}^N\norm u^{q+r_1}_{L^\infty(B_{r_0}\setminus B_{\frac {r_0}2})}\norm {\nabla \tilde\eta}^q_{L^\infty(B_{r_0})}.
$$
This leads to the following inequality
\bel {Rem4}\BA{lll}\dsps
mc_{N,q}\left(\frac{q}{r_1+q}\right)^q\norm{\tilde\eta u}^{q+r_1}_{L^{\frac{N(q+r_1)}{N-q}}(B_{r_0})}-mK_1\left(\frac{q}{r_1+q}\right)^q\leq \norm{\tilde\eta^\frac {(N-q)(q+r_1)}{Nq}u}^{\frac {Nq}{N-q}}_{L^{\frac {Nq}{N-q}}(B_{r_0})}\\[4mm]
\phantom{\dsps
mc_{N,q}\left(\frac{q}{r_1+q}\right)^q\norm{\tilde\eta u}^{q+r_1}_{L^{\frac{N(q+r_1)}{N-q}}(B_{r_0})}-mK_1\left(\frac{q}{r_1+q}\right)^q}
\dsps\leq \norm{\tilde\eta u}^{\frac {Nq}{N-q}}_{L^{\frac {Nq}{N-q}}(B_{r_0})},
\EA\ee
since $\frac {(N-q)(q+r_1)}{Nq}>1$ from $(\ref{Rem2})$ and $q>p$ combined with the fact that $\tilde\eta\leq 1$.

Next we proceed by induction, setting
\bel {Rem4*}r_{j+1}=\frac{N(q+r_j)}{N-q}-p\quad\text{for }j\geq 1,
\ee
with explicit value 
\bel {Rem4**}
r_{j+1}=\left(\left(\frac{N}{N-q}\right)^{j+1}-1\right)\frac{(N-q)r_1}{q}.
\ee
Taking $\eta^{q+r_{j+1}}_nT_k(u^{r_{j+1}})$ for test function and letting successively $n\to\infty$ and $k\to\infty$ we obtain
\bel {Rem5}
m\int_{B_{r_0}} u^{r_{j+1}}|\nabla u|^q \tilde\eta^{q+r_{j+1}} dx\leq \int_{B_{r_0}}u^{\frac {N(q+r_j)}{N-q}} \tilde\eta^{q+r_{j+1}} dx\leq \int_{B_{r_0}}(\tilde\eta u)^{\frac {N(q+r_j)}{N-q}}dx.
\ee
Note that for the right-hand side we have used $q+r_{j+1}\geq \frac {N(q+r_j)}{N-q}$ and $\tilde\eta\leq 1$.  Moreover
\bel{+2}\BA{lll}\dsps\int_{B_{r_0}} u^{r_{j+1}}|\nabla u|^q \tilde\eta^{q+r_{j+1}}dx\geq \left(\frac{q}{r_{j+1}+q}\right)^q\int_{B_{r_0}} | \tilde\eta^{1+\frac{r_{j+1}}{q}}\nabla (u^{1+\frac{r_{j+1}}{q}})|^qdx.\EA\ee
Writing
$$\tilde\eta^{1+\frac{r_{j+1}}{q}}\nabla (u^{1+\frac{r_{j+1}}{q}})=\nabla (\tilde\eta u)^{1+\frac{r_{j+1}}{q}}-\frac{q+r_{j+1}}{q}u^{1+\frac{r_{j+1}}{q}}\tilde\eta^{\frac{r_{j+1}}{q}}\nabla\tilde\eta,
$$
we have, since $\tilde\eta=1$ in $B_{\frac {r_0}2}$ and $0\leq\tilde\eta\leq 1$, and using Sobolev inequality,
\bel {Rem6}\BA{lll}\dsps \norm{\tilde\eta^{1+\frac{r_{j+1}}{q}}\nabla (u^{1+\frac{r_{j+1}}{q}})}_{L^{q}(B_{r_0})}
 \geq \norm{\nabla (\tilde\eta u)^{1+\frac{r_{j+1}}{q}}}_{L^{q}(B_{r_0})}\\[4mm]\phantom{----------------} \dsps
 -\frac{q+r_{j+1}}{q}\norm{\nabla\tilde\eta}_{L^\infty}\norm{u^{1+\frac{r_{j+1}}{q}}}_{L^{q}(B_{r_0}\setminus B_{\frac {r_0}2})}\\[4mm]
 \phantom{---------}\dsps\geq c_{N,q}\norm{\tilde\eta u}^{\frac{q+r_{j+1}}{q}}_{L^{\frac{N(q+r_{j+1})}{N-q}}(B_{r_0})}-\frac{q+r_{j+1}}{q}\norm{\nabla\tilde\eta}_{L^\infty}\norm{u}^{\frac{q+r_{j+1}}{q}}_{L^{q+r_{j+1}}(B_{r_0}\setminus B_{\frac {r_0}2})}.
\EA \ee
Let us assume now that $u\notin L^\infty(B_{r_0})$, otherwise the result follows, then 
\bel {Rem7}
\lim_{j\to\infty}\norm{\tilde\eta u}_{L^{\frac{N(q+r_{j+1})}{N-q}}(B_{r_0})}=\infty,
\ee
and there exists $j_0\geq 1$ such that for any $j\geq j_0$,
\bel {Rem7*}
\norm{\tilde\eta u}_{L^{\frac{N(q+r_{j+1})}{N-q}}(B_{r_0})}\geq 2\norm{\nabla\tilde\eta}^{\frac{q}{q+r_{j+1}}}_{L^\infty}\norm{u}_{L^{q+r_{j+1}}(B_{r_0}\setminus B_{\frac {r_0}2})};
\ee
as a consequence the right-hand side of $(\ref{Rem6})$ is bounded from below by
\bel {Rem8}
\left(c_q-2^{-\frac{q+r_{j+1}}{q}}\frac{q+r_{j+1}}{q}\right)\norm{\tilde\eta u}^{\frac{q+r_{j+1}}{q}}_{L^{\frac{N(q+r_{j+1})}{N-q}}(B_{r_0})}
\geq \frac{c_{N,q}}{2}\norm{\tilde\eta u}^{\frac{q+r_{j+1}}{q}}_{L^{\frac{N(q+r_{j+1})}{N-q}}(B_{r_0})}
\ee
for $j\geq j_1\geq j_0$. Combining $(\ref{Rem5})$, $(\ref{+2})$ and $(\ref{Rem8})$ we derive
\bel {Rem9}
\frac 1m\int_{B_{r_0}}(\tilde\eta u)^{\frac{N(q+r_j)}{N-q}}dx\geq \left(\frac{qc_{N,q}}{2(r_{j+1}+q}\right)^q\norm{\tilde\eta u}^{q+r_{j+1}}_{L^{\frac{N(q+r_{j+1})}{N-q}}(B_{r_0})}.
\ee
We obtain finally 
\bel {Rem10}
\norm{\tilde\eta u}_{L^{\frac{N(q+r_{j+1})}{N-q}}(B_{r_0})}\leq\left(\frac{2(r_{j+1}+q)}{qc_{N,q}m^{\frac 1q}}\right)^{\frac{q}{q+r_{j+1}}}
\norm{\tilde\eta u}^{\frac{N(q+r_j)}{(N-q)(q+r_{j+1})}}_{L^{\frac{N(q+r_{j})}{N-q}}(B_{r_0})}.
\ee
Put
$$X_j=\ln\left(\norm{\tilde\eta u}_{L^{\frac{N(q+r_{j})}{N-q}}(B_{r_0})}\right).
$$
Since
\bel {Rem11}
\frac{N(q+r_j)}{(N-q)(q+r_{j+1})}=\frac{p+r_{j+1}}{q+r_{j+1}}<1,
\ee
we deduce 
\bel {Rem12}
X_{j+1}\leq \frac{q}{q+r_{j+1}}\ln\left(\frac{2(r_{j+1}+q)}{qc_qm^{\frac 1q}}\right)+X_j,
\ee
which implies that 
\bel {Rem13}\dsps
\ln\left(\norm u_{L^\infty(B_{\frac {r_0}2})}\right)\leq \limsup_{j\to\infty}X_{j+1}\leq X_1+q\sum_{j=1}^\infty\frac{1}{q+r_{j+1}}\ln\left(\frac{2(r_{j+1}+q)}{qc_qm^{\frac 1q}}\right)<\infty,
\ee
by $(\ref{Rem4**})$. This is a contradiction with $(\ref{Rem7})$, which ends the proof.
\qeda
\subsection{Removable singular sets }

In the following theorem we combine the technique of \rth {Theorem G} with the geometric approach based upon the construction of tubular neighbourhoods used in \cite{Ve1} to prove the removability of singular sets contained into a smooth 
submanifold. The next result proves and completes \rth{Theorem H}.

\bth{remov-2}Let $\Gw\subset\BBR^N$ be a bounded smooth domain with $N\geq 3$ and $\Gs\subset \Gw$ be a $k$-dimensional compact complete smooth submanifold with $1\leq k\leq N-2$. If $1\leq p< q$ and $q\geq \frac{N-k}{N-1-k}$, any nonnegative solution $u\in C^2(\Gw\setminus\Gs)$ of  $(\ref{I-1})$ in $\Gw\setminus\Gs$ can be extended as a weak solution of the same equation in $\Gw$ which belongs to $L^\infty_{loc}(\Gw)\cap W^{1,q}_{loc}(\Gw)\cap H^1_{loc}(\Gw)$.
\es
\Proof {\it Step 1: We claim that there exists $r_0>0$ and $C=C(N,p,q,m,r_0,\Gs)>0$ such that 
\bel{R1}
|\nabla u(x)|\leq C(\dist (x,\Gs))^{-\frac 1{q-1}}\quad\text{for all $x$ s.t. }\dist (x,\Gs)\leq r_0.
\ee}
For $\gd>0$ we set 
$$TUB_{\gd}(\Gs)=\{x\in\BBR^N:\dist (x,\Gs)<\gd\}.
$$
If $\gd\leq \inf\{dist (x,\Gs):x\in\Gw^c\}$, we have that $TUB_{\gd}(\Gs)\subset\Gw$. Since $\Gs$ is smooth with no boundary, there exixts $\gd_0>0$ such that the sets 
$\prt TUB_{\gd}(\Gs)=\{x\in\Gw:\dist (x,\Gs)=\gd\}$ are $k$-dimensional compact complete smooth submanifolds of $\Gw$. We use the ideas of the proof of \rth {Theorem E} adapting it to the peculiar geometric configuration. By rescaling we can assume that 
$\gd_0=1$ and for $0<\gth<\tfrac 14$, we set
$\Gth_\gth=TUB_{1-\gth}(\Gs)\setminus TUB_{\gth}(\Gs)$. For any $0<\ge<\frac 12$ we have by $(\ref{Es9*})$,
\bel{R2}\BA{lll}\dsps
\max_{\overline\Gth_\gth}|\nabla u|\leq c_{1}\left((\ge\gth)^{-\frac{1}{q-1}}+\max_{\overline\Gth_{\frac{\gth}{1+\ge}}}\left(u^p+u^{\frac{p-1}{2(q-1)}}\right)^{\frac{1}{q}}\right)\leq c_{2}\left((\ge\gth)^{-\frac{1}{q-1}}+1+\max_{\overline\Gth_{\frac{\gth}{1+\ge}}}u^{\frac pq}\right).\EA\ee
In order to obtain an upper bound on $u(x)$ for $x\in \overline\Gth_{\frac{\gth}{1+\ge}}$, we join it to some $x_\ge\in \prt TUB_{1}(\Gs)$ by a smooth curve $\gw$ such that $\gw(0)=x$, $\gw(1)=x_\ge$. We can choose $\gw$ such that $|\gw'(t)|\leq 2$ for all $t\in [0,1]$ and 
$$ 2^{-1}\dist(tx+(1-t)x_\ge,\Gs)\leq \dist(\gw(t),\Gs)\leq 2\dist(tx+(1-t)x_\ge,\Gs).$$
Then
\bel{R3}\BA{lll}\dsps u(x)\leq u(x_\ge)+\left|\int_0^1\nabla u(\gw(t)).\gw'(t) dt\right|\leq u(x_\ge)+2\int_0^1|\nabla u(\gw(t))|dt\\[2mm]
\phantom{u(x)}\dsps \leq \norm{u}_{L^\infty(TUB_1(\Gs)\setminus TUB_\frac12(\Gs))}+2\max_{\Gw_{\frac{\gth}{1+\ge}}}|\nabla u|.\EA\ee
Therefore
\bel{R4}\BA{lll}\dsps
\max_{\overline\Gth_{\frac{\gth}{1+\ge}}}u^{\frac pq}\leq c_{3}\left(\norm{u}^\frac pq_{L^\infty(TUB_1(\Gs)\setminus TUB_\frac12(\Gs))}+\max_{\overline\Gth_{\frac{\gth}{1+\ge}}}|\nabla u|^{\frac pq}\right)\\[6mm]
\phantom{\dsps\max_{\overline\Gth_{\frac{\gth}{1+\ge}}}u^{\frac pq}}\dsps\leq c_{3}\left(\norm{u}^\frac pq_{L^\infty(TUB_1(\Gs)\setminus TUB_\frac12(\Gs))}+\max_{\overline\Gth_{(1-\ge)\gth}}|\nabla u|^{\frac pq}\right).
\EA\ee
We put 
$$\dsps B(\gth)=\max_{\overline\Gth_\gth}\gth^{\frac{1}{q-1}}|\nabla u(z)|\,\text{ and }\; F(\gth)=1+B(\gth),
$$
and we obtain from $(\ref{R2})$ and $(\ref{R4})$
\bel{R5}
F(\gth)\leq c_{4}\ge^{-\frac 1{q-1}}F^{\frac pq}((1-\ge)\gth),
\ee
where $c_{4}$ depends on the structural constants and of $\norm{u}_{L^\infty(TUB_1(\Gs)\setminus TUB_\frac12(\Gs))}$.
It follows from \rlemma {ST0} that $B(\gth)$ is bounded independently of $\gth$, which implies $(\ref{R1})$. \smallskip

\nind In order to derive the upper estimate on $u$ we set $\gm= \sup\{u(y):y\in \prt TUB_1(\Gs)\}$. If $0<\dist(x,\Gs)=t\leq 1$ there exists 
 $z_x\in\Gs$ and $\xi \in \prt TUB_1(\Gs)$ such that 
$$2^{-1}|tx+(1-t)\xi-z_x|\leq \dist(tx+(1-t)\xi,\Gs)\leq 2|tx+(1-t)\xi-z_x|.
$$
Since $\dist(\xi,\Gs)=1$, 
$$\BA{lll}\dsps u(x)\leq \gm+c_{5}\int_0^1|tx+(1-t)\xi-z_x|^{-\frac 1{q-1}} dt\\[4mm]
\phantom{ u(x)}\dsps
\leq \gm+c_{5}\int_0^1\left(t\dist(x,\Gs)+(1-t)\dist(\xi,\Gs)\right)^{-\frac 1{q-1}}=\gm+c_{5}\int_0^1\left(t\dist(x,\Gs)+1-t\right)^{-\frac 1{q-1}}\\[4mm]
\phantom{ u(x)}\dsps
\leq \gm+c_{5}\frac{q-1}{2-q}\left(1-\dist(x,\Gs)\right)\left((\dist(x,\Gs))^{\frac{2-q}{q-1}}-1\right),
\EA$$
if $q\neq 2$, with an obvious modification if $q=2$. At end we deduce
\bel{R4*}\dsps
u(x)\leq c_{6}\left\{\BA{lll}(\dist(x,\Gs))^{\frac{2-q}{q-1}}+C'\quad&\text{for all }x\in TUB_1(\Gs)&\text{ if }q\neq 2\\[2mm]
|\ln(\dist(x,\Gs))|+C'\quad&\text{for all }x\in TUB_1(\Gs)&\text{ if }q= 2.
\EA\right.
\ee
\smallskip

\nind {\it Step 2: We claim that  $u\in L^p(TUB_{1}(\Gs))$ and $|\nabla u|\in L^q(TUB_{1}(\Gs))$}. For such a task we consider test functions 
$\eta_n\in C^\infty_0(TUB_1(\Gs))$ with value in $[0,1]$ vanishing in $TUB_{1/(2n)}(\Gs)\cup TUB^c_{2/3}(\Gs)$, with value $1$ in  $TUB_{1/2}(\Gs)\setminus TUB_{1/n}(\Gs)$
and such that $$|\nabla\eta_n(x)|\leq c_{7}n{\bf 1}_{TUB_{1/n}(\Gs)\setminus TUB_{1/2n}(\Gs)},$$
where the constant $c_7>0$ depends on the geometry of $\Gs$. If $q>2$, $u$ is bounded thus $u^p\in L^1(TUB_{1}(\Gs))$. If $\frac{N-k}{N-k-1}\leq q\leq 2$ we have for $1>\ge>\frac 1n$ 
$$\BA{lll}\dsps
\int_{TUB_{\ge}(\Gs)}\eta_nu^pdx\leq \int_{TUB_{\ge}(\Gs)\setminus TUB_{1/2n}(\Gs)}u^pdx\\[4mm]
\phantom{\dsps\int_{TUB_{\ge}(\Gs)}\eta_nu^pdx}\dsps
\leq c_8\int_{1/2n}^\ge\gt^{-\frac{(2-q)p}{q-1}}\frac{d}{d\gt}Vol(TUB_\gt(\Gs))d\gt
\\[4mm]
\phantom{\dsps\int_{TUB_{\ge}(\Gs)}\eta_nu^pdx}\dsps
\leq c_8\ge^{-\frac{(2-q)p}{q-1}}Vol(TUB_\ge(\Gs))
+c_8\frac{(2-q)p}{q-1}\int_{1/2n}^\ge\gt^{-\frac{(2-q)p}{q-1}-1}Vol(TUB_\gt(\Gs))d\gt.
\EA$$
By Weyl's formula \cite{We}
\bel{R5W}\dsps
Vol(TUB_\gt(\Gs)) =\sum_{i=0}^{[k/2]}a_i \gt^{N-k+2i}
\ee
where the $a_i$ are smooth bounded functions near $\Gs$ and $[k/2]$ is the integer part of $k/2$. Therefore 
$$\int_{1/(2n)}^\ge\gt^{-\frac{(2-q)p}{q-1}}\frac{d}{d\gt}Vol(TUB_\gt(\Gs))d\gt\leq C(\ge) +c_9n^{\frac{(2-q)p}{q-1}-N+k}.
$$
Since $\frac{(2-q)p}{q-1}<\frac q{q-1}\leq N-k$, we have that $\frac{(2-q)p}{q-1}-N+k<0$. Letting $n\to\infty$ we obtain that $u^p\in L^1(TUB_1(\Gs))$.\smallskip

\nind For the second assertion we have with the same test function $\eta_n$, 
$$\int_{TUB_{1}(\Gs)}\nabla u.\nabla\eta_n dx+m\int_{TUB_{1}(\Gs)}|\nabla u|^q\eta_n dx=\int_{TUB_{1}(\Gs)}u^p\eta_n dx.
$$
Using $(\ref{R1})$ and $(\ref{R5W})$,
$$\left|\int_{TUB_{1}(\Gs)}\nabla u.\nabla\eta_n dx\right|\leq Cn^{\frac{q}{q-1}}Vol(TUB_\gt(1/n))=C'n^{\frac{q}{q-1}+k-N}.
$$
By assumption $\frac{q}{q-1}\leq N-k$. Since $u\in L^p(TUB_{1}(\Gs))$ we conclude that $|\nabla u|\in L^q(TUB_{1}(\Gs))$ by Fatou's lemma.\smallskip

\nind {\it Step 3: We claim that  $u\in L^\infty(TUB_{1}(\Gs))$}. The proof that $u$ is a weak solution of $(\ref{I-1})$ is similar to the one in \rth{Theorem G}. For obtaining that 
$u\in L^\infty(TUB_{1}(\Gs))$ we use the same test functions $\eta_n$ as in Step 2, the same sequence $\{r_j\}$ defined by $(\ref{Rem4*})$ and derive $(\ref{Rem8})$ where $B_R$ is replaced by $TUB_{1}(\Gs)$ under the assumption $(\ref{Rem7})$.
And similarly $(\ref{Rem13})$, again replacing $B_R$ by $TUB_{1}(\Gs)$ holds in the same way, we obtain a contradiction.

\qeda
\medskip

The next theorem extends a previous result of Brezis and Nirenberg \cite{BrNi} that they proved in the case $q=2$. The technique is completely different from the one used in \rth{remov-2} and based upon capacity theory.

\bth{remov-3}Let $\Gw\subset\BBR^N$ $N\geq 2$,  be a bounded smooth domain. Assume $p$ and $q$ are real numbers such that $0< p\leq\max\{2,p\}\leq q$ and $m>0$. Let $K\subset\Gw$ be a compact set and $u\in C^1(\overline \Gw\setminus K)$
be a positive function satisfying
\bel{e1}
-\Gd u+m|\nabla u|^q-u^p\leq 0
\ee
in $\Gw\setminus K$ and such that $u\geq \gd>0$. If $cap_{1,q'}(K)=0$, then $u\in L^\infty(\Gw)$. 
\es

\nind\Proof If $cap_{1,q'}(K)=0$, then $|K|=0$ and  there exists a sequence $\{\gz_k\}\subset C^\infty_c(\Gw)$ such that $0\leq \gz_k\leq 1$, $\gz_k=1$ in a neighborhoood of 
$K$ such that  
\bel{e4}\dsps
\lim_{k\to\infty}\norm{|\nabla\gz_k|}_{L^{q'}(\Gw)}=0.
\ee
Furthermore $\gz_k\to 0$ a.e. in $\Gw$, and we set $\eta_k=1-\gz_k$. For $\gth>0$ let $j_\gth$ be a $C^\infty(\BBR)$ nondecreasing  function
with value $0$ on $(-\infty,0]$ and $1$ on $[\gth,\infty)$. We set  
$$\gl(t)=meas\{x\in\Gw:u(x)\geq t\}$$ 
for $t\geq t_0$ where $t_0=\sup_{\prt\Gw} u\geq\gd$. Taking $\eta^{q'}_kj_\gth(u-t)u^{-p}$ 
as a test function, we have 
$$\BA{lll}\dsps
q'\int_{\Gw}\eta_k^{q'-1}j_\gth(u-t)u^{-p}\nabla u.\nabla \eta_kdx+\int_{\Gw}j'_{\gth}(u-t)u^{-p}|\nabla u|^2\eta^{q'}_kdx
\\[4mm]\phantom{---}
\dsps-p\int_{\Gw}\eta_k^{q'}j_\gth(u-t)u^{-p-1}|\nabla u|^2dx
+m\int_\Gw \eta^{q'}_kj_\gth(u-t)u^{-p}|\nabla u|^qdx\leq \int_\Gw  \eta^{q'}_kj_\gth(u-t)dx.
\EA$$
Since $j'_\gth\geq 0$, it follows
\bel{e5}\BA{lll}\dsps
q'\int_{\Gw}\eta_k^{q'-1}j_\gth(u-t)u^{-p}\nabla u.\nabla \eta_kdx-p\int_{\Gw}\eta^{q'}_kj_\gth(u-t)u^{-p-1}|\nabla u|^2dx
\\[4mm]\phantom{--------}
\dsps
+m\int_\Gw \eta^{q'}_kj_\gth(u-t)u^{-p}|\nabla u|^qdx\leq \int_\Gw \eta^{q'}_kj_\gth(u-t)dx\leq\gl(t).
\EA
\ee
\medskip

\nind {\it Step 1: the basic inequality.} We set
\bel{e5-1}\BA{lll}\dsps
S(t)=\left\{\BA{lll}\frac{q}{q-p}t^{\frac{q-p}{q}}\qquad&\text{if }p<q\\[3mm]
\ln t &\text{if }p=q.
\EA\right.
\EA
\ee
Then $u^{-p}|\nabla u|^q=|\nabla S(u)|^q$ and 
\bel{e6}\BA{lll}\dsps
m\int_\Gw \eta^{q'}_kj_\gth(u-t)|\nabla(S(u))|^q dx\leq \gl(t)+q'\int_{\Gw}\eta_k^{q'-1}j_\gth(u-t)u^{-p}|\nabla u||\nabla \eta_k|dx
\\[4mm]\phantom{------------------}\dsps
+p\int_{\Gw}\eta^{q'}_kj_\gth(u-t)u^{-p-1}|\nabla u|^2dx.
\EA
\ee
We take $t\geq t_1\geq t_0$ for some $t_1$ to be fixed, then
\bel{e6-1}\BA{lll}\dsps
q'\int_{\Gw}\eta_k^{q'-1}j_\gth(u-t)u^{-p}|\nabla u||\nabla \eta_k|dx=q'\int_{\Gw}\eta_k^{q'-1}j_\gth(u-t)u^{-\frac{p(q-1)}{q}}u^{-\frac pq}|\nabla u||\nabla \eta_k|dx
\\[4mm]
\phantom{\dsps
q'\int_{\Gw}\eta_k^{q'-1}j_\gth(u-t)u^{-p}|\nabla u||\nabla \eta_k|dx}\dsps
\leq q't_1^{-\frac{p(q-1†)}{q}}\int_{\Gw}\eta_k^{q'}j_\gth(u-t)|\nabla S(u)|\frac{|\nabla\eta_k|}{\eta_k}dx\\[4mm]
\phantom{-----}
\dsps
\leq q't_1^{-\frac{p(q-1)}{q}}\left(\frac{\ge^q}{q}\int_{\Gw}\eta_k^{q'}j_\gth(u-t)|\nabla S(u)|^qdx+\frac{1}{q'\ge^{q'}}\int_{\Gw}j_\gth(u-t)|\nabla \eta_k|^{q'}dx\right).
\EA\ee
We recall that $\gs=(p+1)q-2p$. Since $q\geq 2$ we have that $\gs\geq 2$, with strict inequality if $q>2$. Therefore
\bel{e7}\BA{lll}\dsps
p\int_{\Gw}\eta^{q'}_kj_\gth(u-t)u^{-p-1}|\nabla u|^2dx=p\int_{\Gw}\eta^{q'}_kj_\gth(u-t)u^{-\frac{\gs}{q}}u^{-\frac{2p}{q}}|\nabla u|^2dx
\\[4mm]\phantom{\dsps p\int_{\Gw}\eta^{q'}_kj_\gth(u-t)u^{-p-1}|\nabla u|^2dx}
\dsps \leq pt_1^{-\frac\gs q}\int_{\Gw}\eta^{q'}_kj_\gth(u-t)|\nabla S(u)|^2dx.
\EA\ee
We first consider the case  $q> 2$. We have by H\"older's inequality,
\bel{e8}\BA{lll}\dsps
p\int_{\Gw}j_\gth(u-t)u^{-p-1}|\nabla u|^2\eta^{q'}_kdx\leq pt_1^{-\frac\gs q}
\left(\frac{2\ge^q}{q}\int_{\Gw}j_\gth(u-t)|\nabla S(u)|^q\eta^{q'}_kdx\right.\\[4mm]
\phantom{\dsps----
p\int_{\Gw}j_\gth(u-t)u^{-p-1}|\nabla u|^2\eta^{q'}_kdx}
\left.\dsps+\frac{q}{(q-2)\ge^{\frac{q}{q-2}}}\int_{\Gw}j_\gth(u-t)\eta^{q'}_kdx\right).
\EA\ee 
We then deduce that
\bel{e9}\BA{lll}\dsps
\left(m-\ge^q\left(\frac{2p}{q}t_1^{-\frac\gs q}+\frac{1}{q-1}t_1^{-\frac{p(q-1)} q}\right)\right)\int_\Gw \eta^{q'}_kj_\gth(u-t)|\nabla(S(u))|^q dx
\\[4mm]\phantom{--------}
\dsps\leq \left(1+\frac{pq}{(q-2)\ge^{\frac{q}{q-2}}}\right)\gl(t)+\frac{t_1^{-\frac{p(q-1)}{q}}}{\ge^{q'}}\int_\Gw j_\gth(u-t)|\nabla\eta_k|^{q'}dx\\[4mm]\phantom{--------}
\dsps\leq \left(1+\frac{pq}{(q-2)\ge^{\frac{q}{q-2}}}\right)\gl(t)+\frac{t_1^{-\frac{p(q-1)}{q}}}{\ge^{q'}}\int_\Gw|\nabla\eta_k|^{q'}dx.
\EA
\ee
Since $cap_{1,q'}(K)=0$ and $\eta_k\to 1$, we let $k\to\infty$ and obtain 
\bel{e10}\BA{lll}\dsps
\left(m-\ge^q\left(\frac{2p}{q}t_1^{-\frac\gs q}+\frac{1}{q-1}t_1^{-\frac{p(q-1)} q}\right)\right)\int_\Gw j_\gth(u-t)|\nabla(S(u))|^q dx
\leq \left(1+\frac{pq}{(q-2)\ge^{\frac{q}{q-2}}}\right)\gl(t),
\EA
\ee
having fixed $t_1\geq t_0$ and $\ge>0$ small enough such that 
$$m-\ge^q\left(\frac{2p}{q}t_1^{-\frac\gs q}+\frac{1}{q-1}t_1^{-\frac{p(q-1)} q}\right)\geq\frac m2.
$$
We set 
$$\gn(s)=meas\{x\in\Gw:\,S(u(x))\geq s\}.$$ 
By letting $\gth\to 0$ we infer that there exists a constant $C_1>0$ such that,
 for $s\geq s_1=S(t_1)$, 
\bel{e11}\BA{lll}\dsps
\int_{\Gw}|\nabla(S(u)-s)_+|^qdx\leq C_1\gn(s).
\EA
\ee
Before continuing on this inequality, we can look at the case $q=2$ (which is actually the case considered by Brezis and Nirenberg \cite{BrNi}).
Then $\gs=2$ and (\ref{e9}) is replaced by 
\bel{e12}\BA{lll}\dsps
\left(m-\left(2pt_1^{-1}-\ge^2t_1^{-\frac{p} 2}\right)\right)\int_\Gw \eta^{2}_kj_\gth(u-t)|\nabla(S(u))|^2 dx
\leq \gl(t)+\frac{t_1^{-\frac{p}{2}}}{\ge^{2}}\int_\Gw|\nabla\eta_k|^{2}dx.
\EA
\ee
By choosing $\ge$ and $t_1$ we obtain (\ref{e11}) with $q=2$ and a specific constant $C_1$.\medskip

\nind {\it Step 2: end of the proof.}
 We set $w=S(u)$ and by H\"older's inequality since $q>2$,
\bel{e13}\BA{lll}\dsps
\int_{\Gw}|\nabla(w-s)_+|^{q'}dx\leq \left(\int_{\Gw}|\nabla(w-s)_+|^{q}dx\right)^{\frac{q'}{q}}\left(meas\left\{|\nabla (w-s)_+>0|\right\}\right)^{1-\frac{q'}{q}}\\[2mm]
\phantom{\dsps
\int_{\Gw}|\nabla(w-s)_+|^{q'}dx}\dsps
\leq c_1^{\frac {q'}{q}}(\gn(s))^{\frac {q'}{q}}\left(meas\left\{|\nabla (w-s)_+>0|\right\}\right)^{1-\frac{q'}{q}}\\[2mm]
\phantom{\dsps
\int_{\Gw}|\nabla(w-s)_+|^{q'}dx}\dsps
\leq c_1^{\frac {q'}{q}}\gn(s),
\EA
\ee
since $\nabla (w-s)_+=0$ a.e. on the set where $(w-s)_+=0$. This implies that, up to a set of zero measure, we have  $\left\{|\nabla (w-s)_+>0|\right\}\subset \left\{ (w-s)_+>0\right\}$,
thus $meas\left\{|\nabla (w-s)_+>0|\right\}\leq\gn(s)$. Note that this also holds if $q=2$. By Sobolev inequality, 
\bel{e14}\BA{lll}\dsps
\left(\int_\Gw(w-s)_+^{q'^*}dx\right)^{\frac{q'}{q'^*}}\leq c(N,q)\int_{\Gw}|\nabla(w-s)_+|^{q'}dx\quad \text{with }\;q'^*=\frac{Nq'}{N-q'},
\EA
\ee
if $q'<N$ which is always satisfied except in the case $q=2=N$ in which case the modifications are straightforward and left to the reader. Furthermore
$$\int_\Gw(w-s)_+dx\leq \left(\int_\Gw(w-s)_+^{q'^*}dx\right)^{\frac{1}{q'^*}}(\gn(s))^{1-\frac{1}{q'^*}}.
$$
This yields
\bel{e15}\BA{lll}\dsps
\int_\Gw(w-s)_+dx\leq c_2\gn(s))^{1+\frac{1}{N}}\quad\text{for any }s\geq s_1,
\EA
\ee
since $1+\frac{1}{q'}-\frac{1}{q'^*}=1+\frac 1N$. Set 
$$\gf(s)=\int_\Gw(w-s)_+dx=\int_s^\infty\gn(\gt)d\gt,\,\text{ hence}\; -\gf'(s)=\gn(s),$$
and (\ref{e15}) leads to 
$\gf(s)\leq c_2(-\gf'(s))^{\frac{N+1}{N}}$ and we finally obtain the following differential inequality
\bel{e16}\BA{lll}\dsps 
\gf'+c_2^{\frac{N}{N+1}}\gf^{\frac{N}{N+1}}\leq 0\quad\text{on } [s_1,\infty).
\EA
\ee
The solution is explicit:
\bel{e17}\dsps
\gf(s)\leq\left\{\BA{lll}\left((\gf(s_1))^\frac{1}{N+1}-\frac{c_2^\frac{N}{N+1}}{N}(s-s_1)\right)^{N+1}\quad&\text{if }s_1\leq s\leq s_2,\\
0&\text{if } s> s_2
\EA\right.
\ee
where 
$$\dsps s_2=s_1+Nc_2^{-\frac{N}{N+1}}(\gf(s_1))^\frac{1}{N+1}.
$$
Hence $(w-s)_+=0$ if $s\geq s_2$ which implies the claim.\qeda\medskip

\nind{\it Proof of \rth {Theorem I}}. If $u$ is a solution the assumption that $u\geq \gd >0$ can be replaced by $u\geq 0$ since $u+\gd$ is a subsolution. 
It is standard that if $u$ is bounded and $cap_{1,q'}(K)$ is zero then it is a weak solution. \qeda\medskip

Motivated by the result of \rth{Theorem G} when $K$ is a single point, we have the following conjecture.\medskip

\nind{\bf Conjecture}. {\it Let $\Gw\subset\BBR^N$ be a bounded smooth domain. Assume $p,q$ are such that $1\leq p\leq q<2$ and $m>0$. Let $K\subset\Gw$ be a compact set and $u\in C^1(\overline \Gw\setminus K)$
be a nonnegative solution of 
\bel{f1}
-\Gd u+m|\nabla u|^q-u^p= 0
\ee
in $\Gw\setminus K$. If $cap_{1,q'}(K)=0$, then $u$ is a weak solution of $(\ref{f1})$ in $\Gw$ and it belongs to $L^\infty(\Gw)$. }

\mysection{Asymptotics of solutions}
The natural way for studying the singular or asymptotic behaviour of solutions of $(\ref{I-1})$ is to use  the spherical coordinates $(r,\gth)\in [0,\infty)\ti S^{N-1}$. Denoting $u(x)= u(r,\gth)$, equation $(\ref{I-1})$ endows the form 
\bel{S1}
-u_{rr}-\frac{N-1}{r}u_r-\frac{1}{r^2}\Gd'u+m\left(u_r^2+\frac{1}{r^2}|\nabla 'u|^2\right)^{\frac q2}-u^p=0,
\ee
where $\Gd'$ and $\nabla '$ represent respectively the Laplace Beltrami operator and the covariant gradient identified with the tangential derivative on the unit sphere. 
This equation admits separable solutions i.e. solutions under the form $u(r,\gth)=r^{-a}\gw(\gth)$ if and only if $q=\frac {2p}{p+1}$, in which case 
$$a=\ga=\gb=\gg.
$$
Then $\gw$ is a nonnegative solution of 
\bel{S2}
-\Gd'\gw-\ga\left(\ga +2-N\right)\gw+m\left(\ga^2\gw^2+|\nabla'\gw|^2\right)^{\frac {p}{p+1}}-\gw^p=0\quad\text{in } S^{N-1}.
\ee
When $q\neq \frac{2p}{p+1}$, one nonlinear term could dominate the other thus the asymptotics can be described either by the separable solutions of the Lane-Emden equation $(\ref{I-3})$ or the Riccatti equation $(\ref{I-5})$. For the Lane-Emden equation the separable solutions have the form $u(r,\gth)=r^{-\ga}\gw(\gth)$ where $\gw$ is a positive solution of 
  \bel{S3}
-\Gd'\gw-\ga\left(\ga+2-N\right)\gw-\gw^p=0\quad\text{in } S^{N-1},
\ee
while for the Riccatti equation the separable solutions are under the form $u(r,\gth)=r^{-\gb}\gw(\gth)$ where $\gw$ is a positive solution of 
  \bel{S4}
-\Gd'\gw-\gb\left(\gb+2-N\right)\gw +m\left(\gb^2\phi^2+|\nabla'\gw|^2\right)^{\frac {q}{2}}=0\quad\text{in } S^{N-1}.
\ee
Separable nonnegative solutions of the eikonal equation $(\ref{I-6})$ have the form $u(r,\gth)=r^{-\gg}\gw(\gth)$ and $\gw$ satisfies
  \bel{S5}
 m\left(\gg^2\gw^2+|\nabla'\gw|^2\right)^{\frac {q}{2}}-\gw^p=0\quad\text{in } S^{N-1}.
\ee
We recall below some results concerning these equations.
\bth{Sep1} Let $N\geq 2$, $p,q>1$ and $m\geq 0$. \smallskip

\nind 1- Suppose $q=\frac{2p}{p+1}$. \\
1-a If $N\geq 3$, $p\geq\frac N{N-2}$ and $m>0$ there exists a unique positive constant solution $x_m$ to $(\ref{S2})$.\\
1-b If $N=2$ and $p>1$, or $N\geq 3$ and $1<p<\frac{N}{N-2}$ there exists no positive constant solution to $(\ref{S2})$ if 
$0\leq m <\gm^*$, a unique positive constant solution $x_{\gm^*}$ if $m=\gm^*$ and two positive constant solutions 
$x_{1,m}<x_{2,m}$ if $m>\gm^*$, where
  \bel{S6}
\gm^*:=(p+1)\left(\frac{N-(N-2)p}{2p}\right)^{\frac{p}{p+1}}.
\ee
\smallskip

\nind 2-  There exist positive solutions to $(\ref{S3})$ if and only if $p>\frac{N}{N-2}$. Furthermore, if 
$\frac{N}{N-2}<p<\frac{N+1}{N-3}$, the positive solutions are constant and therefore unique with value 
  \bel{S7-0}
\gw_0=\left(\ga(N-2-\ga)\right)^{\frac 1{p-1}}=\left(\ga\frac{(N-2)p-N}{p-1}\right)^{\frac{1}{p-1}}.
\ee
\smallskip

\nind 3- If $m>0$ and $1<q<\frac{N}{N-1}$ there exists a unique positive solution to $(\ref{S4})$. This solution is constant with value
  \bel{S7}
\xi_m=\frac 1\gb\left(\frac{(N-1)q-N}{m(q-1)}\right)^{\frac 1{q-1}}.
\ee
If $q\geq \frac{N}{N-1}$ there exists no  positive solution to $(\ref{S4})$.\smallskip

\nind 4- If $m>0$ and $p,q>1$, $p\neq q$, any positive solution to $(\ref{S5})$ is constant with value
  \bel{S8}
X_m=\left(m|\gg|^q\right)^{\frac{1}{p-q}}.
\ee
\es\medskip

\nind\Remark Assertion 1 is proved in \cite[Proposition 6.1]{BiGaVe2}, assertion 2 in \cite{GS}, assertions 3 and 4 are easy consequences of the study of the extrema of a positive smooth solution. 
\subsection{Isolated singularities}

In this Section we obtain the precise behaviour of positive singular solutions of $(\ref{I-1})$ in $B_{r_0}\setminus\{0\}$.\smallskip
\subsubsection{Proof of \rth{Theorem J}}
The proof is a delicate combination of various techniques, some new and some other already which have already been used by the authors in several different contexts. \\
Up to change of scale we assume that ${r_0}=1$. Set 
\bel{S00}u(r,\gth)=r^{-\ga}v(t,\gth)\quad\text{with }t=\ln r,\; t\leq 0.
\ee
The function $v$ satisfies
\bel{SI2}\BA{lll}
v_{tt}+(N-2-2\ga )v_t+\ga\left(\ga +2-N\right)v+\Gd'v \\[3mm]
\phantom{---------)--}-me^{-\frac{\gs t}{p-1}}\left((v_t-\ga v)^2+|\nabla'v|^2\right)^{\frac {q}{2}}+v^p=0,
\EA\ee
in $(-\infty,0]\ti S^{N-1}$, recalling that $\gs=(p+1)q-2p$. By Theorem B the functions $v$, $v_t$ and $|\nabla ' v|$ is bounded in $(-\infty,0]\ti S^{N-1}$. By standard regularity estimates and Ascoli-Arzela theorem the limit set 
 at $-\infty$ of the trajectory of $v$ in $C^2(S^{N-1})$,
$$\dsps\CT_-[v]=\bigcup _{t\leq0 }\{v(t,.)\},
$$
is a non-empty compact connected subset $\Gg_{-}$ of $C^2(S^{N-1})$. Set
$$\CE[v](t)=\frac{1}{2}\int_{S^{N-1}}\left(v_t^2-|\nabla 'v|^2+\ga\left(\ga +2-N\right)v^2+\frac{2}{p+1}|v|^{p+1}\right)dS,
$$
then 
$$\frac{d}{dt}\CE[v](t)=-(N-2-2\ga)\int_{S^{N-1}}v_t^2dS-me^{-\frac{\gs t}{p-1}}\int_{S^{N-1}}\left((v_t-\ga v)^2+|\nabla'v|^2\right)^{\frac {q}{2}}v_tdS.
$$
Therefore, for any $t<0 $,
\bel{SI3-1}\BA{lll}\dsps
\CE[v](t)-\CE[v](0 )=(N-2-2\ga )\int_t^0 \int_{S^{N-1}}v_t^2dSd\gt\\[4mm]
\phantom{---------}\dsps+m\int_t^0  e^{-\frac{\gs \gt}{p-1}}\int_{S^{N-1}}\left((v_t-\ga v)^2+|\nabla'v|^2\right)^{\frac {q}{2}}v_tdSd\gt.
\EA\ee
Since $\CE[v](t)$ and $\left((v_t-\ga v)^2+|\nabla'v|^2\right)^{\frac {q}{2}}$ are uniformly bounded, $N-2-2\ga\neq 0$ because $p\neq\frac{N+2}{N-2}$ and $\gs<0$, this implies that 
\bel{SI3-2}\BA{lll}\dsps
\int_{-\infty}^0 \int_{S^{N-1}}v_t^2dSd\gt<\infty.\EA\ee
Since $v_t$ is uniformly continuous on $(-\infty,0 ]\ti S^{N-1}$, it implies in turn that 
$$\dsps\lim_{t\to-\infty}\int_{S^{N-1}}v_t^2(t)dS=0.
$$
Multiplying the equation $(\ref{SI2})$ by $v_{tt}$, using the $C^2$ estimate on $v$ and $(\ref{SI3-2})$ we obtain that 
\bel{SI4}\BA{lll}\dsps
\int_{-\infty}^0 \int_{S^{N-1}}v_{tt}^2dSd\gt<\infty,\EA\ee
which implies in turn 
$$\dsps\lim_{t\to-\infty}\int_{S^{N-1}}v_{tt}^2(t)dS=0. 
$$
Letting $t\to-\infty$ in $(\ref{SI2})$ we conclude that $\Gg_{-}$ is a a non-empty compact connected subset of the set on nonnegative solutions of $(\ref{S3})$.\\
 If $1<p\leq \frac N{N-2}$ we have 
\bel{SI5}\BA{lll}\dsps
\lim_{t\to-\infty}v(t,.)=0\quad\text{uniformly on }S^{N-1}.
\EA\ee
If $\frac N{N-2}<p<\frac {N+2}{N-2} $,
\bel{SI6}\BA{lll}\dsps
\text{either }\;\lim_{t\to-\infty}v(t,.)=0\quad\text{or }\lim_{t\to-\infty}v(t,.)=\gw_0\quad\text{uniformly on }S^{N-1}.
\EA\ee
where $\gw_0$ is defined by $(\ref{S7-0})$.\smallskip

\nind The remaining problem is to analyse the case where $\dsps \lim_{t\to-\infty}v(t,.)=0$. This is delicate and presented in the following lemmas.

\blemma{PSing1} Let $N\geq 3$, $p\in \left(1,\infty\right)\setminus \left\{\frac {N}{N-2},\frac {N+2}{N-2}\right\}$ and $1<q<\frac{2p}{p+1}$. If $u$ is a nonnegative solution of $(\ref{I-1})$ in $B_2\setminus\{0\}$, such that 
\bel{SI7}\BA{lll}\dsps
\dsps\lim_{x\to 0}|x|^\ga u(x)=0,
\EA\ee
then there exists $\ge>0$ such that 
\bel{SI8}\BA{lll}\dsps
\dsps u(x)\leq C|x|^{-\ga+\ge}\qquad\text{for all }x\in B_{1}\setminus\{0\}.
\EA\ee
Furthermore 
\bel{SI8^*}\BA{lll}\dsps
\dsps |\nabla u(x)|\leq C'|x|^{-\ga-1+\ge}\qquad\text{for all }x\in B_{1}\setminus\{0\}.
\EA\ee
\es
\Proof The key point is the proof is that under the assumptions on $p$ the coefficients $\ga(\ga+2-N)$ and $N-2-2\ga$ in the equation $(\ref{SI2})$ satisfied by the function $v$ defined before are not zero.  
 We note that $(\ref{SI8})$ is equivalent to 
\bel{SI9}v(t,\gth)\leq Ce^{\ge t}\qquad\text{for all }(t,\gth)\in (-\infty,0]\ti S^{N-1}.\ee
If $(\ref{SI9})$ does not hold we have that
$$
\limsup_{t\to-\infty}e^{-\ge t}\gr(t)=+\infty\quad\text{for all }\ge>0,
$$
where $\gr(t)=\sup\{v(t,\gth):\gth\in S^{N-1}\}$. We use now a technique introduced in \cite[Lemma 2.1]{CMV}: it is proved that there exists a function $\eta\in C^\infty\big((-\infty,0]\big)$ such that 
  \bel{BB4}\BA{lll}
&(i)\qquad\displaystyle &\displaystyle\eta>0,\,\eta'>0,\,\lim_{t\to-\infty}\eta(t)=0;\qquad\qquad\qquad\qquad\qquad\qquad\qquad\qquad\quad\quad\quad\\[3mm]
&(ii)\qquad\displaystyle &\displaystyle0<\limsup_{t\to-\infty}\frac{\gr(t)}{\eta(t)}<+\infty;\\[3mm]
&(iii)\qquad\displaystyle &\displaystyle\lim_{t\to-\infty}e^{-\vge t}\eta(t)=+\infty\quad\text{for all }\vge>0;\\[3mm]
&(iv)\qquad\displaystyle &\displaystyle\left(\frac{\eta'}{\eta}\right)',\, \left(\frac{\eta''}{\eta}\right)'\in L^1((-\infty,0));\\[3mm]
&(v)\qquad\displaystyle &\displaystyle\lim_{t\to-\infty}\frac{\eta'(t)}{\eta(t)}=\lim_{t\to-\infty}\frac{\eta''(t)}{\eta(t)}=0.
\EA\ee
We define $\psi$ by $v(t,\cdot)=\eta(t)\psi(t,.)$, then 
\bel{SI10}\BA{lll}\dsps
\psi_{tt}+K_1\psi_t+K_2\psi+\Gd'\psi -me^{-\frac{\gs t}{p-1}}\eta^{q-1}\left(\left(\psi_t-\ga\frac{\eta_t}{\eta} \psi\right)^2+|\nabla'\psi|^2\right)^{\frac {q}{2}}\\[4mm]
\phantom{------------------}+\eta^{p-1}\psi^p=0\quad\text{in }(-\infty,0]\ti S^{N-1},
\EA\ee
where 
$$K_1(t)= N-2-2\ga+2\frac{\eta'}{\eta}\quad{\rm and}\quad K_2(t)= \ga(\ga +2-N)+(N-2-2\ga )\frac{\eta'}{\eta}+\frac{\eta''}{\eta}.$$
The function $\psi$ is bounded and by standard regularity estimates it is uniformly bounded in the $C^2$-topology of $(-\infty,0 ]\ti S^{N-1}$.
We set
$$\tilde\CE[\psi](t)=\frac{1}{2}\int_{S^{N-1}}\left(\psi_t^2-|\nabla '\psi|^2-\ga\left(\ga +2-N\right)\psi^2\right)dS,
$$
then 
\bel{BB4-1}\BA{lll}\dsps
\frac{d}{dt}\tilde\CE[\psi](t)=-\left(N-2-2\ga+2\frac{\eta'}{\eta}\right)\int_{S^{N-1}}\psi_t^2dS+\left((N-2-2\ga)\frac{\eta'}{\eta}+\frac{\eta''}{\eta}\right)\int_{S^{N-1}}\psi\psi_t dS\\[4mm]
\dsps -\eta^{p-1}\int_{S^{N-1}}\psi^p\psi_t dS
+me^{-\frac{\gs t}{p-1}}\eta^{q-1}\int_{S^{N-1}}\left(\left(\psi_t-\ga\frac{\eta_t}{\eta} \psi\right)^2+|\nabla'\psi|^2\right)^{\frac {q}{2}}\psi_tdS.
\EA\ee
We analyse the different terms in the right-hand side of $(\ref{BB4-1})$:
$$\int_{S^{N-1}}\psi^p\psi_t dS=\frac{1}{p+1}\frac{d}{dt}\int_{S^{N-1}}\psi^{p+1}\eta^{p-1}-\frac{p-1}{p+1}\eta'\eta^{p-2}\int_{S^{N-1}}\psi^{p+1}dS.
$$
By the mean value theorem, for any $t<0 $ there exists $t^*\in (t,0 )$ such that 
$$\BA{lll}\dsps\int_t^0 \int_{S^{N-1}}\eta^{p-1}\int_{S^{N-1}}\psi^p\psi_t dSd\gt=\frac{1}{p+1}\left[\int_{S^{N-1}}\psi^{p+1}\eta^{p-1}\right]_t^0 \\[4mm]
\phantom{-----------}-\dsps
\frac{1}{p+1}\left(\eta^{p-1}(0)-\eta^{p-1}(t)\right)\int_{S^{N-1}}\psi^{p+1}(t^*,.) dS,
\EA$$
and this expression is bounded independently of $t<0 $. Also
$$\BA{lll}\dsps
\left((N-2-2\ga)\frac{\eta'}{\eta}+\frac{\eta''}{\eta}\right)\int_{S^{N-1}}\psi\psi_t dS=\frac{1}{2}\frac{d}{dt}\left(\left((N-2-2\ga)\frac{\eta'}{\eta}+\frac{\eta''}{\eta}\right)
\int_{S^{N-1}}\psi^{2}dS\right)\\[4mm]
\phantom{-----------------}\dsps
-\frac{1}{2}\left((N-2-2\ga)\left(\frac{\eta'}{\eta}\right)'+\left(\frac{\eta''}{\eta}\right)'\right)\int_{S^{N-1}}\psi^{2}dS.
\EA
$$
The term involving the gradient is clearly integrable on $(-\infty,0 )$. Hence we obtain for any $t<0 $, 
\bel{SI11}\BA{lll}\dsps
\tilde\CE[\psi](0 )-\tilde\CE[\psi](t)=-\int_{t}^0 \left(N-2-2\ga+2\frac{\eta'}{\eta}\right)\int_{S^{N-1}}\psi_t^2dSd\gt +A(t)
\EA
\ee
where $A(t)$ is bounded independently of $t<0 $. Because the left-hand side of $(\ref{SI11})$ is bounded independently of $t<0 $, $\frac{\eta'}{\eta}(\gt)\to 0$ when $\gt\to-\infty$ and $N-2-2\ga \neq 0$ as $p\neq \frac{N+2}{N-2}$, we infer that
\bel{SI12}\BA{lll}\dsps
\int_{-\infty}^0 \int_{S^{N-1}}\psi_t^2dSd\gt<\infty.
\EA
\ee
By uniform continuity, this implies that $\psi_t(t)\to 0$ in $L^2(S^{N-1})$ when $t\to-\infty$. Multiplying the equation satisfied by $\psi_{tt}$ we obtain similarly, using 
the previous estimate and $(\ref{BB4})$-(iv)-(v) that 
\bel{SI13}\BA{lll}\dsps
\int_{-\infty}^0 \int_{S^{N-1}}\psi_{tt}^2dSd\gt<\infty;
\EA
\ee
in turn this implies that $\psi_{tt}(t)\to 0$ in $L^2(S^{N-1})$ when $t\to-\infty$. The limit set at $-\infty$ of the trajectory $\CT_-[\psi]$ is a connected and compact subset of the set of nonnegative solutions of 
\bel{SI14}\BA{lll}\dsps
\ga(\ga+2-N)\gw+\Gd'\gw=0\quad\text{in } S^{N-1}.
\EA
\ee
Since $\ga(\ga+2-N)$ is not an eigenvalue of $-\Gd'$ in $W^{1,2}(S^{N-1})$, it follows that $\gw=0$, which contradicts the fact that by $(\ref{BB4})$-(ii)  the limit set contains at least one non-zero positive element. Hence $(\ref{SI8})$ holds, as for  $(\ref{SI8^*})$ it is a consequence of \rth{pre}. This ends the proof.\qeda
\medskip


\blemma{PSing2} Let the assumptions of \rth {Theorem J} hold, then\\
\nind 1- If  $N\geq 3$ and $1<p<\frac{N}{N-2}$ (resp. $N=2$ and $p>1$) there exists $k\geq 0$ such that $|x|^{N-2}u(x)$   (resp. $-u(x)/\ln |x|$) converges to $k$  when $x\to 0$. Furthermore $u$ satisfies $(\ref{SI1})$.

\nind 2- If $N\geq 3$ and $\frac{N}{N-2}<p<\frac{N+2}{N-2}$, \\
 2-(i) either $|x|^\ga u(x)$ converges to $\gw_0$ when $x\to 0$,\\
2-(ii) or $u$ is a classical solution of $(\ref{I-1})$ in $B_{r_0}$.
\es
\Proof Since $|x|^\ga u(x)+|x|^{\ga+1}|\nabla u(x)|$ remains bounded and $q\leq \frac{2p}{p+1}$, we have
\bel{SI15}\BA{lll}\dsps
|x|^2u^{p-1}(x)+|x||\nabla u(x)|^{q-1}\leq c_1\quad\text{for all }x\in B_{r_0}.
\EA
\ee
Hence Harnack inequality is valid uniformly on any sphere with center $0$ (see e.g. \cite{GT}) in the sense that
\bel{SI16}\BA{lll}
\dsps \max_{|y|=r}u(y)\leq c_2\min_{|y|=r}u(y)\quad\text{for all }0<r\leq \tfrac {r_0}2.
\EA
\ee
{\it Step 1: first estimate on the average of $v$}. The second order linear equation
\bel{SI17}\BA{lll}
X''+(N-2-2\ga)X'+\ga(\ga+2-N)X=0
\EA
\ee
admits the two linearly independent solutions
$$X_1(t)=e^{\gl_1t} \quad\text{and } X_2(t)=e^{\gl_2t},
$$
where the $\gl_j$ are the roots of $P(\gl)=\gl^2+(N-2-2\ga)\gl+\ga(\ga+2-N)$. Note that  these roots are explicit:
 \bel{SI18}\BA{lll}
\gl_1=\ga>\gl_2=\ga+2-N,
\EA
\ee
and $\gl_2>0$ (resp. $\gl_2<0$) if $1<p<\frac{N}{N-2}$ (resp. $p>\frac{N}{N-2}$). We set
 \bel{SI1H}\BA{lll}
H(t,.)=me^{-\frac{\gs t}{p-1}}\left((v_t-\ga v)^2+|\nabla'v|^2\right)^{\frac {q}{2}}-v^p.
\EA
\ee
Since $\norm{v(t,.)}_{L^{\infty(S^{}N-1)}}+\norm{\nabla' v(t,.)}_{L^{\infty(S^{}N-1)}}\leq Ce^{\ge t}$ by $(\ref{SI8})$-$(\ref{SI8^*})$, there holds
 \bel{SI1H*}\BA{lll}
\norm {H(t,.)}_{L^{\infty(S^{}N-1)}}\leq c_3e^{\gd_1t}
\EA
\ee
where
\bel{gd1}\gd_1=\min\left\{\ge p,\ge q-\tfrac{\gs}{p-1}\right\},
\ee
and $\gs=(p+1)q-2p<0$. Let $\bar v(t)$ and $\overline  H(t)$ be the average respectively of $v(t,.)$ and $H(t,.)$ on $S^{N-1}$. Then $|\overline H(t)|\leq Ce^{\gd_1t}$. Since 
 \bel{SI18*} \bar v''+(N-2-2\ga)\bar v'+\ga(\ga+2-N)\bar v=\overline H(t).
\ee
 Assuming that $\gd_1\neq \gl_1,\gl_2$ (which can always be assume up to changing $\ge$)
the function $\bar v$ endows the general form
\bel{SI19}\BA{lll}\dsps
\bar v(t)=Ae^{\gl_1t}+Be^{\gl_2t}+C(t)e^{\gd_1t},
\EA
\ee
for some constants $A$ and $B$ and for some particular solution $C(t)e^{\gd_1t}$ where $C$ is bounded on $(-\infty,0]$. This can be checked by the so-called method of "the variation of constants".
Therefore, since $v(t,.)\to 0$ when $t\to-\infty$, 
\bel{SI20}\dsps
\bar v(t)=\left\{\BA{lll} Ae^{\gl_1t}+Be^{\gl_2t}+C(t)e^{\gd_1t}\qquad&\text{if }1<p<\frac{N}{N-2}\\[2mm]
Ae^{\gl_1t}+C(t)e^{\gd_1t}\qquad&\text{if }p>\frac{N}{N-2}.
\EA\right.\ee
This leads us to the second decay estimate (besides the one given by \rlemma{PSing1})
\bel{SI21}\BA{lll}\dsps
\bar v(t)\leq c_4e^{\gth_1t}
\EA
\ee
where $\gth_1=\min\left\{\gl_2,\gd_1\right\}$  if $1<p<\tfrac{N}{N-2}$ and $\gth_1=\min\left\{\gl_1,\gd_1\right\}$ if $p>\tfrac{N}{N-2}$. \smallskip

\nind {\it Step 2: first a priori estimate on $v$}. The global estimate on $v$ is obtained by using an iterative method based upon the integral representation of the solutions introduced in \cite{BoVe}. We set 
\bel{SI22}\BA{lll}\dsps
\BBL=-\left(-\Gd'+\tfrac{(N-2)^2}{4}I\right)^{\frac 12},
\EA
\ee
 and let $S(t)=e^{t\BBL}$ be the semigroup of contraction generated by $\BBL$ in $L^2(S^{N-1})$. Introducing the standard Hilbertian decomposition of $H^1(S^{N-1})$ associated to 
 the operator $-\Gd'$, it is classical that
 the space $\BBH=\{\gf\in L^2(S^{N-1}):\bar\gf=0\}$ is invariant by $\BBL$, since $\bar\gf$ is the orthogonal projection in $H^1(S^{N-1})$ onto $(\ker(-\Gd'))^{\perp}=\BBH$. Because
 $$\inf\gs(\BBL\lfloor_\BBH)=\frac {N^2}4,
$$
we have 
\bel{SI23}\BA{lll}\dsps
\norm{S(t)\gf}_{L^2(S^{N-1})}\leq e^{-\frac {Nt}2}\norm\gf_{L^2(S^{N-1})}\quad\text{for all }t>0\text{ and }\gf\in\BBH,
\EA
\ee
and
\bel{SI24}\BA{lll}\dsps
\norm{S(t)\gf}_{L^\infty(S^{N-1})}\leq Ce^{-\frac {Nt}{2}}\norm\gf_{L^\infty{(S^{N-1})}}\quad\text{for all }t>0\text{ and }\gf\in\BBH\cap L^\infty(S^{N-1}).
\EA
\ee
 for some $C>0$. Note that this last inequality is easily obtained by using the Hilbertian decomposition with spherical harmonics. The following representation formula for $v^*=v-\bar v$ is proved in \cite{BoVe}:
 \bel{SI25}\BA{lll}\dsps
v^*(t,.)=e^{\frac{2\ga+2-N}{2}t}S(-t)v^*(0,.)-\int_t^0e^{\frac{2\ga+2-N}{2}s}S(-s)\int_\infty^0 e^{\frac{N-2\ga-2}{2}\gt}S(-\gt)H^*(-t-\gt+s,\gs)d\gt ds
\EA
\ee
 where $ H^*(t,.)=H(t,.)-\overline  H(t)$.
 Since
  \bel{SI27}\norm{H^*(t,.)}_{L^\infty(S^{N-1})}\leq c_3e^{\gd_1t}
 \ee
by $(\ref{SI1H*})$ where $\gd_1$ is defined in $(\ref{gd1})$, we get
  \bel{SI26}\BA{lll}\dsps
\norm{v^*(t,.)}_{L^\infty(S^{N-1})}\leq c_5e^{(\ga+1)t}+c_6e^{\gd_1t}\quad\text{for all }\,t\leq 0.
\EA
\ee
Writing $v(t,.)=\bar v(t)+v^*(t,.)$ we deduce
  \bel{SI29}\BA{lll}\dsps
\norm{v(t,.)}_{L^\infty(S^{N-1})}\leq c_7e^{(\ga+1)t}+c_8e^{\gd_1t}+c_9e^{\gth_1t}\leq c_{10}e^{\gth_1t}\quad\text{for all }\,t\leq 0,
\EA
\ee
where we use the value of $\gth_1$ defined in $(\ref{SI21})$ and $\gl_1,\gl_2$ given in $(\ref{SI18})$. This leads us to an improvement of the decay estimate given by $(\ref{SI9})$. Notice also that if $\gth_1=\gl_2=\ga+2-N$ (resp. $\gth_1=\gl_1=\ga$) when $1<p<\frac{N}{N-2}$ (resp. $\frac{N}{N-2}<p<\frac{N+2}{N-2}$) we deduce from the definition of $v$ that the function $u$ is smaller that $c_{10}|x|^{2-N}$ (resp. is bounded by $c_{10}$). 
\smallskip

\nind{\it Step 3: a priori estimate on $v$ by iterations}. For the sake of understanding we will distinguish two cases according to the sign of $p-\frac{N}{N-2}$.\\
(i) Let $1<p<\frac{N}{N-2}$. Since $v(t,.)\leq c_{10}e^{\gth_1t}$, then by \rth{pre} that $v(t,.)+|\nabla v(t,.)|\leq c_{11}e^{\gth_1t}$. Therefore
$$ \norm {H(t,.)}_{L^\infty(S^{N-1})}\leq c_{12}e^{\delta_2t}
$$
with 
$$\delta_2=\min\left\{\gth_1 p,\gth_1 q-\tfrac{\gs}{p-1}\right\}.$$
Since
 $({\ref{SI18*}})$ holds with $H$ satisfying $(\ref{SI1H*})$ with $\gd_1$ replaced by $\gd_2$, we deduce 
that
 $$\bar v(t)=Ae^{\gl_1t}+Be^{\gl_2t}+C(t)e^{\gd_2t}
 $$
 where $A,B$ are constants and $C$ is bounded which implies
$\gth_2=\min\{\gl_2,\gd_2\}$. Since
 $({\ref{SI18*}})$ holds with $H$ satisfying $(\ref{SI1H*})$ with $\gd_1$ replaced by $\gd_2$
   \bel{SI30}\BA{lll}\dsps
\bar v(t)\leq c_{13}e^{\gth_2t},
\EA
\ee
with $\gth_2=\min\{\gl_1,\gl_2,\gd_2\}=\min\{\gl_2,\gd_2\}$. The integral representation $(\ref{SI25})$ is satisfied by $v^*=v-\bar v$ and we obtain as in the previous step that $(\ref{SI26})$ holds with $\gd_1$ replaced by $\gd_2$ and finally
  \bel{SI31-2}\BA{lll}\dsps
\norm{v(t,.)}_{L^\infty(S^{N-1})}\leq c_{14}e^{(\ga+1)t}+c_{15}e^{\gd_2t}+c_{16}e^{\gth_2t}\leq c_{17}e^{\gth_2t}\quad\text{for all }\,t\leq 0.
\EA
\ee
If $\gth_2=\ga+2-N$ we have the desired estimate, otherwise we iterate. We define the sequences 
   \bel{SI31*}\BA{lll}(i)\quad\quad\qquad\qquad&\gd_1=\min\left\{p\ge,q\ge-\frac{\gs}{p-1}\right\}\,\text{and } \gth_1=\min\{\gl_2,\gd_1\}\\[2mm]
(ii)\qquad&\gd_n=\min\left\{p\gth_{n-1},q\gth_{n-1}-\frac{\gs}{p-1}\right\}\,\text{and } \gth_n=\min\{\gl_2,\gd_n\},\quad\qquad\qquad
\EA\ee
for all the integers $n$ such that $\gd_n<\gl_2$. Then $\gd_n,\gth_n>0$ and the function $v$ satisfies 
  \bel{SI31-n}\BA{lll}\dsps
\norm{v(t,.)}_{L^\infty(S^{N-1})}\leq c_{1,n}e^{(\ga+1)t}+c_{2,n}e^{\gd_nt}+c_{3,n}e^{\gth_nt}\leq c_{4,n}e^{\gth_nt}\quad\text{for all }\,t\leq 0.
\EA
\ee
Furthermore
   \bel{SI32}\gth_n-\gth_{n-1}=\min\left\{\gl_2-\gth_{n-1},\min\left\{(p-1)\gth_{n-1},(q-1)\gth_{n-1}-\frac{\gs}{p-1}\right\}\right\}.
\ee
We assume first that there exists a largest integer  $n_0$ such that $\gth_n< \gl_2$. Then 
$\gth_1<\gth_2<...<\gth_n<...\gth_{n_0}$ and $\gth_{n_0+1}=\gl_2$.\\
If such a largest integer does not exist, then $\{\gth_n\}$ is increasing with limit $\gth_\infty\leq\gl_2$. By $(\ref{SI32})$, $\gth_\infty$ and $\gl_2$ coincide. 
By $(\ref{SI31*})$-(ii),  $\{\gd_n\}$ is increasing. For any $\ge>0$ there exists $n_\ge\in\BBN$ such that $\gl_2-\ge\gth_n<\gl_2$ for $n\geq n_\ge$, hence 
$$\gd_{n_\ge}>\min\left\{p(\gl_2-\ge),q\gl_2-\ge)-\frac{\gs}{p-1}\right\}>\gl_2
$$
if $\ge$ is small enough. This implies that $\gth_{n_\ge}=\gl_2$, contradiction. Therefore inequality $(\ref{SI31-n})$ with $n=n_\ge$ becomes
  \bel{SI33}\BA{lll}\dsps
\norm{v(t,.)}_{L^\infty(S^{N-1})}\leq c_{18}e^{(\ga+2-N)t}\quad\text{for all }\,t\leq 0.
\EA
\ee
(ii) Let $\frac{N}{N-2}<p<\frac{N+2}{N-2}$. The proof differs from the previous one  only with very little modifications. Since $\gl_2<0$,   $(\ref{SI31*})$ is replaced by
\bel{SI34}\BA{lll}(i)\quad\quad\qquad\qquad&\gd_1=\min\left\{p\ge,q\ge-\frac{\gs}{p-1}\right\}\,\text{and } \gth_1=\min\{\gl_1,\gd_1\}\\[2mm]
(ii)\qquad&\gd_n=\min\left\{p\gth_{n-1},q\gth_{n-1}-\frac{\gs}{p-1}\right\}\,\text{and } \gth_n=\min\{\gl_1,\gd_n\}.\quad\qquad\qquad
\EA\ee
Inequality $(\ref{SI31-n})$ holds with the $\gth_n$ defined above, and there exists an integer $n_\ge$ such that $\gth_n=\gl_1=\ga$. Hence 
  \bel{SI35}\BA{lll}\dsps
\norm{v(t,.)}_{L^\infty(S^{N-1})}\leq c_{19}e^{\ga t}\quad\text{for all }\,t\leq 0.
\EA
\ee

\smallskip

\nind {\it Step 4: convergence}. (i) When $1<p<\frac N{N-2}$, the function $H$ defined $(\ref{SI1H})$ satisfies 
  \bel{SI36}\BA{lll}\dsps
\norm{H(t,.)}_{L^\infty(S^{N-1})}\leq c_{20}e^{\tilde\gd t}\quad\text{for all }\,t\leq 0.
\EA
\ee
with $\tilde \gd=\min\{\gl_2p,\gl_2q-\frac{\gs}{p-1}\}$. Hence $|\overline H(t)|$ satisfies the same estimate and $\bar v$ can be written as in $(\ref{SI19})$ with new coefficients $ A$, $B$ and $C(.)$ under the form 
  \bel{SI37-1}\BA{lll}\dsps
\bar v(t)=Ae^{\gl_1t}+Be^{\gl_2t}+C(t)e^{\tilde\gd t}=Be^{\gl_2t}+o(e^{\gl_2t})\quad\text{as }t\to-\infty.
\EA
\ee
Since formulas $(\ref{SI25})$, $(\ref{SI27})$ and $(\ref{SI26})$ holds with $\gd_1$ replaced by $\gd$ we conclude that 
  \bel{SI38}\norm{v^*(t,.)}_{L^\infty(S^{N-1})}=o(e^{\gl_2t})\quad\text{as }t\to-\infty,
\ee
and finally 
  \bel{SI39-1}\dsps
  \lim_{t\to-\infty}e^{(N-2-\ga)t}v(t,.)=B\quad\text{uniformly on }S^{N-1}.
\ee
Equivalently
  \bel{SI40-1}\dsps
  \lim_{x\to 0}|x|^{N-2}u(x)=B.
\ee
Therefore $u\in L^p(B_{r_0})$. We use the same type of cut-off function $\eta_n$ used in the proof of \rth{Theorem G},  except that we assume also that 
$|\Gd\eta_n|\leq cn^2{\bf 1}_{B_{1/n}\setminus B_{1/(2n)}}$,
and we obtain 
  \bel{SI41-1}\dsps
  -\int_{B_{r_0}}u\Gd\eta_n dx+m \int_{B_{r_0}}|\nabla u|^q\eta_n dx=\int_{B_{r_0}}u^p\eta_n dx.
\ee
The right-hand side of $(\ref{SI41-1})$ is bounded from above by $\norm u^p_{L^p(B_{\frac {2r_0}3})}$.  We have also
$$\left|\int_{B_{r_0}}u\Gd\eta_n dx\right|\leq c_{21}n^{2-N-2+N}\leq c_{22}.
$$
By Fatou's lemma we deduce that $\nabla u\in L^q(B_{\frac {2r_0}3})$. Therefore, by the Brezis-Lions Lemma \cite{BrLi} we conclude that there exists $k$ such that 
$(\ref{SI1})$ holds. \smallskip

If $k=0$, then $B=0$ and $(\ref{SI37-1})$ yields
  \bel{SI37-2}\BA{lll}\dsps
\bar v(t)\leq c_{23}e^{\tilde\gth_1t},
\EA
\ee
with $\tilde\gth_1=\min\left\{\gl_1,\tilde\gd\right\}$. Using again the representation $(\ref{SI25})$ combined with $(\ref{SI36})$ we obtain
  \bel{SI29-1}\BA{lll}\dsps
\norm{v(t,.)}_{L^\infty(S^{N-1})}\leq c_{24}e^{(\ga+1)t}+c_{25}e^{\tilde\gd t}+c_{26}e^{\tilde\gth_1t}\leq c_{27}e^{\tilde\gth_1t}\quad\text{for all }\,t\leq 0,
\EA
\ee
We define now the sequence 
\bel{SI34-1}\BA{lll}(i)\quad\quad\qquad\qquad&\tilde\gd_1:=\tilde\gd\,\text{and } \tilde\gth_1=\min\{\gl_1,\tilde \gd_1\}\\[2mm]
(ii)\qquad&\tilde\gd_n=\min\left\{p\tilde\gth_{n-1},q\tilde\gth_{n-1}-\frac{\gs}{p-1}\right\}\,\text{and } \tilde\gth_n=\min\{\gl_1,\tilde\gd_n\},\quad\qquad\qquad
\EA\ee
and we have
  \bel{SI39-n}\BA{lll}\dsps
\norm{v(t,.)}_{L^\infty(S^{N-1})}\leq Ce^{\tilde \gth_nt}\quad\text{for all }\,t\leq 0.
\EA
\ee
By the construction of Step 3-(ii) there exists $n^*$ such that $\tilde\gth_n=\gl_1$ which means that inequality $(\ref{SI35})$ holds and $\bar v$ satisfies
  \bel{SI41-2}\BA{lll}\dsps
\bar v(t)=Be^{\gl_1t}+C(t)e^{\tilde\gd_{n^*} t})=Be^{\gl_1t}+o(e^{\gl_1t})\quad\text{as }t\to-\infty,
\EA
\ee
and
  \bel{SI41-3}\BA{lll}\dsps
\norm{v^*(.,t)}_{L^\infty(S^{N-1})}=o(e^{\gl_1t})\quad\text{as }t\to-\infty,
\EA
\ee
Hence
  \bel{SI42-1}\dsps
  \lim_{t\to-\infty}e^{-\ga t}v(t,.)=A\quad\text{uniformly on }S^{N-1},\,\text{ equivalently }\;  
  \lim_{x\to 0}u(x)=A.
\ee
Using again the same type of cut-off function $\eta_n$ as in the proof of \rth{Theorem G} we obtain successively that $|\nabla u|\in L^q(B_{r_0})$ and that $u$ is a classical solution. 
\smallskip

\nind (ii) When $\frac N{N-2}<p<\frac {N+2}{N-2}$, $(\ref{SI36})$ is valid with $\gd=\tilde \gd=\min\{\gl_1p,\gl_1q-\frac{\gs}{p-1}\}$. Hence the proof of 
(i) when $A=0$ applies and we obtain that $u$ is a bounded classical solution.

 \qeda
\blemma{PSing3} Let the assumptions of \rth {Theorem J} holds with $N\geq 3$ and $p=\frac{N}{N-2}$, then\\
 (i) either $|x|^{N-2}(-\ln |x|)^{\frac{N-2}{2}} u(x)$ converges to $\left(\frac{N-2}{\sqrt 2}\right)^{N-2}$ when $x\to 0$,\\
(ii) or $u$ is a classical solution of $(\ref{I-1})$ in $B_{r_0}$.
\es
\Proof The proof is based upon a combination of several techniques introduced in \cite{Ve2} for analysing the exterior problem
 \bel{SI39-}\BA{lll}\dsps
-\Gd u+|u|^{\frac 2{N-2}}u=0\qquad\text{in }B_{r_0}^c,
\EA
\ee
and adapted in \cite{Avi} to characterise the isolated singularities of 
  \bel{SI39}\BA{lll}\dsps
-\Gd u=u^{\frac N{N-2}}.
\EA
\ee
\nind {\it 1- We claim that $u$ satisfies

  \bel{SI41}\BA{lll}\dsps
u(x) \leq C|x|^{2-N}(-\ln |x|)^{\frac{2-N}{2}}
\EA
\ee
for $0<|x|\leq r_1$ where $r_1<\min\left\{1,\frac{r_0}{2}\right\}$}.
\smallskip

\nind The function $v$ which is defined by $(\ref{S00})$ with $\ga=N-2$ here is bounded and satisfies
  \bel{SI40}\BA{lll}\dsps
v_{tt}+(2-N)v_t+\Gd'v-me^{-\frac{\gs t}{p-1}}\left((v_t+(2-N)v)^2+|\nabla 'v|^2\right)^{\frac q2}+v^{\frac{N}{N-2}}=0
\EA
\ee
in $(-\infty,0]\ti S^{N-1}$. By (\ref{SI5}), $v(t,.)\to 0$ uniformly when $t\to-\infty$. The average $\bar v$ satisfies 
$$\BA{lll}\dsps\bar v_{tt}+(2-N)\bar v_{t}-
\CH(t)=0,
\EA$$
where 
$$\CH(t)=\frac{1}{|S^{N-1}|}\int_{S^{N-1}}\left(me^{-\frac{\gs t}{p-1}}\left((v_t+(2-N)v)^2+|\nabla 'v|^2\right)^{\frac q2}-v^{\frac{N}{N-2}}\right)dS.
$$
Set $s=e^{(N-2)t}$, $z(s,.)=v(t,.)$ and $\bar z(s)=\bar v(t)$, then there holds
  \bel{SI41+}\BA{lll}\dsps
s^2\bar z_{ss}-Z_1(s)+Z_2(s)=0\quad\text{in }(0,e^{2-N})
\EA
\ee
where
$$Z_1(s)=\frac{ms^{-\frac{\gs}{(p-1)(N-2)}}}{(N-2)^2|S^{N-1}|}\int_{S^{N-1}}\left[(N-2)^2(sz_s-z)^2+|\nabla 'z|^2\right]^{\frac q2}dS
$$
and
$$Z_2(s)=\frac{1}{(N-2)^2|S^{N-1}|}\int_{S^{N-1}}z^{\frac{N}{N-2}}dS.
$$
Using the energy method as in \rlemma{PSing1} and $(\ref{SI5})$ we obtain that
  \bel{SI42}\norm{z(s,.)}_{L^\infty(S^{N-1})}+\norm{sz_s(s,.)}_{L^\infty(S^{N-1})}\to 0\quad\text{as }s\to 0.
\ee
If $0<\gd<1$ the function $s\mapsto w(s):= \bar z(s)+s^\gd$ satisfies
  \bel{SI43}\BA{lll}\dsps
s^2w_{ss}=s^2\bar z_{ss}+\gd(\gd-1)s^\gd=Z_1(s)-Z_2(s)+\gd(\gd-1)s^\gd.
\EA
\ee
We set 
$$\gd_0=\frac{-\gs}{(N-2)(p-1)}=\frac{2p-q(p+1)}{(N-2)(p-1)}=\frac{N-q(N-1)}{N-2},$$
 then $0<\gd_0<1$ since $1<q<\frac{N}{N-1}$. We take $0<\gd<\min\left\{\gd_0,\frac N{N-2}\right\}$.
 Then there exists $s_0>0$ such that  for $0<s\leq s_0$ there holds $Z_1(s)<\frac{\gd(1-\gd)}{2}s^\gd$ which implies
  \bel{SI44}\BA{lll}\dsps
s^2w_{ss}+\frac{\gd(1-\gd)}{2}s^\gd+Z_2(s)\leq 0\quad\text{in }(0,s_0].
\EA
\ee
The function $w$ is therefore concave. Since it vanishes for $s=0$, it is increasing. We now adapt the proof of \cite[Lemma 1]{Avi0} and integrate $(\ref{SI44})$ on $(s,s_0)$. Using the fact that $Z_2(s)\geq \frac{1}{(N-2)^2}\bar z^{\frac{N}{N-2}}(s)$, we obtain
  \bel{SI45}\BA{lll}\dsps 
  w_s(s_0)=w_s(s)+\int_s^{s_0}w_{ss} d\gt\leq w_s(s)-\int_s^{s_0}\left(\frac{\gd(1-\gd)}{2}\gt^{\gd-2}+\frac{Z_2(\gt)}{\gt^{2}}\right)d\gt
 \\[4mm]
 \phantom{w_s(s_0)}\dsps
 \leq w_s(s)-\int_s^{s_0}\left(\frac{\gd(1-\gd)}{2}\gt^{\gd-2}+\frac{\bar z^{\frac{N}{N-2}}(\gt)}{(N-2)^2\gt^{2}}\right)d\gt.\EA
\ee
Since
$$w^\frac{N}{N-2}\leq 2^\frac{2}{N-2}\left(\bar z^\frac{N}{N-2}+s^\frac{N\gd}{N-2}\right),
$$
we infer that
  \bel{SI46}\BA{lll}\dsps w_s(s_0)\leq w_s(s)+\frac{1}{(N-2)(N-2-N\gd)}\left(s^{\frac{N\gd}{N-2}-1}-s_0^{\frac{N\gd}{N-2}-1}\right)\\[4mm]
  \phantom{w_s(s_0)-------------}\dsps-\frac{1}{2^{\frac{2}{N-2}}(N-2)^2}\int_s^{s_0}\frac{w^{\frac{N}{N-2}}(\gt)}{\gt^{2}}d\gt
  \\[5mm]
  \phantom{w_s(s_0)}\dsps
  \leq w_s(s)-C_1\frac{w^{\frac {N}{N-2}}(s)}{s}+C_2s^{\frac{N\gd}{N-2}-1}+C_1\frac{w^{\frac {N}{N-2}}(s)}{s_0}-C_2s_0^{\frac{N\gd}{N-2}-1}
\EA
\ee
for some $C_1,C_2>0$.\\
We claim that
  \bel{SI48}\BA{lll}\dsps 
  w_s(s)-C_1\frac{w^{\frac{N}{N-2}}(s)}{s}+C_2s^{\frac{N\gd}{N-2}-1}\geq 0.
\EA
\ee
Actually, if it were not true there would exist a sequence $\{s_n\}\subset (0,s_0]$ decreasing to $0$ such that 
$$ w_s(s_n)-C_1\frac{w^{\frac{N}{N-2}}(s_n)}{s_n}+C_2s_n^{\frac{N\gd}{N-2}-1}< 0,
$$
which would imply
  \bel{SI49}\dsps w_s(s_0)<C_1\frac{w^{\frac{N}{N-2}}(s_n)}{s_0}-C_2s_0^{\frac{N\gd}{N-2}-1}.
\ee
Since $w(s_n)\to 0$, it would follow that $w_s(s_0)<0$, contradiction.\\
Next we set 
$$\gr(s)=w(s)+cs^{\frac{N\gd}{N-2}},
$$
for some $c>0$ which will be fixed later on. Then, from $(\ref{SI48})$
$$\gr_s(s)\geq C_1\frac{w^{\frac{N}{N-2}}(s)}{s}+\left(c\frac{N\gd}{N-2}-C_2\right)s^{\frac{N\gd}{N-2}-1}.
$$
Now 
$$\gr^\frac{N}{N-2}(s)\leq 2^{\frac{2}{N-2}}\left(w^{\frac{N}{N-2}}(s)+c^{\frac{N}{N-2}}s^{\left(\frac{N}{N-2}\right)^2\gd}\right).
$$
Therefore
$$\BA{lll}
\gr_s(s)\geq C_12^{-\frac{2}{N-2}}\frac{\gr^{\frac{N}{N-2}}(s)}{s}+\left(c\frac{N\gd}{N-2}-C_2\right)s^{\frac{N\gd}{N-2}-1}
-C_12^{-\frac{2}{N-2}}C^{\frac{N}{N-2}}s^{\left(\frac{N}{N-2}\right)^2\gd-1}.
\EA$$
Fixing $c=2C_2\frac{N-2}{N\gd}$, we deduce that for $s$ small enough, 
  \bel{SI51}\dsps 
  \gr_s(s)\geq C_12^{-\frac{2}{N-2}}\frac{\gr^{\frac{N}{N-2}}(s)}{s},
\ee
which implies by integration, 
  \bel{SI51+}\dsps 
  \gr(s)\leq C_3\left(-\ln s\right)^{\frac{2-N}{2}}\quad\text{on }(0,s_1].
\ee

\medskip

\nind {\it 2- End of the proof}. Set $h(t,.)=(-t)^{\frac{N-2}{2}}v(t,.)$, then $h$ is bounded and it satisfies
  \bel{SI57}\BA{lll}\dsps 
h_{tt}+\left(N-2)(1+t)\right)h_t-\frac{1}{t}\left(h^{\frac{2}{N-2}}-\frac{(N-2)^2}{2}\right)h+\frac{N(N-2)}{4t^2}h\\[4mm]
  \phantom{------}\dsps
  -me^{\frac{\gs t}{p-1}}(-t)^{\frac{(2-N)q}{2}}\left(\left(h_t-(N-2)\left(1+\frac 1t\right)h\right)^2+|\nabla'h|^2\right)^{\frac q2}=0.
\EA\ee
Using methods introduced in \cite{Ve2}, it is proved in \cite[Corollary 4.2]{BiRa} that $\norm{h(t,.)-\bar h(t)}_{L^\infty(S^{N-1})}$ tends to $0$ as $t\to\infty$ and consequently that $h(t,.)$ converges in $C^2(S^{N-1})$ to some limit $\ell$ and necessarily
  \bel{SI58}\BA{lll}\dsps 
\ell\in\left\{0,\left(\frac{N-2}{\sqrt 2}\right)^{N-2}\right\}.
\EA\ee
This ends the proof of \rlemma{PSing3} and consequently of \rth {Theorem J}.\qeda

\medskip

\nind {\it Remark 1.} The convergence result {\it 3} of \rth{Theorem G} can be extended to the case $p\in \left(\frac{N}{N-2},\frac{N+1}{N-3}\right)\setminus\{\frac{N+2}{N-2}\}$ for every positive solution $u$ such that 
$|x|^\ga u(x)$ is bounded. 

\medskip

\nind {\it Remark 2.} When $p=\frac{N}{N-2}$, the proof of the existence of solutions of $(\ref{I-1})$ satisfying 
$$\dsps\lim_{x\to 0}|x|^{N-2}\left(-\ln|x|\right)^{\frac{N-2}{2}}=\left(\frac{N-2}{\sqrt 2}\right)^{N-2}
$$
is obtained in the radial case in \cite{BiVe4} using techniques from dynamical systems theory such as the central manifold.

\medskip

\nind {\it Remark 3.} The description of the behaviour in the case $q=\frac{2p}{p+1}$ exhibits a remarkable complexity which appears out of reach in the general case. The treatment of radial solutions is performed in \cite{BiGaVe3} and shows this complexity.

\subsubsection{Proof of \rth{Theorem K}}

Before proving the result we recall that if $q\geq\frac{N}{N-1}$ and $1<p<q$ any nonnegative solution  $u$ of  $(\ref{I-1})$ in $B_{r_0}\setminus\{0\}$ is a bounded weak solution of $(\ref{I-1})$ in $B_{r_0}$ by \rth{Theorem G}.\smallskip

\nind {\it Proof.} Next we assume $p<q<\frac{N}{N-1}$. By \rth {Theorem E} $u$ satisfies
  \bel{SII1}\BA{lll}\dsps 
|x|u(x)+|\nabla u(x)|\leq c_1|x|^{-\frac1{q-1}},
\EA\ee
for $0<|x|\leq {r_0}$. Since $q>\frac{2p}{p+1}$, this implies that $(\ref{SI15})$ holds and therefore $u$ satisfies a uniform Harnack inequality in $B_{\frac {r_0}2}$ in the sense that 
  \bel{SII1*}\BA{lll}\dsps 
u(x)\leq  c_2u(y)\qquad\text{for all }x,y\in B_{\frac {r_0}2}\setminus\{0\}\;\text{ s.t. }|x|=|y|.
\EA\ee
{\it Case 1}. Assume that $|x|^{N-2}u(x)$ is bounded. We cannot apply directly the result of \rth{pre} since $q>\frac{2p}{p+1}$ and we define $u_\ell$ by 
$$u_\ell(x)=\ell^{N-2}u(\ell x)\quad\text{for }\ell>0.
$$
Then $u_k$ satisfies
$$-\Gd u_\ell+m\ell^{N-q(N-1)}|\nabla u_\ell|^q-\ell^{N-p(N-2)}u^p_\ell=0\qquad\text{in }B_{\frac {r_0}\ell}.
$$
Since $q<\frac{N}{N-1}$, $N-q(N-1)> 0$, therefore we deduce as in the proof of \rth{pre} that $\nabla u_\ell$ satisfies estimate $(\ref{Es18})$ with $k$ replaced by $\ell$, which implies 
  \bel{SII1**}\BA{lll}\dsps 
|\nabla u(x)|\leq  c_3|x|^{1-N}\qquad\text{for all }x\in B_{\frac {r_0}2}\setminus\{0\}.
\EA\ee
then 
$$|\nabla u|^{q}\in L^{\frac {N}{N-1}-\ge}(B_{r_0})\quad\text{and }\;u^p\in L^{1}(B_{r_0}),
$$
for any $\ge>0$. By the Brezis-Lions Lemma \cite {BrLi} there exists $k\geq 0$ such that $u$ satisfies 
\bel{SII2}
-\Gd u+m|\nabla u|^q=u^p+k\gd_0\quad\text{in }\CD'(B_{r_0}).
\ee
Furthermore, $u$ verifies
  \bel{SII1***}\BA{lll}\dsps 
\lim_{r\to 0}r^{N-2}u(r,.)=c_Nk
\EA\ee
in $L^1(S^{N-1})$ and actually uniformly. By comparing $u$ with the radial solution $\tilde u_k$ of the Riccatti equation $(\ref{I-5})$ 
  \bel{SII2*}-\Gd u+m|\nabla u|^q=k\gd_0\qquad\text{in }\CD'(B_{{r_0}})
\ee
vanishing on $\prt B_{r_0}$ (see \cite{BiGaVe}), we obtain by the maximum principle that $u\geq \tilde u_k$. The solution $u_k^*$ of $(\ref{SII2*})$ with ${r_0}=\infty$ and vanishing at infinity is explicit and given in \cite[Theorem 3.13]{BiGaVe} by 
  \bel{SII3}\BA{lll}\dsps 
u_k^*(x)=\int_{|x|}^\infty s^{1-N}\left(\frac{q-1}{N-q(N-1)}s^{N-q(N-1)}+c_Nk^{1-q}\right)^{-\frac{1}{q-1}}ds.
\EA\ee
Therefore we easily obtain that the solution $u$ verifies
  \bel{SII3*}\BA{lll}\dsps 
u_k^*(x)-C({r_0})\leq \tilde u_k\leq u(x)\quad\text{for all }x\in B_{r_0}\setminus\{0\},
\EA\ee
for some constant $C({r_0})>0$. \\
If $k=0$, we proceed as in the proof of \rlemma{PSing2}-{\it Step 4} with the same sequences $\{\tilde\gd_n\}$ and $\{\tilde\gth_n\}$. With the notations therein, we obtain $(\ref{SI41-3})$ and $(\ref{SI42-1})$ and derive that $u$ is a bounded regular solution.
\smallskip

\nind {\it Case 2}. Assume that $|x|^{N-2}u(x)$ is unbounded near $x=0$. Then there exists a sequence $\{r_n\}$ decreasing to $0$ such that 
$$\dsps\lim_{r_n\to 0}\sup_{|x|=r_n}r_n^{N-2}u(x)=\infty. 
$$
By $(\ref{SII1*})$ there holds
$$\dsps\lim_{r_n\to 0}\inf_{|x|=r_n}r_n^{N-2}u(x)=\infty. 
$$
Let $k>0$, since $|x|^{N-2}\tilde u_k(x)=c_Nk$, where $\tilde u_k$ has been defined in $(\ref{SII2*})$, for $r_n\leq r_{n_k}$, one has $\tilde u_k\leq u$ in $B_{r_0}\setminus B_{r_n}$
by the maximum principle, which implies that the same inequality holds in  $B_{r_0}\setminus\{0\}$. Let $k\to \infty$ implies that 
$$\dsps\lim_{k\to\infty}\tilde u_k:=\tilde u_\infty\leq u\quad\text{in }B_{r_0}\setminus\{0\}.$$
Since $(\ref{SII3*})$ still holds with $k=\infty$ and combining with \cite[Theorem 3.13]{BiGaVe} we obtain that 
  \bel{SII4}\BA{lll}\dsps 
\xi_m|x|^{-\gb}-C({r_0})\leq \tilde u_\infty\leq u(x)\quad\text{for all }x\in B_{r_0}\setminus\{0\},
\EA\ee
where $\xi_m$ is expressed by $(\ref{S7})$; indeed it is proved in the above mentioned article that $\dsps\lim_{k\to\infty} u^*_k:=u^*_\infty(x)=\xi_m|x|^{-\gb}$. This yields
  \bel{SII5-}\BA{lll}\dsps 
\liminf_{x\to 0}|x|^\gb u(x)\geq\xi_m.
\EA\ee
In order to obtain the sharp estimate from above, we define, for $\ell>0$, $S_\ell[u](x)=\ell^{\gb}u(\ell x)=u_\ell(x)$ in $B_{\frac{{r_0}}{\ell}}\setminus\{0\}$, where $u_\ell$ satisfies
  \bel{SII5}\BA{lll}\dsps 
-\Gd u_\ell+m|\nabla u_\ell|^q=\ell^{\gb(p-1)-2}u_\ell^p.
\EA\ee
Let 
$$\phi^*=\limsup_{|x|\to 0}|x|^\gb u(x)=\lim_{r_n\to 0}r_n^\gb u(r_n,\gth_n),
$$
for some sequence $\{(r_n,\gth_n)\}\to (0,\gth_*)$ and set $u_n(x):=u_{r_n}(x)$. Then $\phi^*\geq \xi_m$ by $(\ref{SII5-})$.
The function $u_n$ satisfies 
  \bel{SII5*}-\Gd u_n+m|\nabla u_n|^q=r_n^{2-\gb(p-1)}u_n^p
\ee
in $B_{\frac{{r_0}}{r_n}}\setminus\{0\}$ and 
  \bel{SII6}|x|u_n(x)+|\nabla u_n(x)|\leq  c_4|x|^{-\frac1{q-1}} \quad\text{if }0<|x|\leq\tfrac{{r_0}}{2r_n}.
\ee
Since $q>p>\frac{2p}{p+1}$, we have $2-\gb(p-1)>0$ and by standard regularity result (see e.g. \cite{GT}), there exists a subsequence, still denoted by $ \{u_{r_n}\}$, and a $C^2$ function $u^*$ such that $u_{r_n}\to u^*$ in the 
$C^2_{loc}$ topology of $\BBR^N\setminus\{0\}$. The function $u^*$ is a nonnegative solution of the Riccatti equation $(\ref{I-5})$ in $\BBR^N\setminus\{0\}$ and it tends to $0$ at $\infty$. By \cite[Theorem 3.13]{BiGaVe}, either $u^*\equiv 0$, either there exists $k>0$ such that $u^*$ verifies $(\ref{SII1***})$,
or 
  \bel{SII6+}\BA{lll}\dsps 
u^*(x)=\xi_m|x|^{-\gb},
\EA\ee
where $\xi_m$ is expressed by $(\ref{S7})$. Note that $\xi_m|x|^{-\gb}$ is the maximal positive solution of $(\ref{I-5})$ in $\BBR^N\setminus\{0\}$ which  tends to $0$ at infinity.  Since $u^*(1,\gs_*)=\phi^*\geq \xi_m$, we obtain that  $\phi^*= \xi_m$ which implies 
  \bel{SII7}\BA{lll}\dsps 
\lim_{x\to 0}|x|^\gb u(x)=\xi_m.
\EA\ee
\qeda
\medskip

\nind \Remark The existence of solutions of $(\ref{SII2})$ for any $k>0$ is proved in the radial case in \cite{BiVe4}. We can observe that if 
$k>0$ is small enough the existence is straightforward since there exists a solution $\hat u_k$ of 
  \bel{SII8}\BA{lll}-\Gd u-u^p=k\gd_0\qquad&\text{in }\CD'(B_{{r_0}})\\
  \phantom{-\Gd -u^p}
  u=0\qquad&\text{in }\prt B_{{r_0}},
\EA\ee
see \cite{Li0}. The function  $\hat u_k$ is a supersolution of $(\ref{I-1})$. Since the solution $\tilde u_k$ of $(\ref{SII2*})$  is a subsolution, and both $\hat u_k$ and $\tilde u_k$ are ordered and have the same behaviour at $0$ given by $(\ref{SII1***})$ it follows that there exists a solution $u_k$ of $(\ref{I-1})$ which vanishes on $\prt B_{r_0}$ and satisfies $\tilde u_k\leq u_k\leq \hat u_k$. Hence it satisfies $(\ref{SII1***})$ and it is easy to check that it is a solution of $(\ref{SII2})$.


\subsection{Behaviour at infinity}

The asymptotic behaviour of positive solutions of $(\ref{I-1})$ in an exterior domain is obtained in some particular cases by using the energy method. Here we make more precise the results contained in \rth{Theorem F}.
\bth{AST1} Let $N\geq 3$, $\frac{N}{N-2}<p<\frac{N-1}{N+3}$, $p\neq \frac{N+2}{N-2}$, $q>\frac{2p}{p+1}$ and $m>0$. If $u$ is a positive solution of $(\ref{I-1})$ in $B_{r_0}^c$ satisfying 
$(\ref{AS2})$ the following alternative holds.\\
(i)  Either 
\bel{Asy1}\dsps
\lim_{|x|\to\infty}|x|^\ga u(x)=\gw_0
\ee
where $\gw_0$ is given by $(\ref{S7-0})$.\\
(ii) Or there exists $k>0$ such that 
\bel{Asy2}\dsps
\lim_{|x|\to\infty}|x|^{N-2} u(x)=k.
\ee
\es
\Proof We recall that estimate $(\ref{AS2})$ holds when  $\frac{N}{N-2}<p< \frac{N+2}{N-2}$ by the doubling method. As in the proof of \rth {Theorem J} we set $u(r,\gth)=r^\ga w(t,\gth)$ with $t=\ln r>0$ (we can assume that ${r_0}<1$) and $w$ is a bounded solution of $(\ref{SI2})$ in $(0,\infty)\ti S^{N-1}$. Notice that $\gs>0$. 
The omega-limit set of the trajectory 
$$\dsps\CT_+[v]=\bigcup _{t\geq0 }v(t,.)
$$
is a non-empty compact connected subset $\Gg_{+}$ of $C^2(S^{N-1})$. The energy method used in the proof of \rth {Theorem J} applies because $p\neq\frac{N+2}{N-2}$, hence
$$\dsps\lim_{t\to\infty}\norm{v_t(t,.)}_{L^2(S^{N-1})}=\lim_{t\to\infty}\norm{v_{tt}(t,.)}_{L^2(S^{N-1})}=0.
$$
This implies that $\Gg_{+}$ is a compact and connected subset of the set of nonnegative solutions of $(\ref{S3})$. Since $\frac{N}{N-2}<p< \frac{N+1}{N-3}$, $\Gg_{+}=\{0,X_0\}$ by 
\cite{GS}, hence if $X_0\in \Gg_{+}$, then $(\ref{Asy1})$ holds, otherwise 
\bel{Asy3}\dsps
\lim_{|x|\to\infty}|x|^\ga u(x)=0.
\ee
In such a case, we obtain by changing $t$ into $-t$ as in the proof of \rlemma{PSing1}, that there exists $\ge>0$ such that 
\bel{Asy4}\dsps
v(t,\gth)\leq  c_1e^{-\ge t}\quad\text{in }(0,\infty)\ti S^{N-1}\Longrightarrow u(x)\leq  c_1|x|^{-\ga-\ge}\quad\text{in }B_{r_0}\setminus\{0\}.
\ee
The computations of \rlemma{PSing2} are still valid, but since $t\to\infty$ the results therein have to be re-interpreted. Since the spherical average $\bar v(t)$ of $v(t,.)$ satisfies $(\ref{SI18*})$, in this equation the right-hand side $\overline H(t)$ which satisfies $\overline H(t)\leq c_2e^{-\gd_1t}$ and $\gd_1$ expressed by $(\ref{gd1})$. By the same standard method of "the variation of constants"  the expression $(\ref{SI19})$ which expressed all the solutions of under the form
\bel{Asy5}\BA{lll}\dsps
\bar v(t)=Ae^{\gl_1t}+Be^{\gl_2t}+C(t)e^{-\gd_1t},
\EA
\ee
where $A$ and $B$ are constant and $C(t)$ is a bounded function. The exponents $\gl_1$ and $\gl_2$ are given by $(\ref{SI18})$. It is important to notice that $\gl_2<0<\gl_1$. Thus, $\bar v(t)\to 0$ when $t\to\infty$ implies $A=0$ and 
 \bel{Asy5-1}\BA{lll}\dsps
\bar v(t)\leq c_3e^{-\gd_1t}\quad\text{for } t>0
\EA
\ee
with $\gd_1$ given by $(\ref{SI31*})$-(i).
The representation formula (\ref{SI25}) valid for $v^*=v-\bar v$ is replaced by
 \bel{Asy6}\BA{lll}\dsps
v^*(t,.)=e^{\frac{2\ga+2-N}{2}t}S(t)v^*(0,.)-\int_0^te^{\frac{2\ga+2-N}{2}s}S(s)\int_0^\infty e^{\frac{N-2\ga-2}{2}\gt}S(\gt)H^*(t+\gt-s,\gs)d\gt ds
\EA
\ee
 see \cite[(1.14)]{BoVe}, where 
 $$\BA{lll}\dsps H^*(t,.)=me^{-\frac{\gs t}{p-1}}\left((v_t-\ga v)^2+|\nabla'v|^2\right)^{\frac {q}{2}}-v^p\\[2mm]
 \phantom{----------}
 \dsps-\myfrac{1}{|S^{N-1}|}\int_{S^{N-1}}\left(me^{-\frac{\gs t}{p-1}}\left((v_t-\ga v)^2+|\nabla'v|^2\right)^{\frac {q}{2}}-v^p\right)dS.
 \EA$$
 Since 
 $$\norm {H(t,.)}_{L^\infty(S^{N-1}}\leq c_4e^{-\gd_1t},
 $$
and $(\ref{SI24})$ holds, we deduce that 
  \bel{Asy7}\BA{lll}\dsps
\norm{v^*(t,.)}_{L^\infty(S^{N-1})}\leq C_1e^{-(N-\ga-1)t}+C_2e^{-\gd_1t}\quad\text{for all }\,t\leq 0.
\EA
\ee
Since $v(t,.)=\bar v(t)+v^*(t,.)$ we deduce
  \bel{Asy8}\BA{lll}\dsps
\norm{v(t,.)}_{L^\infty(S^{N-1})}\leq C_1e^{-(N-\ga-1)t}+C_2e^{-\gd_1t}+C_3e^{-\gth_1t}\leq C_4e^{-\gth_1t}\quad\text{for all }\,t\leq 0,
\EA
\ee
with $\gth_1$ from $(\ref{SI31*})$-(i). We iterate the process and, defining $\gd_n$ and $\gth_n$ by $(\ref{SI31*})$, we obtain, as long as 
$\gth_n<\gl_2$, 
  \bel{Asy8-1}\BA{lll}\dsps
\norm{v(t,.)}_{L^\infty(S^{N-1})}\leq C_1e^{-(N-\ga-1)t}+C_2e^{-\gd_nt}+C_3e^{-\gth_nt}\leq C_4e^{-\gth_nt}\quad\text{for all }\,t\geq 0,
\EA
\ee
Then there exists $n^*$ such that $\gth_{n^*}=\gl_2=\ga+2-N$ and this implies that 
  \bel{Asy8-12}\BA{lll}\dsps
v(t,.)\leq C_5e^{(\ga+2-N)t}.
\EA
\ee
This implies 
$$\bar v(t)=Be^{\gl_2t}(1+o(1)) \quad\text{as }t\to\infty.
$$
Since 
$$\norm{v^*(t,.)}_{L^\infty(S^{N-1})}:=\norm{v(t,.)-\bar v(t)}_{L^\infty(S^{N-1})}\leq C_1e^{-(N-\ga-1)t}+C_2e^{-\gd_{n^*}t}
$$
and $\dsps\gd_{n^*}=\min\left\{p\gth_{n^*},q\gth_{n^*}+\frac{\gs}{p-1}\right\}>\gth_{n^*}$, we conclude that 
  \bel{Asy8-13}\BA{lll}\dsps
\lim_{t\to\infty}e^{(N-2-\ga)t}v(t,.)=B\quad\text{uniformly on } S^{N-1},
\EA
\ee
 which is $(\ref{Asy2})$ with $k=B$. By \rcor{Subcor} we have necessarily $k>0$.
 \qeda\medskip
 
 \nind\Remark The existence of radial solutions in $B_{r_0}^c$ satisfying (\ref{Asy2})  with $k>0$ is proved in \cite{ABGMQ}.\smallskip
 
 The next result completes \rth{Theorem E'}.
 
 \bth{AST2} Let $N\geq 3$, $1<q<\min\{\frac{2p}{p+1},\frac{N}{N-1}\}$ and $m>0$. Let $u$ be a positive solution of $(\ref{I-1})$ in $B_{r_0}^c$.\smallskip
 
 \nind 1- Then
 \bel{Asy9-}\dsps
\liminf_{|x|\to\infty}|x|^{\gb} u(x)\geq \xi_m.
\ee
2- If $|x|^\gb u(x)$ is bounded, then 
\bel{Asy9}\dsps
\lim_{|x|\to\infty}|x|^{\gb} u(x)=\xi_m.
\ee
\es

\nind\Proof For $\ell\geq 1$ the function $u_\ell(x)=\ell^\gb u(\ell x)$ satisfies $(\ref{SII5})$ in $B_{\frac{r_0}{\ell}}^c$ and is bounded therein. Since $q<\frac{2p}{p+1}$, $\gb(p-1)-2<0$, thus we deduce by regularity techniques that  
\bel{Asy9+}\dsps
|x|u(x)+|\nabla u(x)|\leq C|x|^{-\frac{1}{q-1}}.
\ee
This implies that $|x|^2u^{p-1}(x)+|x||\nabla u(x)|^{q-1}\leq C$ in $B_{r_0}^c$, and therefore Harnack inequality holds uniformly in $B_{r_0}^c$ in the sense that 
\bel{Asy10}\dsps
\max_{|x|=r}u(x)\leq C\min_{|x|=r}u(x)\qquad\text{for all }r\geq r_0.
\ee
 Set 
 $\dsps \gm=\min_{|z|=1}u(z)$ and define $k_\gm$ by 
\bel{Asy11}
\gm=u^*_{k_\gm}(1)=\int_1^\infty\left(\frac{q-1}{N-q(N-1)}s^{N-q(N-1)}+k_\gm^{1-q}\right)^{-\frac{1}{q-1}}s^{N-1}ds.
\ee
Then for any $\ge>0$, $u\geq (u^*_k-\ge)_+$ which is a subsolution of the Riccatti equation in $B_1^c$. This implies that $u\geq u^*_{k_\gm}$ in $B_1^c$. Since 
\bel{Asy12}\dsps
\lim_{|x|\to\infty}|x|^\gb u^*_{k_\gm}(x)=\lim_{|x|\to\infty}\int_{|x|}^\infty\left(\frac{q-1}{N-q(N-1)}s^{N-q(N-1)}+k_\gm^{1-q}\right)^{-\frac{1}{q-1}}s^{N-1}ds=\xi_m,
\ee
actually this limit is independent of $k_\gm$, it follows that 
$$\liminf_{|x|\to\infty}|x|^\gb u(x)\geq \xi_m.
$$
This implies $(\ref{Asy9})$.\\
Set 
$$\psi^*=\limsup_{|x|\to\infty}|x|^\gb u(x)=\lim_{r_n\to\infty} r_n^\gb u(r_n,\gth_n)
$$
 where $\gth_n\in S^{N-1}$ and we can assume that $\gth_n\to\gth^*\in S^{N-1}$. Then $\psi^*\geq\xi_m$.
 The function $u_{r_n}:x\mapsto r^\gb_nu(r_nx)$ satisfies 
 \bel{Asy12-n}\dsps
-\Gd u_{r_n}+m|\nabla u_{r_n}|^q=r_n^{2-\gb(p-1)}u^p_{r_n}=r_n^{\frac{\gs}{q-1}}u^p_{r_n}
\ee
in $B_{\frac{r_0}{r_n}}^c$. Since $\gs<0$, we have that $r_n^{\frac{\gs}{q-1}}\to 0$. By  the local regularity a priori estimates inherited from $(\ref{Asy9+})$ implies that, 
up to a subsequence still denoted by $\{r_n\}$, $u_{r_n}$ converge in the $C^2$-local topology of $\BBR^N\setminus\{0\}$ to a positive solution $w$ of 
 \bel{Asy12-w}\dsps
-\Gd w+m|\nabla w|^q=0\qquad\text{in } \BBR^N\setminus\{0\}.
\ee
Because of $(\ref{Asy9+})$
and similarly to the proof of 
 \rth{Theorem K} we can use Arzela-Ascoli theorem to infer that up to a subsequence still denoted by $\{r_n\}$, $u_{r_n}$  converges in the $C^2_{loc}$ topology of $\BBR^N\setminus\{0\}$ 
 to a positive solution of the Riccatti equation $(\ref{I-5})$ in $\BBR^N\setminus\{0\}$ which is a function $u_k^*$ ($0<k\leq\infty$) given by the expression given by $(\ref{SII3})$. Because $\dsps \psi^*=w(1)\geq\xi_m=\lim_{k\to\infty}u_k^*(1)$. Hence $\psi^*=\xi_m$ which conclude the proof.
\qeda \medskip

\mysection {Appendix}
In this Section we prove a technical result concerning the existence of positive radial solutions of 
 \bel{sub2}\dsps
-v''-\frac{N-1}{r}v'+m|v'|^q=0
\ee
on $(r_0,\infty)$ satisfying non-homogeneous Dirichlet conditions at $r=r_0$ and at infinity.
\blemma{Ap1} Let $q>1$, $0<r_0<\gt$ and $a,b>0$. Then there exists a solution $v$ of $(\ref{sub2})$ on $(r_0,\gt)$ satisfying $v(r_0)=a$ and $v(\gt)=b$ if and only if $a=b$, or, if $a\neq b$:\smallskip

\nind 1- When $a< b$, for any $1<q\leq 2$ and $\gt>r_0$.\smallskip

\nind 2- When $a< b$, for any $q> 2$ and $\gt\geq\gt^*>r_0$ where $\gt^*$ depends on $b-a$.

\nind 3- When $a> b$, for any $1<q\leq 2$  and $\gt>r_0$\smallskip

\nind 4- When $a> b$, for any $q> 2$  and $\gt>r_0$ if and only if 
  \bel{sub0}\dsps 
a-b<\left(\frac{q(N-1)-N}{m(q-1)}\right)^{\frac{1}{q-1}}\!\!\!\!r_0^{2-N}\int_1^{\frac{\gt}{r_0}}\!t^{1-N}\left(1-t^{N-q(N-1)}\right)^{-\frac{1}{q-1}}dt.
\ee
\es
\Proof If $a=b$ the constant function $v\equiv a$ is a solution. 
If $v_1$ and $v_2$ are solutions of $(\ref{as4})$ and if there exists $\gth>r_0$ such that $v'_1(\gth)=v'_2(\gth)$, then $v_1=v_2+v_1(\gth)-v_2(\gth)$ by the Cauchy-Lipschitz theorem. This implies in particular that if $v_1$ and $v_2$ are solution either on $(r_0,\gt)$ with $v_1(r_0)=v_2(r_0)$ and $v_1(\gt)=v_2(\gt)$,  or on $(r_0,\infty)$ with $v_1(r_0)=v_2(r_0)$ and $\lim_{r\to\infty}(v_1(r)-v_2(r))=0$, then $v_1=v_2$. We first consider the problem on $(r_0,\gt)$ for some $\gt>r_0$ and if $a, b>0$ we denote by $v:=v_{a,b}$ the solution of 
$(\ref{sub2})$ on $(r_0,\gt)$ such that $v(r_0)=a$ and $v(\gt)=b$. Solutions are explicit by setting $w(r)=r^{N-1}v'(r)$, then 
 \bel{sub2'}w'-mr^{(1-q)((N-1))}|w|^q=0.\ee

\nind {\it Case 1: $a<b$}. If a solution exists it is increasing and we can replace $v$ by $\tilde v=v-a$, thus $\tilde v(r_0)=0$ and $\tilde v'(r)\geq 0$
$$r^{N-1}\tilde v'(r)=\left\{\BA{lll}
\left[(r_0^{N-1}\tilde v'(r_0))^{1-q}-\frac{m(q-1)}{N-q(N-1)}\left(r^{N-q(N-1)}-r_0^{N-q(N-1)}\right)\right]^{-\frac{1}{q-1}}\;&\text{if }\,q\neq\frac{N}{N-1}\\[4mm]
\left[(r_0^{N-1}\tilde v'(r_0))^{1-q}-m(q-1)\ln\frac r{r_0}\right]^{-\frac{1}{q-1}}\;&\text{if }\,q=\frac{N}{N-1}.
\EA\right.$$ 
We set $X:=\tilde v'(r_0)$ and we study the mapping  $r\mapsto \CT_X(r)$ defined by
 \bel{sub3}\dsps
\CT_X(r)=\int_{r_0}^rs^{1-N}\left[(r_0^{N-1}X)^{1-q}-\frac{m(q-1)}{N-q(N-1)}\left(s^{N-q(N-1)}-r_0^{N-q(N-1)}\right)\right]^{-\frac{1}{q-1}}ds
\ee
if $q\neq\frac{N}{N-1}$, and 
 \bel{sub4}\dsps
\CT^*_X(r)=\int_{r_0}^rs^{1-N}\left[(r_0^{N-1}X)^{1-q}-m(q-1)\ln\frac r{r_0}\right]^{-\frac{1}{q-1}}ds
\ee
$q=\frac N{N-1}$. \\
(i) If $N-q(N-1)>0$, $\CT_X$ is defined for $r_0\leq r<r_X:=\left[\frac{N-q(N-1)}{m(q-1)}(r_0^{N-1}X)^{1-q}+r_0^{N-q(N-1)}\right]^\frac{1}{N-q(N-1)}$.
(ii)  If $q=\frac N{N-1}$, $\CT^*(X)$ is defined for $r_0\leq r<r^*_X:=r_0e^{\frac{1}{m(q-1)}(r_0^{N-1}X)^{1-q}}$.\\
 (iii) If $N-q(N-1)<0$, $\CT_X$ is defined for any $r\geq r_0$ if $X\leq X_0:=\left[\frac{N(q-1)-N}{m(q-1)r_0}\right]^\frac{1}{q-1}$, and for 
 $$r<\tilde r_X:=\left[1-\frac{q(N-1)-N}{m(q-1)r_0X^{q-1}}\right]^{-\frac{1}{q(N-1)-N}}r_0
 $$
 if $X>X_0$.
 \smallskip
 
\nind In case (i) (resp. (ii)), we fix $\gt>r_0$ then the mapping $X\mapsto \CT_X(\gt)$ (resp. $X\mapsto \CT^*_X(\gt))$ is continuous, increasing and defined provided $\gt<r_X$ (resp. $\gt<r^*_X$), that is 
 \bel{sub5}\dsps X<X_\gt:=r_0^{1-N}\left[\frac{m(q-1)}{N-q(N-1)}\left(\gt^{N-q(N-1)}-r_0^{N-q(N-1)}\right)\right]^{-\frac{1}{q-1}},\ee
 in case (i) and
 \bel{sub6}\dsps X<X^*_\gt:=r_0^{1-N}\left[m(q-1)\ln\frac \gt{r_0}\right]^{-\frac{1}{q-1}}
\ee
in case (ii). Furthermore $\CT_0(\gt)=\CT^*_0(\gt)=0$ and $\lim_{X\uparrow X_\gt}\CT_X(\gt)=\lim_{X\uparrow X^*_\gt}\CT^*_X(\gt)=\infty$ since $q\leq 2$. As a consequence there exists a unique 
$\tilde X\in (0,X_\gt)$ (resp. $\tilde X\in (0,X^*_\gt)$) such that $\CT_{\tilde X}(\gt))=b-a$ (resp. $\CT^*_{\tilde X}(\gt))=b-a$).\smallskip

\nind In case (iii) we have in the case $X\leq X_0$, 
 \bel{sub7}\dsps\lim_{r\to\infty}\CT_X(r)=\left\{\BA{lll}\infty \quad&\text{if }N=2\\[2mm]
C_1(X):=\frac{r_0X}{N-2}\left[1-\left(\frac{X}{X_0}\right)^{q-1}\right]^{-\frac{1}{q-1}}\quad&\text{if }N\geq 3.
\EA\right.
\ee
Since $C_1(0)=0$ and $C_1(X)\to\infty$ when $X\uparrow X_0$, $C_1$ is a continuous increasing function from $[0,X_0]$ onto $[0,\infty]$. If $X> X_0$, 
 \bel{sub8}\dsps\lim_{r\to \tilde r_X}\CT_X(r)=\left\{\BA{lll}\infty \qquad&\text{if }\,\frac{N}{N-1}<q\leq 2\\[2mm]
C_2(X)\qquad&\text{if }\,q>2,
\EA\right.
\ee
where 
 \bel{sub9}\dsps C_2(X)=\left(\frac{q(N-1)-N}{m(q-1)}\right)^{\frac{1}{q-1}}\tilde r_X^{\frac{q-2}{q-1}}\int_{\frac{r_0}{\tilde r_X}}^1\left(t^{N-q(N-1)}-1\right)^{-\frac{1}{q-1}}t^{1-N}dt.\ee
For $\gt>r_0$, we introduce again the mapping $X\mapsto \CT_X(\gt)$.  In view of the last relation in the case $\frac{N}{N-1}<q\leq 2$ then for any $b>a$ and $\gt>r_0$ there exists a unique $\tilde X>X_0$ such that $\gt<r_{\tilde X}$ and $ \CT_{\tilde X}(\gt)=b-a$.\\
 If $q>2$ and $N\geq 3$, for any $b>a$ there exists $\gt^*>r_0$, depending on $b-a$,  such that for any 
$\gt\geq \gt^*$ there exists $X\leq X_0$ such that $\CT_X(\gt)=b-a$. We can explicit $\gt^*$  by  $\gt^*=\tilde r_{X^*}$ where $X^*$ is characterized by 
$C_2(X^*)=b-a$.
\smallskip

\nind {\it Case 2}: $a>b$. Then $v$ is decreasing
and the method has to be slightly modified in order to obtain a positive solution of $-v''-\frac{N-1}{r}v'+m|v'|^q=0$ on $(r_0,\gt)$ such that $v(r_0)=a$ and $v(\gt)=b$. By replacing $v$ by $\tilde v:=v-b$ we look for a solution $\tilde v$ vanishing at $\gt$ and positive on $(r_0,\gt)$. Let $X=\tilde v'(r_0)$ then 
$$
-r^{N-1}\tilde v'(r)=\left\{\BA{lll}
\left[ (-r_0^{N-1}X)^{1-q}+\frac{m(q-1)}{N-q(N-1)}\left(r^{N-q(N-1)}-r_0^{N-q(N-1)}\right)\right]^{-\frac{1}{q-1}}\quad &\text{if }q\neq\frac{N}{N-1}\\[2mm]
\left[(-r_0^{N-1}X)^{1-q}+m(q-1)\ln\frac r{r_0}\right]^{-\frac{1}{q-1}}\quad &\text{if }q=\frac{N}{N-1}.
\EA\right.
$$
We study the mapping $r\mapsto\CS_X(r)$ defined by
 \bel{sub10}\dsps 
 \CS_X(r)=a-b-\int_{r_0}^rs^{1-N}\left[ (-r_0^{N-1}X)^{1-q}+\frac{m(q-1)}{N-q(N-1)}\left(s^{N-q(N-1)}-r_0^{N-q(N-1)}\right)\right]^{-\frac{1}{q-1}}ds
 \ee
if $q\neq\frac{N}{N-1}$ and 
 \bel{sub11}\dsps 
 \CS^*_X(r)=a-b-\int_{r_0}^rs^{1-N}\left[(-r_0^{N-1}X)^{1-q}+m(q-1)\ln\frac s{r_0}\right]^{-\frac{1}{q-1}}ds
 \ee
if $q=\frac{N}{N-1}$.
If $q\leq \frac N{N-1}$, these two functions are defined on $(r_0,\gt)$. A solution $\tilde v$ satisfying the boundary conditions at $r=r_0$ and $r=\gt$ corresponds to the fact that  $\CS_X(\gt)=0$ if  $q\neq\frac{N}{N-1}$ or $\CS^*_X(\gt)=0$ if  $q=\frac{N}{N-1}$. \\
(i) If $q<\frac{N}{N-1}$ we have 
 \bel{sub12}\dsps 
\lim_{X\uparrow 0} \CS_X(\gt)=a-b\,\text{ and }\,\lim_{X\to -\infty} \CS_X(\gt)=-\infty,
 \ee
 because $q<2$ implies that $\myint{r_0}{\gt} s^{1-N}\left[ \frac{m(q-1)}{N-q(N-1)}\left(s^{N-q(N-1)}-r_0^{N-q(N-1)}\right)\right]^{-\frac{1}{q-1}}ds=\infty.$\\
 (ii) If $q=\frac{N}{N-1}$ we have also
  \bel{sub13}\dsps 
\lim_{X\uparrow 0} \CS^*_X(\gt)=a-b\,\text{ and }\,\lim_{X\to -\infty} \CS^*_X(\gt)=-\infty.
 \ee
 This implies that in these two cases for any $\gt>0$ there exists a unique $X<0$ such that $\CS_X(\gt)=0$ or $\CS^*_X(\gt)=0$.\\
\nind (iii) If $q>\frac{N}{N-1}$, $\CS_X(r)$ is defined for any $X\leq 0$ and any $r\in (r_0,\gt)$. We write it under the form
 \bel{sub14}\dsps 
 \CS_X(\gt)=a-b-\int_{r_0}^\gt s^{1-N}\left[ (-r_0^{N-1}X)^{1-q}+\frac{m(q-1)}{q(N-1)-N}\left(r_0^{N-q(N-1)}-s^{N-q(N-1)}\right)\right]^{-\frac{1}{q-1}}ds
 \ee
 We have that $\dsps\lim_{X\uparrow 0} \CS_X(\gt)=a-b$ and $\dsps\lim_{X\to -\infty} \CS_X(\gt)=-\infty$ if $\frac{N}{N-1}<q\leq 2$; in such  case there exists $X_\gt<0$ such that  $\CS_{X_\gt}(\gt)=0$. On the contrary, if $q> 2$, we have
  \bel{sub15}\dsps 
\lim_{X\to -\infty} \CS_X(\gt)=
a-b-\left(\frac{q(N-1)-N}{m(q-1)}\right)^{\frac{1}{q-1}}\!\!\!\!r_0^{2-N}\int_1^{\frac{\gt}{r_0}}\!t^{1-N}\left(1-t^{N-q(N-1)}\right)^{-\frac{1}{q-1}}dt.
\ee
In that case  we can find some $X=X(\gt)<0$ (actually always unique) such that $\CS_{X_\gt}(\gt)=0$ if and only if 
  \bel{sub16}\dsps 
a-b<\left(\frac{q(N-1)-N}{m(q-1)}\right)^{\frac{1}{q-1}}\!\!\!\!r_0^{2-N}\int_1^{\frac{\gt}{r_0}}\!t^{1-N}\left(1-t^{N-q(N-1)}\right)^{-\frac{1}{q-1}}dt.
\ee
Letting $\gt\to\infty$ we can find $\gt>r_0$ such that (\ref{sub16}) holds if and only if $(\ref{sub0})$ holds.
\qeda
\medskip

\end{document}